\documentstyle[11pt]{article}
\title{On some aspects of the Deligne-Simpson problem
\footnote{Research partially supported by INTAS grant 97-1644}}
\author{Vladimir Petrov Kostov\\ \\ \hspace{7cm}{\sl To prof. V.I.Arnold}} 
\date{}
\newtheorem{tm}{Theorem}
\newtheorem{lm}[tm]{Lemma}
\newtheorem{cor}[tm]{Corollary}
\newtheorem{prop}[tm]{Proposition}
\begin{document}
\maketitle 

\begin{abstract}
The Deligne-Simpson problem in the multiplicative version 
is formulated like this: {\em give necessary and 
sufficient conditions for the choice of the conjugacy classes 
$C_j\in SL(n,{\bf C})$ so that there exist 
irreducible $(p+1)$-tuples of matrices $M_j\in C_j$  
satisfying the equality $M_1\ldots M_{p+1}=I$}. 

We solve the problem for generic eigenvalues in the 
case when all the numbers $\Sigma _{j,m}(\sigma )$ of 
Jordan blocks of a given matrix $M_j$, with a given 
eigenvalue $\sigma$ and of a given size $m$ (taken over all $j$, $\sigma$, 
$m$) are divisible by $d>1$. 
Generic eigenvalues are defined by explicit algebraic inequalities of the 
form $a\neq 0$. For such eigenvalues there exist no reducible 
$(p+1)$-tuples.  

The matrices $M_j$ are interpreted as monodromy operators of 
regular linear systems on Riemann's 
sphere. \\  

{\bf Key words:} generic eigenvalues, (poly)multiplicity vector, 
corresponding Jordan normal forms, monodromy operator.

{\bf AMS classification index:} 15A30, 20G05
\end{abstract}
\tableofcontents 

\section{Introduction}
\subsection{Formulation of the problem\protect\label{formulation}}

In the present paper we consider the multiplicative version of the 
{\em Deligne-Simpson problem:}\\ 

{\em Give necessary and sufficient conditions for the choice of the 
conjugacy classes $C_j\in SL(n,{\bf C})$ so that there exist irreducible 
$(p+1)$-tuples of matrices $M_j\in C_j$ such that} 

\begin{equation}\label{M_j}
M_1\ldots M_{p+1}=I 
\end{equation}
In the additive version the conjugacy classes $c_j$ belong to 
$sl(n,{\bf C})$ and the matrices $A_j\in c_j$ satisfy the condition
\begin{equation}\label{A_j}
A_1+\ldots +A_{p+1}=0 
\end{equation}
The matrices $A_j$ and $M_j$ are interpreted respectively as 
matrices-residua and as monodromy operators of fuchsian linear systems, see 
the next section. Both versions of the problem were considered in \cite{Ko1} which is 
the first part of this paper. 

The problem (in the multiplicative version) was stated by P.Deligne and 
C.Simpson obtained the first results towards its solution, see \cite{Si}. 

Further in the text we consider sometimes $C_j$ ($c_j$) as conjugacy classes 
from $GL(n,{\bf C})$ (from $gl(n,{\bf C})$) instead of $SL(n,{\bf C})$ 
(instead of $sl(n,{\bf C})$) because when solving the problem 
there appear such matrices. The passage from  
the problem for $M_j\in SL(n,{\bf C})$ (or $A_j\in sl(n,{\bf C})$) to the 
one for $M_j \in GL(n,{\bf C})$ (or $A_j\in gl(n,{\bf C})$) and vice versa 
is trivial.

We presume that there hold the necessary conditions 
$\sum _{j=1}^{p+1}$Tr$(c_j)=0$ and $\prod _{j=1}^{p+1}\det (C_j)=1$. 
This means that the eigenvalues $\lambda _{k,j}$ ($\sigma _{k,j}$) of the 
matrices $A_j$ ($M_j$) satisfy the conditions 

\[ \sum _{j=1}^{p+1}\sum _{k=1}^n\lambda _{k,j}=0~~,~~
\prod _{j=1}^{p+1}\prod _{k=1}^n\sigma _{k,j}=1\] 
Here $k=1,\ldots ,n$ and for $j$ fixed the eigenvalues are not presumed 
distinct. By definition, the {\em multiplicity} of the eigenvalue 
$\lambda _{k,j}$ (or $\sigma _{k,j}$) is the number of eigenvalues 
$\lambda _{i,j}$ (or $\sigma _{i,j}$; $j$ is fixed) equal to it including 
the eigenvalue itself. 

Define as {\em generic} any set of eigenvalues $\lambda _{k,j}$  
or $\sigma _{k,j}$ which satisfy none of the equalities 

\[ \begin{array}{rllcrlll}
\sum _{j=1}^{p+1}\sum _{k\in \Phi _j}\lambda _{k,j}&=&0&,&
\prod _{j=1}^{p+1}\prod _{k\in \Phi _j}\sigma _{k,j}&=&1~~~~~~&(\gamma )
\end{array}\]
where the sets $\Phi _j$ contain one and the same number $\kappa$ of indices 
($1<\kappa <n$) for all $j$. Reducible $(p+1)$-tuples exist only for 
non-generic eigenvalues -- if a $(p+1)$-tuple is block upper-triangular, 
then the eigenvalues of each diagonal block satisfy some relation ($\gamma$). 

{\bf Definition.} Call {\em Problem (I)} the Deligne-Simpson problem like it 
is formulated above and {\em Problem (TC)} the same problem in which the 
requirement of irreducibility is replaced by the requirement the centralizer 
of the $(p+1)$-tuple to be trivial. 

It is clear that Problem (TC) is weaker than Problem (I) and that for 
generic eigenvalues the answers to both problems coincide. Part of the 
results from this paper concern Problem (TC).

\subsection{The quantities $q$, $d$, $\xi$ and $m_0$\protect\label{qdxm}}

{\bf Definition.} A {\em multiplicity vector (MV)} is a vector whose 
components are non-negative integers whose sum equals $n$. We always 
interpret a MV as the vector of the multiplicities of the eigenvalues of 
a matrix $M_j$ or $A_j$. A {\em polymultiplicity vector (PMV)} is a 
$(p+1)$-tuple or $(p+2)$-tuple of MVs (depending on whether we deal with a 
$(p+1)$- or $(p+2)$-tuple of matrices $A_j$ or $M_j$). A PMV is called 
{\em simple} (resp. {\em non-simple}) if the greatest common divisor $q$ of 
all the components of all its MVs equals 1 (if not). 

{\bf Notation.} Denote by 
$\Sigma _{j,m}(\sigma )$ the number of Jordan blocks of size $m$, of a given 
matrix $M_j$ and corresponding to its eigenvalue $\sigma$. Denote 
by $d$ the greatest common divisor of all numbers $\Sigma _{j,m}(\sigma )$ 
(over all $\sigma$, $j$ and $m$).

{\bf Remark:} The condition $q=1$ implies $d=1$, and $d>1$ implies $q>1$ but 
the inverse implications are not true. It is true that $d$ divides $q$ and 
that $q$ divides $n$. 

In \cite{Ko1} we considered the case of generic eigenvalues: the 
additive version was considered completely, in the multiplicative one we 
considered only the situation when $d=1$. In the present second part of 
\cite{Ko1} we consider the 
multiplicative version for generic eigenvalues and $d>1$ (we call it 
{\em the case~}$d>1$).

{\bf Notation.} Let $q>1$. Denote by $\xi$ the 
product of the eigenvalues $\sigma _{k,j}$ when their multiplicities are 
reduced $q$ times. Hence, $\xi$ is a root of unity of order $q$. Set 
$\xi =\exp (2\pi im_0/q)$. If $(m_0,q)=1$, then $\xi$ is a primitive root. 
Another equivalent definition of $m_0$ is given in 5) of the remarks in 
Subsection~\ref{Leveltsresult}. 

{\bf Example:} If $p=n=2$ and if each of the three Jordan normal forms 
consists of one Jordan block of size 2, then there are two possibilities -- 
either $\sigma _{1,1}\sigma _{1,2}\sigma _{1,3}=1$ or 
$\sigma _{1,1}\sigma _{1,2}\sigma _{1,3}=-1$ (because 
$\sigma _{1,j}=\sigma _{2,j}$ and 
$(\sigma _{1,1}\sigma _{1,2}\sigma _{1,3})^2=1$). In the first (in the 
second) case there are no (there are) irreducible triples of matrices $M_j$ 
satisfying (\ref{M_j}). In the first (in the second) case the eigenvalues 
are not (are) generic. In both cases one has $q=2$, $d=1$. In the first 
case $\xi =1$, in the second case $\xi =-1$. 

{\bf Remark:} If $q>1$ and if $\xi$ is not a primitive root of unity of 
order $q$, then the eigenvalues $\sigma _{k,j}$ are not generic. 
However, $\xi$ can be a primitive root of unity and the 
eigenvalues can be non-generic. 

{\bf Example:} Let $n=4$, $p=3$ and let all Jordan normal forms be diagonal. 
Let for all $j$ one have $\sigma _{1,j}=\sigma _{2,j}$
$\neq \sigma _{3,j}=\sigma _{4,j}$. Hence, $d=q=2$. Let 
$\sigma _{1,1}=\sigma _{1,2}=\sigma _{1,3}=1$,  
$\sigma _{1,4}=i$, $\sigma _{3,1}=\sigma _{3,2}=\exp (\pi i/4)$, 
$\sigma _{3,3}=\exp (-\pi i/8)$, 
$\sigma _{3,4}=\exp (\pi i/8)$. Hence,  
$\xi =-1$. One has $\sigma _{1,1}\sigma _{1,2}\sigma _{3,3}\sigma _{3,4}=1$ 
which is a non-genericity relation.  

In the additive version, when $q>1$, the eigenvalues $\lambda _{k,j}$ 
satisfy a non-genericity relation -- if their multiplicities are reduced 
$q$ times, then their sum is 0. Denote this relation by $(\gamma _1)$.

In the multiplicative version, when $\xi$ is not a primitive root of unity, 
the eigenvalues $\sigma _{k,j}$ satisfy a non-genericity relation 
$(\gamma ^*)$ -- when their multiplicities are reduced $(m_0,q)$ times, then 
their product equals 1. The {\em multiples} of $(\gamma ^*)$ are obtained 
(by definition) 
when after the reduction $(m_0,q)$ times the multiplicities are increased 
$s$ times with $1\leq s<(m_0,q)$. The multiples are denoted by 
$s(\gamma ^*)$. (The multiples of $(\gamma _1)$ are defined by analogy.)

{\bf Definition.} If $q>1$, if in the multiplicative version $(m_0,q)>1$ 
and if the only non-genericity relation satisfied by the eigenvalues 
$\lambda _{k,j}$ (resp. $\sigma _{k,j}$) is $(\gamma _1)$ (resp. 
$(\gamma ^*)$) and its multiples, then the eigenvalues are called 
{\em relatively generic}. 

\subsection{The results obtained up to now}

{\bf Definition.} Call {\em Jordan normal form of size $n$} a family 
$J^n=\{ b_{i,l}\}$ ($i\in I_l$, $I_l=\{ 1,\ldots ,s_l\}$, $l\in L$) of 
positive integers $b_{i,l}$ 
whose sum is $n$. Here $L$ is the set of eigenvalues (all distinct) and 
$I_l$ is the set of Jordan blocks with eigenvalue $l$, $b_{i,l}$ is the 
size of the $i$-th block with this eigenvalue. We presume that the 
following inequalities hold (for each $l$ fixed):

\[ b_{1,l}\geq b_{2,l}\geq \ldots \geq b_{s_l,l}\] 
An $n\times n$-matrix 
$Y$ has the Jordan normal form $J^n$ (notation: $J(Y)=J^n$) if to its 
distinct eigenvalues $\lambda _l$, $l\in L$, there belong Jordan blocks of 
sizes $b_{i,l}$.

For a conjugacy class $C$ in $GL(n,{\bf C})$ or $gl(n,{\bf C})$ denote by 
$d(C)$ its dimension and for a matrix $Y$ from $C$ set 
$r(C):=\min _{\lambda \in {\bf C}}{\rm rk}(Y-\lambda I)$. The integer 
$n-r(C)$ is the maximal number of Jordan blocks of $J(Y)$ with one and the 
same eigenvalue. Set $d_j:=d(C_j)$ (resp. $d(c_j)$), $r_j:=r(C_j)$ 
(resp. $r(c_j)$). The quantities 
$r(C)$ and $d(C)$ depend only on the Jordan normal form $J(Y)=J^n$, not 
on the eigenvalues, so we write sometimes $r(J^n)$ and $d(J^n)$. 

{\em The following two conditions are necessary for the 
existence of irreducible $(p+1)$-tuples of matrices $M_j$ satisfying 
(\ref{M_j}) (proved in \cite{Si}):}

\[ \begin{array}{rlllll}d_1+\ldots +d_{p+1}&\geq &2n^2-2&&&~~~~(\alpha _n)\\
{\rm for~all~}j~(r_1+\ldots +\hat{r}_j+\ldots +r_{p+1})&\geq &n&&&~~~~
(\beta _n)\end{array}\]

These conditions are also necessary for the existence of irreducible 
$(p+1)$-tuples of matrices $A_j$ satisfying (\ref{A_j}), see \cite{Ko1}. 

Recall the basic result from \cite{Ko1}. Its formulation depends not on 
the conjugacy classes $C_j$ but only on the Jordan normal forms defined by 
them (as long as the eigenvalues remain generic).

For a given $(p+1)$-tuple $J^n$ of Jordan normal forms 
$J_j^n$ (the upper index indicates the size of the 
matrices) define the $(p+1)$-tuples of Jordan normal forms 
$J^{n_{\nu}}$, $\nu =0$, $\ldots$, $s$ as follows. 
Set $n_0=n$. If $J^n$ satisfies the condition 

\[ \begin{array}{rlllll}r_1+\ldots +r_{p+1}&\geq &2n&&&~~~(\omega _n)
\end{array}\] 
or if it doesn't satisfy condition $(\beta _n)$ or if 
$n=1$, then set $s=0$. If not, then set 

\[ n_1=r_1+\ldots +r_{p+1}-n~~{\rm (hence,}~n_1<n{\rm).}\] 
For $j=1,\ldots ,p+1$ the Jordan normal form $J_j^{n_1}$ is obtained from 
$J_j^n$ {\em by finding (one of) the eigenvalue(s) of $J_j^n$ with greatest 
number of Jordan blocks} and by decreasing 
by 1 the sizes of the $n-n_1$ {\em smallest} Jordan blocks with this 
eigenvalue. (Their number is $n-r_j$ which is $\geq n-n_1$ because if $J^n$ 
satisfies condition $(\beta _n)$, then $n_1\geq r_j$.) Denote symbolically 
the construction of the Jordan normal forms $J_j^{n_1}$ by 

\[Ê\Psi :\left\{ \begin{array}{ccc}n&\mapsto &n_1\\ 
(J_1^n,\ldots ,J_{p+1}^n)&\mapsto &(J_1^{n_1},\ldots ,J_{p+1}^{n_1})
\end{array}\right.\] 

Let the $(p+1)$-tuples of Jordan normal forms $J^{n_0}$, $\ldots$, 
$J^{n_{\nu _0}}$ be constructed. 
If $J^{n_{\nu _0}}$ satisfies condition $(\omega _{n_{\nu _0}})$ or if 
it doesn't satisfy condition $(\beta _{n_{\nu _0}})$ or if $n_{\nu _0}=1$, 
then set $s=n_{\nu _0}$. If not, then define $J^{n_{\nu _0+1}}$ after 
$J^{n_{\nu _0}}$ in the same way as $J^{n_1}$ was defined after $J^n$, see 
$\Psi$. One has $n=n_0>n_1>\ldots >n_s$. 

In the additive version and in the case $d=1$ of the multiplicative one 
the following theorem is true (see \cite{Ko1}):

\begin{tm}\label{generic1}
For given conjugacy classes $c_j$ or $C_j$ with Jordan normal forms $J_j^n$ 
with generic 
eigenvalues there exist irreducible $(p+1)$-tuples of matrices $A_j$ or 
$M_j$ with Jordan normal forms $J_j^n$ satisfying respectively (\ref{A_j}) 
or (\ref{M_j}) if and only if the following two conditions hold:

{\em i)} The $(p+1)$-tuple of Jordan normal forms $J_j^n$ satisfies the 
inequalities $(\alpha _n)$ and $(\beta _n)$;

{\em ii)} Either the $(p+1)$-tuple of Jordan normal forms $J_j^{n_s}$ 
satisfies the inequality $(\omega _{n_s})$, or one has $n_s=1$.
\end{tm}

{\bf Remarks:} 1) Set $d_j^{\nu}=d(J_j^{n_{\nu}})$, $d_j^0=d_j=d(J_j^n)$. 
Set $d_1+\ldots +d_{p+1}=$$2n^2-2+\kappa$, $\kappa \in {\bf N}$.  
The number $n$ is divisible by $d$ and the numbers $d_j$ are divisible by 
$d^2$. Hence, if $d>1$ and if 
condition $(\alpha _n)$ holds, then one must have $\kappa \geq 2$. 

The quantity $2-\kappa$ is called {\em index of rigidity} (see \cite{Ka}). 
Irreducible representations with $\kappa =0$ are called {\em rigid}; they 
are unique up to conjugacy, see \cite{Ka} and \cite{Si}. In the last remark 
of this section we explain how to give the exhaustive 
list of $(p+1)$-tuples of Jordan normal forms defining irreducible 
representations of a given index of rigidity and with generic eigenvalues.

2) On the other 
hand, one has (for all $\nu$) $d_1^{\nu}+\ldots +d_{p+1}^{\nu}=$
$2(n_{\nu})^2-2+\kappa$. Indeed, for all $\nu$ one has 
$\sum _{j=1}^{p+1}r(J_j^{n_{\nu}})=n_{\nu}+n_{\nu +1}$ by definition. 
One has as well 
$d_j^{\nu +1}=d_j^{\nu}-2(n_{\nu}-n_{\nu +1})r(J_j^{n_{\nu}})$ (to be 
checked directly). Hence,   

\[ \sum _{j=1}^{p+1}d_j^{\nu +1}=\sum _{j=1}^{p+1}d_j^{\nu}-
2(n_{\nu}-n_{\nu +1})\sum _{j=1}^{p+1}r(J_j^{n_{\nu}})=\] 
\[ =2(n_{\nu})^2-2+\kappa 
-2(n_{\nu}-n_{\nu +1})(n_{\nu}+n_{\nu +1})=2(n_{\nu +1})^2-2+\kappa\] 
Thus in the case 
$d>1$ one can never have $n_s=1$ because for $n_s=1$ one has $\kappa =0$.  

3) It is shown in \cite{Ko2} that if the $(p+1)$-tuple of Jordan normal 
forms $J_j^n$ satisfies condition $(\omega _n)$, then one has 
$d_1+\ldots +d_{p+2}\geq 2n^2$. Hence, if the Jordan normal forms 
$J_j^{n_\nu}$ satisfy condition $(\omega _{n_\nu})$, then 
$d_1^{\nu}+\ldots +d_{p+1}^{\nu}\geq 2(n_{\nu})^2$. 
This together with 1) and 2) implies that 
$n_s=1$ if and only if $(\alpha _n)$ is an equality. 

{\bf Definition.} Call a $(p+1)$-tuple $J^n$ of Jordan normal forms $J_j^n$ 
{\em good} if it satisfies conditions {\em i)} and {\em ii)} of the theorem.

\begin{cor}\label{good}
Let $n>1$. Then the $(p+1)$-tuple $J^n$ of Jordan normal forms $J_j^n$ is 
good when either it satisfies condition $(\omega _n)$ or when the 
$(p+1)$-tuple 
of Jordan normal forms $J_j^{n_1}$ is good (and only in these two cases). 
\end{cor}
 
This follows from the definition of the Jordan normal forms $J_j^{n_{\nu}}$. 
Indeed, if $J^n$ is good, then it either satisfies condition $(\omega _n)$ 
or the $(p+1)$-tuple $J^{n_s}$ satisfies condition $(\omega _{n_s})$ (recall 
that if $d>1$, the possibility $n_s=1$ is excluded). In the second case 
one proves by induction on $k=s-\nu$ that all $(p+1)$-tuples $J^{n_{\nu}}$ 
are good.   

Theorem \ref{generic1} can be reformulated as follows: 

{\em In the additive version and in the case $d=1$ of the multiplicative one 
there exist for generic eigenvalues irreducible $(p+1)$-tuples of 
matrices $A_j$ or $M_j$ satisfying respectively condition (\ref{A_j}) or 
(\ref{M_j}) if and only if the $(p+1)$-tuple of Jordan normal forms $J_j^n$ 
is good.}

\subsection{The basic results of this paper and what still remains to be 
done}

The aim of the present paper is to prove

\begin{tm}\label{generic2}
Conditions {\em i)} and {\em ii)} from Theorem \ref{generic1} are necessary 
and sufficient for the existence of irreducible $(p+1)$-tuples of matrices 
$M_j$ satisfying (\ref{M_j}) and with generic eigenvalues in the case $d>1$.
\end{tm}

The sufficiency of conditions {\em i)} and {\em ii)} from 
Theorem \ref{generic1} follows from 

\begin{tm}\label{suff}
Let $d>1$ and let $\xi$ be a primitive root of unity. Then conditions 
{\em i)} and {\em ii)} from Theorem \ref{generic1} are sufficient 
for the existence of $(p+1)$-tuples of matrices $M_j$ with trivial 
centralizers.
\end{tm}

Theorem \ref{suff} is stronger than the proof 
of the sufficiency in Theorem~\ref{generic2} because the eigenvalues are 
not presumed generic.

\begin{tm}\label{suff1}
If $d>1$ 
and if $d_1+\ldots +d_{p+1}\geq 2n^2+2$, then conditions {\em i)} and 
{\em ii)} from Theorem~\ref{generic1} are necessary and sufficient for the 
existence of $(p+1)$-tuples of matrices $M_j$ satisfying 
(\ref{M_j}) and with trivial centralizers. If the eigenvalues are relatively 
generic, then there exist irreducible such $(p+1)$-tuples. 
\end{tm}

Notice that the theorem does not require primitivity of $\xi$. The last two 
theorems are proved in Section~\ref{sufficiency}.

\begin{lm}\label{2n2}
A) If $q>1$, if $d_1+\ldots +d_{p+1}=2n^2$, if all $(p+1)$ Jordan normal 
forms are diagonal, if the eigenvalues 
$\sigma _{k,j}$ are relatively generic and if $\xi$ is not primitive, then 
such a $(p+1)$-tuple of matrices $M_j$ satisfying (\ref{M_j}) (if it 
exists) is with 
trivial centralizer if and only if it is irreducible.

B) Let $d_1+\ldots +d_{p+1}\geq 2n^2+2$ and let for the rest 
the conditions from A) hold. Then if there exist $(p+1)$-tuples of matrices 
$M_j$ satisfying (\ref{M_j}), then there exist 
irreducible ones as well.
\end{lm}

The lemma is proved in Subsection~\ref{proofof2n2}.

\begin{tm}\label{necessary}
Conditions {\em i)} and 
{\em ii)} from Theorem~\ref{generic1} are necessary for the 
existence of $(p+1)$-tuples of matrices $M_j$ or $A_j$ 
satisfying (\ref{M_j}) or (\ref{A_j}) and with a trivial centralizer.
\end{tm}

The theorem claims necessity in all possible situations. It is proved 
in Section~\ref{necessity}. 

For generic eigenvalues Problem (I) (and Problem (TC) as well) is completely 
solved by Theorems \ref{generic1} and \ref{generic2}. 

For non-generic eigenvalues we focus on Problem (TC). The situations in 
which the answer is known are given by Theorems~\ref{suff} and \ref{suff1}. 
It is shown in \cite{Ko1} that conditions {\em i)} and 
{\em ii)} from Theorem~\ref{generic1} are necessary and sufficient for the 
solvability of Problem (TC) when $d=1$, $(\alpha _n)$ being a strict 
inequality.

The cases in which Problem (TC) remains to be considered for 
non-generic eigenvalues (in both versions -- additive or multiplicative) 
and the conjectures the author makes are:

1) $(\alpha _n)$ is an equality (this implies $d=1$, see the remarks after 
Theorem~\ref{generic1}).

{\bf Conjecture.} {\em Conditions {\em i)} and 
{\em ii)} from Theorem~\ref{generic1} are necessary and sufficient if $q=1$. 
For $q>1$ there are cases in which they are and cases in which they are not 
sufficient.}

2) $d>1$, $d_1+\ldots +d_{p+1}=2n^2$ and in the multiplicative version 
$\xi$ is not primitive.

{\bf Conjecture.} {\em There exist no $(p+1)$-tuples of matrices $A_j$ or 
$M_j$, satisfying (\ref{A_j}) or (\ref{M_j}), with trivial centralizers.} 

{\bf Remark:} It is possible to give an exhaustive list of the 
$(p+1)$-tuples of Jordan normal forms of a given size $n$ admitting generic 
eigenvalues and 
defining irreducible representations of a fixed index of rigidity. Explain 
it first for rigid representations (for them one has $n_s=1$), with diagonal 
Jordan normal forms $J_j^n$. 

Denote by $\Xi (m)$ the set of 
all $(p+1)$-tuples of diagonal Jordan normal forms 
satisfying conditions {\em i)} and {\em ii)} from Theorem~\ref{generic1}, of 
size $\leq m$ (scalar Jordan normal forms are allowed). Find then all 
$(p+1)$-tuples of 
diagonal Jordan normal forms not from $\Xi (m)$ which after applying the 
map $\Psi$ result in a $(p+1)$-tuple from $\Xi (m)$. To this end one has to 
assume that 
every diagonal Jordan normal form of a $(p+1)$-tuple from $\Xi (m)$ has an 
eigenvalue of multiplicity 0 (which eventually was of positive multiplicity 
before applying $\Psi$). Thus one can obtain $\Xi (m+1)$ (one will have to 
exclude the $(p+1)$-tuples of size $>m+1$). 

Having obtained the set $\Xi (n)$, we leave only the $(p+1)$-tuples of 
size exactly $n$. Denote their set by $\Theta$. Explain how to 
obtain the analog of $\Theta$ obtained when the requirement the 
Jordan normal forms to be diagonal to be dropped. One constructs 
the set of $(p+1)$-tuples of Jordan normal forms $J_j^{n,1}$ such that 
the $(p+1)$-tuple of corresponding diagonal Jordan normal forms belongs to 
$\Theta$ -- this defines the set $\Phi$. After this one excludes from 
$\Phi$ the $(p+1)$-tuples not admitting generic eigenvalues -- this gives 
the necessary list for index of rigidity equal to 2.

For index of rigidity $h\leq 0$ one has to give first the list $L_h$ of 
$(p+1)$-tuples of diagonal Jordan normal forms satisfying condition 
$(\omega _n)$; after this one constructs the sets $\Xi (n)$, $\Theta$ and 
$\Phi$ by analogy with the case of index of rigidity equal to 2.  
 
For index of rigidity $h=0$ the list $L_0$ contains only three  
triples and one quadruple of diagonal Jordan normal forms, with 
MVs equal to a) (1,1), (1,1), (1,1), (1,1) b) (1,1,1), (1,1,1), (1,1,1) 
c) (1,1,1,1), (1,1,1,1), (2,2) 
d) (1,1,1,1,1,1), (2,2,2), (3,3). The reader will be able to deduce this 
fact from the beginning of the last section and from Lemma 3 from 
\cite{Ko2}. 

If one wants to obtain the list $L_h$ for $h\leq -2$, 
then one has to assume that $L_0$ contains four series of 
diagonal Jordan normal forms. They are a) (d,d), (d,d), (d,d), (d,d) 
b) (d,d,d), (d,d,d), (d,d,d) c) (d,d,d,d), (d,d,d,d), (2d,2d) d) 
(d,d,d,d,d,d), (2d,2d,2d), (3d,3d) where $d\in {\bf N}$ is such that the 
size of the matrices be $\leq n$. 

For negative indices of rigidity the list $L_h$ can be given by using again 
the ideas of the last section. First of all, we operate instad of 
with diagonal Jordan normal forms with their corresponding Jordan normal 
forms with a single eigenvalue. Having found the list $L_{h-2}$ one obtains 
the list $L_h$ by using {\em operations} $(s,l)$ (defined in {\bf 2} of 
Subsubsection~\ref{plan}) and the inverse of the {\em merging} (defined in 
{\bf 6} of that subsubsection. The details are left for the reader.   

\section{Definitions and notation\protect\label{defnot}}

\subsection{Levelt's result and its corollaries\protect\label{Leveltsresult}}

In \cite{L} Levelt describes the structure of the solution to a regular 
system at a pole:

\begin{tm}\label{Lv}
In a neighbourhood of a pole the solution to the  
regular linear system 

\begin{equation}\label{bs}
\dot{X}=A(t)X 
\end{equation}
can be represented in the form
\begin{equation}
\label{L*}
 X = U_j(t - a_j) (t-a_j)^{D_j} (t - a_j)^{E_j} G_j 
\end{equation}
where $U_j$ is holomorphic in a neighbourhood of the pole $a_j $,
$D_j = \mbox{diag} \, ( \varphi_{1,j} ,\ldots , \varphi_{n,j} )$,
$\varphi_{n,j} \in {\bf Z},$ $G_j\in GL(n,{\bf C})$. The matrix $E_j$ is in  
upper-triangular form and the real parts of its eigenvalues
belong to $[0,1)$ (by definition, $(t-a_j)^{E_j} = e^{E_j \ln
(t-a_j)}$). The numbers $\varphi_{k,j} $ satisfy the 
condition (\ref{L***})  formulated below.

System (\ref{bs}) is fuchsian at $a_j$ if and only if
\begin{equation}
\label{L**}
\det U_j (0) \ \neq 0
\end{equation}
We formulate the condition on $\varphi_{k,j} .$ Let $E_j$ have one and the 
same 
eigenvalue in the rows with indices $s_1<s_2<$$\ldots$$<s_q$. Then we 
have
\begin{equation}
\label{L***}
\varphi_{s_1,j} \geq \varphi_{s_2,j} \geq \ldots \geq \varphi_{s_q,j}
\end{equation}
\end{tm}

{\bf Remarks:} 1) Denote by $\beta _{k,j}$ the diagonal entries (i.e. the 
eigenvalues) of the matrix $E_j$. If the system is fuchsian, then the sums 
$\beta _{k,j}+\varphi _{k,j}$ are the eigenvalues $\lambda _{k,j}$ 
of the matrix-residuum $A_j$, see \cite{Bo1}, Corollary 2.1.

2) The numbers $\varphi _{k,j}$ are defined as valuations in the solution 
eigensubspace for the eigenvalue $\exp (2\pi i\beta _{k,j})$ of the 
monodromy operator, see the details in \cite{L}. These valuations can be 
defined on each subspace invariant for the monodromy operator. 

3) One can assume without loss of generality that equal eigenvalues of 
$E_j$ occupy consecutive positions on the diagonal and that the matrix 
$E_j$ is block-diagonal, with diagonal blocks of sizes equal to their 
multiplicities. The diagonal blocks themselves are upper-triangular.

4) An improvement of Levelt's form (\ref{L*}) can be found in \cite{Ko3}. 
More precisely, the fact that all entries of $E_j$ above the diagonal can be 
made equal to 0 or to 1.

5) Let $q>1$. The rest of the division of 
$\sum _{j=1}^{p+1}\sum _{k=1}^n\beta _{k,j}$ (which is an integer) by $q$ 
equals $m_0$. 

6) In \cite{Bo1} A.Bolibrukh proves the following lemma (see Lemma 3.6 
there):

\begin{lm}\label{bolibr}
The sum of the numbers $\varphi _{k,j}+\beta _{k,j}$ of system 
(\ref{bs}) corresponding to a subspace 
of the solution space invariant for all monodromy operators is a 
non-positive integer.
\end{lm}

\begin{cor}\label{Leveltred1}
Let the sum of the eigenvalues $\lambda _{k,j}$ of the matrices-residua of 
a fuchsian system corresponding to a subspace of the solution space 
of dimension $m$ invariant for all monodromy operators be 0. Denote 
these eigenvalues by $\lambda _{k,j}'$. Then there exists 
a change $X\mapsto RX$, $R\in GL(n,{\bf C})$ after which the system 
becomes block upper-triangular, the left upper block being of size $m$, the 
restrictions of the matrices-residua to it having eigenvalues 
$\lambda _{k,j}'$.
\end{cor}

A proof of the corollary can be found in \cite{Bo2}.

\subsection{Subordinate Jordan normal forms, canonical and strongly generic 
eigenvalues}

{\bf Definition.} Given two conjugacy classes $c'$, $c''$ with one and the 
same eigenvalues, of one and the same multiplicities, we say that $c''$ is 
{\em subordinate} to $c'$ if $c''$ lies in the closure of $c'$, i.e. for any 
matrix $A\in c''$ there exists a deformation $\tilde{A}(\varepsilon )$, 
$\tilde{A}(0)=A$ such that for $0\neq \varepsilon \in ({\bf C},0)$ one has 
$\tilde{A}(\varepsilon )\in c'$. Given two Jordan normal forms $J'$, $J''$, 
we say that $J''$ is {\em subordinate} to $J'$ if there exist conjugacy 
classes $c'\in J'$, $c''\in J''$ such that $c''$ is subordinate to $c'$. 

In the text we denote by $A_j$ matrices from $gl(n,{\bf C})$. They are 
often interpreted as matrices-residua of fuchsian systems on 
Riemann's sphere, i.e. systems of the form 

\begin{equation}\label{fuchsiansystem}
\dot{X}=(\sum _{j=1}^{p+2}A_j/(t-a_j))X
\end{equation}
If this system has no pole at infinity, then one has $A_1+\ldots +A_{p+2}=0$. 

{\bf Remark.} System (\ref{fuchsiansystem}) has $p+2$ poles because in what 
follows we have to realize the monodromy groups by fuchsian systems having 
an additional singularity with scalar local monodromy. When $A_{p+2}=0$ the 
system has $p+1$ poles.

{\bf Definition.} The eigenvalues $\lambda _{k,j}$ of the matrix-residuum 
$A_j$ are called 
{\em canonical} if none of the differences between two of its eigenvalues is 
a non-zero integer. The eigenvalues of the corresponding 
{\em monodromy operators} $M_j$ equal $\exp (2\pi i\lambda _{k,j})$. Hence, 
if the eigenvalues of $A_j$ are canonical, then to equal (to different) 
eigenvalues of the corresponding monodromy operator there correspond equal 
(different) eigenvalues of $A_j$. 

\begin{prop}\label{canonicalevs}
If the eigenvalues of the matrix $A_j$ are canonical, then one has 
$J(A_j)=J(M_j)$ and $M_j$ is conjugate to $\exp (2\pi iA_j)$.
\end{prop}

The proof can be found in \cite{Wa}.

{\bf Definition.} Eigenvalues $\lambda _{k,j}$ satisfying none of the 
equalities ($\gamma$) modulo ${\bf Z}$ (see the previous section) are called 
{\em strongly generic}; this definition is given only for eigenvalues 
$\lambda _{k,j}$; they are strongly generic if and only if the corresponding 
eigenvalues $\sigma _{k,j}=\exp (2\pi i\lambda _{k,j})$ are generic. 

{\bf Definition.} If $q>1$ and if the eigenvalues $\lambda _{k,j}$ satisfy 
none of the equalities ($\gamma$) modulo ${\bf Z}$ except $(\gamma ^*)$ 
and its multiples (see Subsection~\ref{qdxm}), then the eigenvalues are 
called {\em strongly relatively generic}. 

{\bf Definition.} Let the eigenvalues $\lambda _{k,j}$ not be strongly 
generic. Hence, there holds at least one equality of the form 

\[ \sum _{j=1}^{p+2}\sum _{k\in \Phi _j}\lambda _{k,j}=m~,~
m\in {\bf Z}~,~\sharp \Phi _1=\ldots =\sharp \Phi _{p+2}=\kappa ~,~
1\leq \kappa <n~~~~~(\gamma ')\] 
The number $|m|$ is 
called {\em distance of the set of eigenvalues to the relation 
$(\gamma ')$}. The minimal of the numbers $|m|$ (over all relations 
($\gamma '$)) is called {\em distance of the set of eigenvalues to the 
set of non-generic eigenvalues} (or just {\em distance}). For generic 
eigenvalues their distance is by definition equal to $\infty$.

\begin{lm}\label{distance}
{\bf A)} Let the eigenvalues $\sigma _{k,j}$ defined by the conjugacy 
classes $C_j$ be non-generic and either with a simple PMV or with a 
non-simple one $\xi$ 
being a primitive root of unity, and let at least one of the classes 
$C_j$ (say, $C_1$) have at least two different eigenvalues. Then for every 
$h\in {\bf N}$ sufficiently large there exist eigenvalues 
$\lambda _{k,j}$ with zero sum such that 

1) for all $k$, $j$ one has $\exp (2\pi i\lambda _{k,j})=\sigma _{k,j}$; 

2) for $j\leq p$ the eigenvalues $\lambda _{k,j}$ are canonical;

3) for $j=p+1$ if $\lambda _{k_1,j}-\lambda _{k_2,j}\in {\bf Z}$, then 
$\lambda _{k_1,j}-\lambda _{k_2,j}=0$ or $\pm 1$;

4) the distance of the eigenvalues is $\geq h$.

{\bf B)} If $\xi$ is not primitive (the rest of the conditions being like 
in A)), then there exist eigenvalues $\lambda _{k,j}$ satisfying conditions 
1), 2), 3) and 

4') for every relation $(\gamma ')$ which is not a multiple of 
$(\gamma ^*)$ its distance is $\geq h$.  
\end{lm}

Before proving the lemma we deduce from it 

\begin{cor}\label{cordistance}
If the eigenvalues $\lambda _{k,j}$ of the matrices-residua $A_j$ of system 
(\ref{fuchsiansystem}) are like in the lemma, with $h>n$, and if the 
$(p+1)$-tuple of matrices $A_j$ is irreducible, then the monodromy 
group of the system is with trivial centralizer.
\end{cor}

Notice that the irreducibility of the $(p+1)$-tuple of matrices $A_j$ is 
not automatic when $\xi$ is not primitive because there holds $(\gamma ^*)$.

{\em Proof:} 

$1^0$. Suppose first (see $1^0$ -- $3^0$) that $\xi$ is primitive. Assume 
that the monodromy group of system (\ref{fuchsiansystem}) 
satisfying the conditions of the corollary is with non-trivial centralizer 
${\cal Z}$. Then ${\cal Z}$ either 

a) contains a diagonalizable matrix $D$ with exactly two 
distinct eigenvalues or 

b) contains a nilpotent matrix $N$ with $N^2=0$. 

(Indeed, 
if $X\in {\cal Z}$, then every polynomial of the semi-simple or of the 
nilpotent part of $X$ belongs to ${\cal Z}$ which allows to construct such 
matrices $D$ or $N$.) 

$2^0$. In case a) one can assume that $D$ is diagonal -- 
$D=\left( \begin{array}{cc}\alpha I&0\\0&\beta I\end{array}\right)$. Hence, 
the monodromy operators are of the form 
$M_j=\left( \begin{array}{cc}M_j'&0\\0&M_j''\end{array}\right)$. Applying 
Lemma~\ref{bolibr} to the sums $\lambda '$, $\lambda ''$ of the eigenvalues 
$\lambda _{k,j}$ corresponding respectively to the $(p+1)$-tuples of blocks 
$M_j'$, $M_j''$, one obtains the inequalities $\lambda '\leq 0$, 
$\lambda ''\leq 0$. On the other hand, $\lambda '+\lambda ''=0$. 
Hence, $\lambda '=\lambda ''=0$ which contradicts the condition $h>0$. 

$3^0$. In case b) one can assume that 
$N=\left( \begin{array}{ccc}0&0&I\\0&0&0\\0&0&0\end{array}\right)$ and 
$M_j=\left( \begin{array}{ccc}M_j'&\ast &\ast \\0&M_j''&\ast \\
0&0&M_j'\end{array}\right)$ (the form of $N$ can be achieved by conjugation, 
the one of $M_j$ results from $[M_j,N]=0$). 

Denote by $\lambda ^1$, $\lambda ^2$, $\lambda ^3$ the sums of the 
eigenvalues $\lambda _{k,j}$ corresponding respectively to the 
$(p+1)$-tuples of upper blocks $M_j'$, of blocks $M_j''$, of lower blocks 
$M_j'$. One has $\lambda ^1+\lambda ^2+\lambda ^3=0$, $\lambda ^1\leq 0$, 
$\lambda ^1+\lambda ^2\leq 0$, 
$\lambda ^2+\lambda ^3\geq 0$ $\lambda ^3\geq 0$ (Lemma~\ref{bolibr}) and 
$|\lambda ^1-\lambda ^3|<n$ because the eigenvalues 
$\lambda _{k,j}$ are canonical for $j\leq p$ and for $j=p+1$ there holds 
condition 3) of the lemma (and the blocks $M_j'$, $M_j''$ are of sizes $<n$). 
Hence, either $|\lambda ^1+\lambda ^2|<n$ or $|\lambda ^2+\lambda ^3|<n$ 
which is a contradiction with the distance of 
the set of eigenvalues to the set of non-generic eigenvalues being $>n$. 

$4^0$. Suppose now that $\xi$ is not primitive. If the monodromy group of 
the system is with non-trivial centralizer, then in case a) one obtains 
$\lambda '=\lambda ''=0$. This is still not a contradiction with the 
condition $h>0$ because there remains the relation $(\gamma ^*)$. 

By Corollary~\ref{Leveltred1}, one can make a 
change $X\mapsto RX$ and block-triangularize the matrices-residua of the 
system. The change results in $A_j\mapsto R^{-1}A_jR$. This contradicts 
the condition the $(p+1)$-tuple of matrices $A_j$ to be irreducible.

In case b) one has $\lambda ^1+\lambda ^2+\lambda ^3=0$, $\lambda ^1\leq 0$, 
$\lambda ^1+\lambda ^2\leq 0$, 
$\lambda ^2+\lambda ^3\geq 0$ $\lambda ^3\geq 0$ (Lemma~\ref{bolibr}) and 
$|\lambda ^1-\lambda ^3|<n$ (like in $3^0$). This is possible only if 
$\lambda ^1=\lambda ^2=\lambda ^3=0$ and each of the last equalities 
results from $(\gamma ^*)$. By Corollary~\ref{Leveltred1}, one can 
block-triangularize the matrices-residua of the 
system by a change $X\mapsto RX$. This again contradicts the irreducibility 
of the $(p+1)$-tuple of matrices $A_j$.

The corollary is proved.

{\em Proof of Lemma \ref{distance}:} 

$1^0$. Prove A) (see $1^0$ -- $7^0$). Fix some canonical eigenvalues 
$\lambda _{k,j}^0$ satisfying 
condition 1) whose sum eventually is non-zero. Denote by $l$ the rest of 
the division of their sum (which is necessarily integer) by $n$. Decrease 
$l$ of the eigenvalues $\lambda _{k,p+1}^0$ by 1 (and don't change 
the other eigenvalues $\lambda _{k,j}$) -- this defines the 
eigenvalues $\lambda _{k,j}^1$. 

$2^0$. Fix integers $g_{k,j}$ such that 

a) to equal eigenvalues 
$\lambda _{k,j}^0$ there correspond equal integers $g_{k,j}$ and 

b) $\sum _{k=1}^n\sum _{j=1}^{p+1}(\lambda _{k,j}^1+g_{k,j})=0$. 

Set $\lambda _{k,j}^2=\lambda _{k,j}^1+g_{k,j}$. 

$3^0$. Denote by $\lambda '_j$, $\lambda ''_j$ two different eigenvalues 
$\lambda _{k,j}^0$, of multiplicities $m'$, $m''$. Change their 
corresponding numbers $g_{k,j}$ (denoted by $g_j'$, $g_j''$) respectively 
to $g_j'+um''$, $g_j''-um'$, $u\in {\bf N}$. This defines a 
new set of eigenvalues $\lambda _{k,j}^2$. (The other integers $g_{k,j}$ do 
not change.) 

$4^0$. For a given relation ($\gamma '$) (satisfied by the eigenvalues 
$\lambda _{k,j}^0$ instead of $\lambda _{k,j}$) call {\em quasi-multiplicity}  
of an eigenvalue $\lambda _{k,j}^0$ the number of eigenvalues 
$\lambda _{i,j}^0$ equal to $\lambda _{k,j}^0$ such that $i\in \Phi _j$ 
(the sets $\Phi _j$ were defined in Subsection~\ref{formulation}). 
For each given relation 
($\gamma ')$ the quasi-multiplicities of the eigenvalues $\lambda _{k,j}^0$ 
are not proportional to their multiplicities because either the PMV is 
simple or it is not but $\xi$ is a primitive root of unity. 

$5^0$. For each relation ($\gamma '$) and for each couple of eigenvalues 
$\lambda '_j$, $\lambda ''_j$ like in $3^0$ either 

1') for all values of $u\in {\bf N}$ the distance of the set of eigenvalues 
to the relation ($\gamma ')$ remains the same or 

2') only for finitely many of them this distance is $<h$. 

$6^0$. For a given relation ($\gamma '$) there exists a couple of eigenvalues 
$\lambda '_j$, $\lambda ''_j$ like in $3^0$ for which there holds 2'). Having 
chosen this couple, denote by $\Xi$ the set of relations ($\gamma '$) for 
which the chosen couple satisfies 2'). Hence, for the chosen couple one can 
find a finite subset $N_0$ of ${\bf N}$ such that for 
$u\in ({\bf N}\backslash N_0)$ the distance of the set of eigenvalues 
$\lambda _{k,j}^2$ to each of the relations from $\Xi$ is $\geq h$. Fix 
$u=u_0\in ({\bf N}\backslash N_0)$. This means that we change the 
eigenvalues $\lambda _{k,j}^2$.

$7^0$. If $\Xi$ is not the set of all relations ($\gamma '$), then choose 
such a relation not from it and repeat what was done in $6^0$. Every time 
we do this, the distance of the set of eigenvalues to more and more 
relations ($\gamma '$) becomes $\geq h$ and the distance to the rest of them 
does not change. As there are finitely many relations ($\gamma '$), after 
finitely many steps the distance of the set of eigenvalues $\lambda _{k,j}^2$ 
to the set of non-generic eigenvalues becomes $\geq h$, i.e. condition 
4) holds. Conditions 1) -- 3) hold by construction.

$8^0$. The proof of B) is completely analogous. One first divides the 
multiplicities of all eigenvalues by $(m_0,q)$, then constructs 
the eigenvalues $\lambda _{k,j}$ like in case A) and then one multiplies 
the multiplicities by $(m_0,q)$. The distances of $(\gamma ^*)$ and of 
its multiples remain 0, the distances of the other relations are 
$\geq (m_0,q)h>h$.

The lemma is proved.

\subsection{Corresponding Jordan normal forms 
and normalized chains\protect\label{CJNFNC}}

In \cite{Ko1} we define {\em correspondence} between Jordan 
normal forms. Namely, for each Jordan normal form $J_0$ we define its 
{\em corresponding} diagonal Jordan normal form $J_1$ as follows.

Suppose first that $J_0$ has a single eigenvalue 
of multiplicity $n$. 
Replace each Jordan block of $J_0$ of size 
$q'\times q'$ by a diagonal $q'\times q'$-matrix with $q'$ distinct 
eigenvalues where the last 
eigenvalues of the diagonal matrices replacing Jordan blocks of $J_0$ 
are the same, their last but one eigenvalues are the same and different 
from the last ones and for all $s\geq 0$ their last but $s$ eigenvalues 
(denoted 
by $h_s$) are the same and different from the last but $l$ ones for $l<s$. 

If $J_0$ has several eigenvalues, then this procedure is performed for every 
eigenvalue with the requirement eigenvalues of $J_1$ corresponding to 
different eigenvalues of $J_0$ to be different. 

{\bf Remark:} The definition of $J_1$ after $J_0$ can be described 
in equivalent terms like this. 
Call {\em Jordan normal form (JNF) of size $n$} a family 
$J^n=\{ b_{i,l}\}$ ($i\in I_l$, $I_l=\{ 1,\ldots ,s_l\}$, $l\in L$) of 
positive integers $b_{i,l}$ 
whose sum is $n$. Here $L$ is the set of eigenvalues (all distinct) and 
$I_l$ is the set of Jordan blocks with eigenvalue $l$, $b_{i,l}$ is the 
size of the $i$-th block with this eigenvalue. An $n\times n$-matrix 
$Y$ has the JNF $J^n$ (notation: $J(Y)=J^n$) if to its distinct 
eigenvalues $\lambda _l$, $l\in L$, there belong Jordan blocks of sizes 
$b_{i,l}$. 
For a given JNF $J^n=\{ b_{i,l}\}$ define its {\em corresponding} 
diagonal JNF ${J'}^n$. A diagonal JNF is  
a partition of $n$ defined by the multiplicities of the eigenvalues. 
For each $l$ $\{ b_{i,l}\}$ is a partition of $\sum _{i\in I_l}b_{i,l}$ and 
$J'$ is the disjoint sum of the dual partitions. 
Thus if for each fixed $l$ one has 
$b_{1,l}\geq$$\ldots$$\geq b_{s_l,l}$,  
then each eigenvalue 
$l\in L$ is replaced by $b_{1,l}$ new eigenvalues $h_{1,l}$, 
$\ldots$, $h_{b_{1,l},l}$ (hence, ${J'}^n$ has 
$\sum _{l\in L}b_{1,l}$ distinct eigenvalues). For $l$ fixed, set $g_k$ 
for the multiplicity of the eigenvalue $h_{k,l}$. Then the  
first $b_{s_l,l}$ numbers $g_k$ equal 
$s_l$, the next $b_{s_{l-1},l}-b_{s_l,l}$ equal $s_l-1$, $\ldots$, 
the last $b_{1,l}-b_{2,l}$ equal 1. 

{\bf Definition.} Two Jordan normal forms {\em correspond} to one another 
if they correspond to one and the same diagonal Jordan normal form (and 
thus the correspondence of Jordan normal forms is a relation of equivalence).

{\bf Definition.} We say that the eigenvalues of a Jordan matrix $G$ with 
Jordan normal form $J_1$ form a {\em normalized chain w.r.t. the Jordan 
normal form $J_0$} if for every eigenvalue of $J_0$ and for all 
$s$ one has $h_s-h_{s-1}\in {\bf N}^+$ and any two eigenvalues of $G$ 
corresponding to different eigenvalues of $J_0$ are non-congruent modulo 
${\bf Z}$. 

The following properties of corresponding Jordan normal forms are proved 
in \cite{Ko1}:

\begin{tm}\label{correspondingJNF}
1) If the Jordan normal forms $J_0$ and $J_1$ correspond to one another, 
then $r(J_0)=r(J_1)$ and $d(J_0)=d(J_1)$.

2) To each Jordan normal form $J'=\{ b'_{i,l}\}$ there corresponds a single 
Jordan normal form with a single eigenvalue $J''=\{ b''_{k,1}\}$. One has 
$b''_{k,1}=\sum _lb'_{k,l}$, $k=1,\ldots ,s_1$.  

3) Let the two Jordan normal forms ${J^n}'$ and 
${J^n}''$ correspond to one another. Choose in each of them an eigenvalue 
with maximal number of Jordan blocks. By 1) these 
numbers coincide. Denote them by $k'$. 
Decrease by 1 the sizes of the $k$ smallest Jordan blocks with these 
eigenvalues where $k\leq k'$. Then the two Jordan normal 
forms of size $n-k$ obtained in this way correspond to one another. 

4) Denote by $G^0$, $G^1$ two Jordan matrices with Jordan normal 
forms $J_0$, $J_1$ corresponding to each other, where $G^1$ is diagonal 
and $G^0$ is nilpotent and they are defined like $J_1$ and 
$J_0$ from the beginning of the subsection (when $J_1$ and $J_0$ 
are considered like matrices ($J_0$ being Jordan), not like Jordan normal 
forms). Then the orbits of the matrices $\varepsilon G^1$ and 
$G^0+\varepsilon G^1$ are the same for $\varepsilon \in {\bf C}^*$.

5) If $G^0$ is not necessarily nilpotent and if $G^0$, $G^1$ are defined 
by analogy with $J_0$, $J_1$, then the matrix 
$G^0+\varepsilon G^1$ is diagonalizable and its Jordan normal form is $J_1$ 
if $\varepsilon \in {\bf C}^*$ is small enough. 
\end{tm} 

\begin{cor}\label{samesequence}
Consider two $(p+1)$-tuples of Jordan normal forms -- $J_j^n$ and 
${J_j^n}'$ -- where for all $j$ the two Jordan normal forms $J_j^n$ and 
${J_j^n}'$ correspond to one another. Construct for each of the two 
$(p+1)$-tuples the $(p+1)$-tuples of Jordan normal forms 
$J_j^{n_{\nu}}$ and $J_j^{n_{\nu '}}$ like explained before 
Theorem~\ref{generic1}. Then one has $n_{\nu}=n_{\nu '}$ for all $\nu$, and 
for all $\nu$ and all $j$ the Jordan normal forms $J_j^{n_{\nu}}$ and 
$J_j^{n_{\nu '}}$ correspond to one another.
\end{cor}

The corollary follows from 1) and 3) of the above theorem. 

Denote by $J$ and $J'$ an arbitrary Jordan normal form and its corresponding 
Jordan normal form with a 
single eigenvalue. Consider a couple $D$, $D'$ of Jordan matrices with these 
Jordan normal forms. Suppose that the Jordan blocks of $D$ of sizes 
$b_{k,l}$ for $k$ 
fixed are situated in the same rows where is situated the Jordan block 
of size $\sum _lb_{k,l}$ of $D'$ (see 2) of the above theorem). Denote by 
$D_0$ the diagonal (i.e. semi-simple) part of the matrix $D$.

\begin{prop}\label{DD'} 
For all $\varepsilon \neq 0$ the matrix $\varepsilon D_0+D'$ is conjugate to 
$\varepsilon D$.
\end{prop}

{\bf Remark:} Permute the diagonal entries of $D_0$ so that before and after 
the permutation each entry remains in one of the rows of one and the same 
Jordan block of $D'$. Then the proposition holds again and the proof is the 
same, see Proposition 20 from \cite{Ko1}. 

{\bf Notation.} We denote by $a_j$ the poles of the fuchsian system 
(\ref{fuchsiansystem}) and by $\alpha _j$ the quantities $1/(a_j-a_{p+2})$. 
We often use the matrix ${\cal H}=$diag$(-1,\ldots , -1, 
0,\ldots ,0)$ ($m_0$ times $-1$ and $n-m_0$ times 0). 
The matrix with a single entry (equal to 1) in position $(i,j)$ is 
denoted by $E_{i,j}$. In all other cases double subscripts denote matrix 
entries.

\subsection{Plan of the paper}

The main difficulty in the case $d>1$ is the impossibility to realize 
the monodromy groups by 
fuchsian systems with canonical eigenvalues. If the eigenvalues 
$\lambda _{k,j}$ are represented in the form $\beta _{k,j}+\varphi _{k,j}$ 
with $\varphi _{k,j}\in {\bf Z}$ and with Re$(\beta _{k,j})\in [0,1)$, then 
the sum $\sum _{j=1}^{p+1}\beta _{k,j}$ is an integer (the sum of all 
eigenvalues $\lambda _{k,j}$ is 0, see (\ref{A_j})). If this integer is not 
divisible by $d$, then there exists no set of canonical 
eigenvalues $\lambda _{k,j}$ (canonical for each index $j$). 

Therefore we realize such monodromy groups by fuchsian systems having  
$(p+2)$-nd singularities with scalar local monodromy. More exactly -- 
equal to id (such singularities are called {\em apparent}. 
(This is why in the 
formulation of the problem we speak about $(p+1)$ matrices $M_j$ and here 
we consider $(p+2)$-tuples of matrices $A_j$.) The eigenvalues 
at the first $(p+1)$ singular points are canonical. The ones at the 
$(p+2)$-nd singularity compensate the rest $m_0$ of the division by $d$ of 
the above sum ($1\leq m_0\leq d-1$). 

The necessity to operate with systems 
with additional singularities explains why we begin the present second part 
of \cite{Ko1} by the study of 
$(p+2)$-tuples of   
matrices $A_j\in gl(n,{\bf C})$ 
with zero sum and satisfying certain linear equalities (see   
Subsection~\ref{S}). For a fuchsian system defined by such a $(p+2)$-tuple 
of matrices for fixed poles these equalities provide the 
necessary and sufficient condition the monodromy around one of the 
singularities of the fuchsian 
system (\ref{fuchsiansystem}) to be scalar. 

Section \ref{prelim} contains also the modification of the 
{\em basic technical tool} 
used in \cite{Ko1}. It is used to construct analytic deformations of 
$(p+2)$-tuples of matrices $A_j$ with trivial centralizer and satisfying 
certain linear equalities.

Section \ref{necessity} contains the proof of the necessity of conditions 
{\em i)} and {\em ii)} for the existence of irreducible $(p+1)$-tuples of 
matrices $M_j$ satisfying (\ref{M_j}) in the case $d>1$; 
Section~\ref{sufficiency} contains 
the proof of their sufficiency. In the proof of the sufficiency we use 
some results concerning the existence of irreducible $(p+1)$-tuples of 
nilpotent matrices with zero sum, see Section~\ref{existencenilp}.

\section{Some preliminaries\protect\label{prelim}}

\subsection{The basic technical tool\protect\label{basicTT}}

The {\em basic technical tool} is a way to deform a $(p+1)$-tuple of 
matrices $M_j$ or $A_j$ with a trivial centralizer into one with either new 
eigenvalues and the same Jordan normal forms of the respective matrices, or 
into one in which some of these Jordan normal forms $J_j^n$ are replaced by 
other Jordan normal forms ${J_j^n}'$ where $J_j^n$ is subordinate to 
${J_j^n}'$, or where ${J_j^n}'$ corresponds to $J_j^n$. Explain how 
it works in the multiplicative version first.

Given a $(p+1)$-tuple of matrices $M_j^1$ satisfying condition (\ref{M_j}) 
and whose centralizer is trivial, look for 
$M_j$ of the form 

\[ M_j=(I+\varepsilon X_j(\varepsilon ))^{-1}(M_j^1+
\varepsilon N_j(\varepsilon ))(I+
\varepsilon X_j(\varepsilon ))\] 
where the given matrices $N_j$ depend analytically on 
$\varepsilon \in ({\bf C},0)$ and one looks for $X_j$ analytic in 
$\varepsilon$. (One can set 
$M_j^1=Q_j^{-1}G_jQ_j$, $N_j=Q_j^{-1}V_j(\varepsilon )Q_j$ where $G_j$ are 
Jordan matrices and $Q_j\in GL(n,{\bf C})$.)
 
The matrices $M_j$ must satisfy 
equality (\ref{M_j}). In first approximation w.r.t. $\varepsilon$ this 
yields  

\[ \sum _{j=1}^{p+1}M_1^1\ldots M_{j-1}^1([M_j^1,X_j(0)]+N_j(0))
M_{j+1}^1\ldots M_{p+1}^1=0\]
or

\begin{equation}\label{P_j} 
\sum _{j=1}^{p+1}P_{j-1}([M_j^1,X_j(0)(M_j^1)^{-1}]+N_j(0)(M_j^1)^{-1})
P_{j-1}^{-1}=0
\end{equation}
with $P_j=M_1^1\ldots M_j^1$, $P_{-1}=I$ (recall that there holds 
(\ref{M_j}), therefore $M_j^1M_{j+1}^1\ldots M_{p+1}^1=P_{j-1}^{-1}$).

Condition (\ref{M_j}) implies that $\det M_1\ldots \det M_{p+1}=1$. 
There holds $\det M_j=\det M_j^1\det (I+\varepsilon (M_j^1)^{-1}N_j)$
$=(\det M_j^1)(1+\varepsilon$tr$((M_j^1)^{-1}N_j(0))+o(\varepsilon ))$. As 
$\det M_1^1\ldots \det M_{p+1}^1=1$, one has 
tr$(\sum _{j=1}^{p+1}(M_j^1)^{-1}N_j(0))$$=0$ (term of first order w.r.t.  
$\varepsilon$ in $\det M_1\ldots \det M_{p+1}$). 

Equation (\ref{P_j}) admits the following equivalent form: 

\begin{equation}\label{S_j} 
\sum _{j=1}^{p+1}([S_j,Z_j]+T_j)=0 
\end{equation}
with $S_j=P_{j-1}M_j^1P_{j-1}^{-1}$, 
$Z_j=P_{j-1}X_j(0)(M_j^1)^{-1}P_{j-1}^{-1}$, 
$T_j=P_{j-1}N_j(0)(M_j^1)^{-1}P_{j-1}^{-1}$. 

The centralizers of the 
$(p+1)$-tuples of matrices $M_j^1$ and $S_j$ coincide (to be checked 
directly), i.e. they are both trivial.  There holds

\begin{prop}\label{[A_j,X_j]}
The $(p+1)$-tuple of matrices $A_j$ is with trivial 
centralizer if and only if the mapping 
$(sl(n,{\bf C}))^{p+1}\rightarrow sl(n,{\bf C})$, 
$(X_1,\ldots ,X_{p+1})\mapsto \sum _{j=1}^{p+1}[A_j,X_j]$ is surjective.
\end{prop}

{\em Proof:} The mapping is not surjective if and only if the images of all  
mappings $X_j\mapsto [A_j,X_j]$ belong to one and the same proper linear 
subspace of $sl(n,{\bf C})$ which can 
be defined by a condition of the form 
tr$(D[A_j,X_j])=0$ for all $X_j\in sl(n,{\bf C})$ where 
$0\neq D\in sl(n,{\bf C})$. 
This amounts to tr$([D,A_j]X_j)=0$ for all 
$X_j\in sl(n,{\bf C})$, i.e. $[D,A_j]=0$ for $j=1,\ldots ,p+1$. 

The proposition is proved.

The mapping 

\[ (sl(n,{\bf C}))^{p+1}\rightarrow sl(n,{\bf C})~,~
(Z_1,\ldots ,Z_{p+1})\mapsto \sum _{j=1}^{p+1}[S_j,Z_j]\]
is surjective (Proposition \ref{[A_j,X_j]}). Recall that  
tr$(\sum _{j=1}^{p+1}(M_j^1)^{-1}N_j(0))=0$, i.e. 
tr$(\sum _{j=1}^{p+1}T_j)=0$. Hence, equation (\ref{S_j}) is solvable  
w.r.t. the unknown matrices $Z_j$ and, hence, equation (\ref{P_j}) is 
solvable w.r.t. the matrices $X_j(0)$. The implicit function 
theorem implies (we use the surjectivity here) that one can find $X_j$ 
analytic in $\varepsilon \in ({\bf C},0)$, i.e. one can find the necessary 
matrices $M_j$. 

In the additive version one has matrices $A_j^1=Q_j^{-1}G_jQ_j$ instead of 
$M_j^1$ and one sets 

\[ \tilde{A}_j=(I+\varepsilon X_j(\varepsilon ))^{-1}Q_j^{-1}
(G_j+\varepsilon V_j(\varepsilon ))Q_j(I+\varepsilon X_j(\varepsilon ))\] 
where $V_j(\varepsilon )$ are given 
matrices analytic in $\varepsilon$; one has 
tr$(\sum _{j=1}^{p+1}V_j(\varepsilon ))\equiv 0$. The matrices $X_j(0)$ 
satisfy the equation (which is of the form (\ref{S_j})) 

\[ \sum _{j=1}^{p+1}[A_j^1,X_j(0)]=-\sum _{j=1}^{p+1}Q_j^{-1}V_jQ_j~.\] 
The existence of $X_j$ analytic in $\varepsilon$ is justified like 
in the multiplicative version. 

\begin{lm}\label{diagifnondiag}
If for given Jordan normal forms $J_j^n$ of the matrices $M_j$ and generic 
eigenvalues there exists irreducible $(p+1)$-tuples of such matrices 
(satisfying (\ref{M_j})), then there exist such $(p+1)$-tuples 
for $M_j$ from the corresponding diagonal Jordan normal forms (and with 
generic eigenvalues). 
\end{lm}

{\em Proof:} $1^0$. Denote by $J_j'$ and $J_j''$ an arbitrary and its 
corresponding diagonal Jordan normal form. 
Let the $(p+1)$-tuple of matrices $M_j^0$ be irreducible, with 
generic eigenvalues and satisfying 
(\ref{M_j}). Set $M_j^0=Q_j^{-1}G_jQ_j$ where $G_j$ are Jordan matrices and 
$J(G_j)=J_j'$. 

$2^0$. Construct a deformation of the $(p+1)$-tuple of the form 

\[ M_j=(I+\varepsilon X_j(\varepsilon ))^{-1}Q_j^{-1}
(G_j+\varepsilon L_j)Q_j(I+\varepsilon X_j(\varepsilon ))\]
where $X_j$ are analytic in $\varepsilon \in ({\bf C},0)$ and $L_j$ are 
diagonal matrices with $J(L_j)=J_j''$. More exactly, assume that 
$G_j$ and $L_j$ are defined respectively like $G^0$ and $G^1$ from 
4) of Theorem~\ref{correspondingJNF}. By 5) of that theorem, one has 
$J(M_j)=J(L_j)$ for $\varepsilon \neq 0$ small enough.

$3^0$. The basic technical tool provides the existence of irreducible 
$(p+1)$-tuples of matrices $M_j$ for $\varepsilon \neq 0$ small enough. 
Their Jordan normal forms are $J_j''$ and their eigenvalues are generic.

The lemma is proved.

\subsection{Proof of Lemma \protect\ref{2n2}\protect\label{proofof2n2}}

$1^0$. Prove A) (see $1^0$ -- $5^0$). Suppose that there exists such a 
$(p+1)$-tuple of matrices $M_j$ with 
trivial centralizer but not irreducible. Then it can be conjugated to 
a block upper-triangular form with irreducible diagonal blocks whose 
eigenvalue satisfy the only non-genericity relation $(\gamma ^*)$ (defined 
in Subsection~\ref{qdxm}) and eventually some of its multiples. 

$2^0$. Consider two diagonal blocks and the representations $\Phi _1$, 
$\Phi _2$ (of sizes $m_1$, $m_2$) which they define. The Jordan normal forms 
of each of the matrices $M_j$ restricted to each of these blocks is a 
multiple of one and the same diagonal Jordan normal form and the ratio of 
the multiplicities of 
one and the same eigenvalues of $M_j$ as eigenvalues of $\Phi _1$ and 
$\Phi _2$ equals $m_1/m_2$ for every eigenvalue. 

$3^0$. For the dimensions $d_{j,i}$ of the conjugacy classes of the 
restrictions of $M_j$ to the two diagonal blocks one has 
$d_{1,i}+\ldots +d_{p+1,i}=2(m_i)^2$. A direct computation shows that 
dim Ext$^1(\Phi _1,\Phi _2)=0$. 

Indeed, consider the case when there are only 
two diagonal blocks and 
$M_j=\left( \begin{array}{cc}M_j^1&F_j\\0&M_j^2\end{array}\right)$. Hence, 
there exist matrices $G_j$ such that $F_j=M_j^1G_j-G_jM_j^2$. 

One has dim$E=2m_1m_2$ where 
$E=\{ (F_1,\ldots ,F_{p+1})\, |\, F_j=M_j^1G_j-G_jM_j^2\}$ (to be checked by 
the reader). The subspace $E'$ of $E$ defined by  

\[ \sum _{j=1}^{p+1}P_j=0~,~
P_j:=M_1^1\ldots M_{j-1}^1F_jM_{j+1}^2\ldots M_{p+1}^2\}\]  
(condition resulting from (\ref{M_j})) is of dimension $m_1m_2$. These 
conditions are linearly independent (the change of variables 
$(F_1,\ldots ,F_{p+1})\mapsto$$(P_1,\ldots ,P_{p+1})$ is bijective 
due to det$M_j^i\neq 0$).

$4^0$. One has dim$E'=$dim$E''$ where 
$E''=\{ (F_1,\ldots F_{p+1})\, |\, F_j=M_j^1G-GM_j^2, 
G\in M_{m_1,m_2}({\bf C})\}$. 
This is the space of blocks $F$ of size $m_1\times m_2$ resulting from 
the simultaneous congugation of the $(p+1)$-tuple of matrices 
$M_j^0=\left( \begin{array}{cc}M_j^1&0\\0&M_j^2\end{array}\right)$ by 
a matrix $\left( \begin{array}{cc}I&G\\0&I\end{array}\right)$. 

Thus dim$(E'/E'')=$dim Ext$^1(\Phi _1,\Phi _2)=0$.

$5^0$. This is true for every couple of diagonal blocks $\Phi _1$, 
$\Phi _2$. Hence, it is possible to conjugate the $(p+1)$-tuple to a 
block-diagonal form which contradicts the triviality of the centralizer. 

$6^0$. Prove B). Let $\Phi _1$ and $\Phi _2$ have the same meaning as above. 
Then dim Ext$^1(\Phi _1,\Phi _2)\geq 2$ 
and one can construct a semi-direct 
sum of $\Phi _1$, $\Phi _2$ which is not a direct one. Suppose that it is 
defined by matrices $M_j$ like in $3^0$. One can assume that the 
representations $\Phi _1$ and $\Phi _2$ are not equivalent (even if 
$m_1=m_2$) because for neither of them neither of conditions 
$(\alpha _{m_1})$, $(\alpha _{m_2})$ (which they satisfy) is an equality and 
there exist small deformations of the representations into 
nearby non-equivalent ones; when $(\alpha _n)$ is an equality, then 
such an irreducible representation is said to be {\em rigid}; it is 
unique up to conjugacy, see \cite{Si} and \cite{Ka}.

$7^0$. There exist infinitesimal conjugations of the 
matrices $M_j$ of the form 

\[ \tilde{M}_j=(I+\varepsilon X_j)^{-1}M_j(I+\varepsilon X_j)\] 
such that in first approximation w.r.t. $\varepsilon$ one has 
$\tilde{M}_1\ldots \tilde{M}_{p+1}=I$ and the conjugations do not 
result from a simultaneous infinitesimal conjugation of the matrices $M_j$. 
This follows from $d_1+\ldots +d_{p+1}\geq 2n^2+2$. 

The details look like 
this: set $X_j=\left( \begin{array}{cc}V_j&W_j\\U_j&S_j\end{array}\right)$. 
One can assume that $W_j=0$ because infinitesimal conjugations of $M_j$ with 
the matrices 
$\left( \begin{array}{cc}I&\varepsilon W_j\\0&I\end{array}\right)$ do not 
change the block upper-triangular form of the $(p+1)$-tuple. Hence, 

\[ \tilde{M}_j=M_j+\varepsilon 
\left( \begin{array}{cc}[M_j^1,V_j]+F_jU_j&F_jS_j-V_jM_j^2\\
M_j^2U_j-U_jM_j^1&[M_j^2,S_j]-U_jF_j\end{array}\right) =o(\varepsilon )\] 
Set 

\[ T=\{ (U_1,\ldots ,U_{p+1})|U_j\in M_{m_2,m_1}({\bf C}), 
\sum _{j=1}^{p+1}Y_j=0\}\]  
where $Y_j=M_1^2\ldots M_{j-1}^2U_jM_{j+1}^1\ldots M_{p+1}^1$ (the 
condition $\sum _{j=1}^{p+1}Y_j=0$ is the condition 
$\tilde{M}_1\ldots \tilde{M}_{p+1}=I$ restricted to the left lower 
$m_2\times m_1$-block and considered in first approximation w.r.t. 
$\varepsilon$). One has dim$T\geq m_1m_2+2$ (this results from 
$d_1+\ldots +d_{p+1}\geq 2n^2+2$). 

The subspace $T'$ of $T$ defined by the condition 
Tr$\sum _{j=1}^{p+1}F_jU_j=0$ is of codimension 1 in it, i.e. of 
dimension $\geq m_1m_2+1$. This is more than the size of the block $U$, i.e. 
more than the dimension of simultaneous infinitesimal conjugations with 
matrices $\left( \begin{array}{cc}I&0\\\varepsilon U&I\end{array}\right)$ 
(this is the subspace of $T$ of the form $\{ (U,\ldots ,U)\}$). Hence, 
one can choose $U_j$ from $T'/\{ (U,\ldots ,U)\}$ and 
after this choose $V_j$ and $S_j$ such that 

\[ \sum _{j=1}^{p+1}([M_j^1,V_j]+F_jU_j)=0~,~
\sum _{j=1}^{p+1}([M_j^2,S_j]-U_jF_j)=0\] 
The matrices $X_j$ define the infinitesimal deformations $\tilde{M}_j$.  

$8^0$. There exists also a true deformation of the form 

\[ \tilde{M}_j=(I+\varepsilon X_j+\varepsilon ^2Y_j(\varepsilon ))^{-1}
M_j(I+\varepsilon X_j+\varepsilon ^2Y_j(\varepsilon ))\]
with $Y_j$ analytic in $\varepsilon$.  The existence is justified by analogy 
with the basic technical tool and we leave the details for the reader. The 
triviality of the centralizer of the $(p+1)$-tuple of matrices $M_j$ 
makes the implicit function theorem applicable. 

Hence, for $\varepsilon \neq 0$ small enough the $(p+1)$-tuple of 
matrices $\tilde{M}_j$ is irreducible.      

The lemma is proved.

\subsection{The set ${\cal S}$\protect\label{S}}

Fix the distinct complex numbers $a_1$, $\ldots$, $a_{p+2}$. 
Consider the set ${\cal S}$ of 
$(p+2)$-tuples $A'$ of matrices 
$(A_1,\ldots ,A_{p+2})$ such that 

1) $A_{p+2}={\cal H}$, ${\cal H}$ was defined at the end of 
Subsection \ref{CJNFNC};

2) $A_1+\ldots +A_{p+2}=0$;

3) $(\sum _{j=1}^{p+1}\alpha _jA_j)|_{\kappa ,\nu }=0$, 
$\kappa =m_0+1,\ldots ,n$; $\nu =1,\ldots ,m_0$; 

4) $n=dn'$, $n'\in {\bf N}$, $d>1$ and $(d,m_0)=1$, see Subsection \ref{qdxm}; 

5) the Jordan normal forms of the matrices $A_j$ are fixed.

\begin{lm}\label{scalarMO}
Let the matrices $A_j$ satisfy conditions 1) and 2). Then 
the monodromy operator at $a_{p+2}$ of the fuchsian system 
(\ref{fuchsiansystem}) is scalar if and only if 
condition 3) holds (in which case it equals $I$).
\end{lm}

{\em Proof:} 

Represent system (\ref{fuchsiansystem}) locally, at $a_{p+2}$, by its 
Laurent series  

\[ \dot{X}=(A_{p+2}/(t-a_{p+2})+B+o(1))X\] 
where $B=-\sum _{j=1}^{p+1}\alpha _jA_j$ (to be checked directly). The 
change of variables (local, at $a_{p+2}$) 
$X\mapsto$diag$((t-a_{p+2})^{-1},\ldots ,(t-a_{p+2})^{-1},1,
\ldots ,1)X$ ($n-m_0$ units) brings system (\ref{fuchsiansystem}) to the 
form 

\[ \dot{X}=(A'_{p+2}/(t-a_{p+2})+O(1))X\]
with $A'_{p+2}=B'$ where the restriction of the matrix $B'$ to the 
left lower $((n-m_0)\times m_0$-block equals the one of $B$ to it and all 
other entries of $B'$ are 0. By Proposition~\ref{canonicalevs}, the 
monodromy operator $M_{p+2}$ is scalar if and only if $B'=0$. 

The lemma is proved.

{\bf Definition.} Call {\em canonical change of the eigenvalues of the 
matrices $A_j$ (CCE)} a 
change under which each eigenvalue changes by an integer, equal eigenvalues 
remain equal (hence, canonical or strongly generic eigenvalues remain such), 
the eigenvalues of $A_{p+2}$ (which are not canonical) 
do not change and the sum of all eigenvalues remains 0.

\begin{prop}\label{MBTT}
Let for given Jordan normal forms of the matrices $A_j$ and given set 
$\tilde{\lambda}$ of generic eigenvalues there exist 
$(p+2)$-tuples from ${\cal S}$. Then for all generic eigenvalues 
sufficiently close to $\tilde{\lambda}$ there exist $(p+2)$-tuples from 
${\cal S}$ with the same Jordan normal forms of the matrices $A_j$. 
\end{prop}

The proposition is proved in the next subsection. Its proof describes a way 
to construct $(p+2)$-tuples from ${\cal S}$ for nearby eigenvalues by 
deforming given $(p+2)$-tuples from ${\cal S}$ (with given eigenvalues). 
This way is called the {\em modified basic technical tool}.

Denote by ${\bf C}'$ the space of eigenvalues of the matrices $A_1$, 
$\ldots$, $A_{p+1}$ when their Jordan normal forms are fixed. 

\begin{cor}\label{Zariski}
If for given Jordan normal forms of the matrices $A_j$ the set ${\cal S}$ 
is not empty, then for all eigenvalues from some Zariski open dense 
subset of ${\bf C}'$ there exist $(p+2)$-tuples from ${\cal S}$.
\end{cor}

The corollary follows directly from the proposition. 

Denote by $\lambda ^0$ a point from ${\bf C}'$ defining for 
$j\leq p+1$ canonical eigenvalues. Consider all points 
from ${\bf C}'$ obtained from $\lambda ^0$ as a result of a CCE. 
Denote their set by $\Sigma (\lambda ^0)$. 

Choose $\lambda ^0$ such that 

a) for $j\leq p+1$ (one of) the eigenvalue(s) of 
$A_j$ of greatest multiplicity is integer; denote it by $\lambda ^0_j$; all 
other eigenvalues are non-integer (for $j\leq p+1$);

b) the only non-genericity relations modulo ${\bf Z}$ satisfied by the 
eigenvalues $\lambda _{k,j}$ are of the form 

\[ k(\sum _{j=1}^{p+1}\lambda ^0_j)+\delta =0~,~\delta \in {\bf Z}\]
(the eigenvalues of $A_{p+2}$ are integer, so we include them 
in $\delta$; $k\in {\bf N}$ does not exceed the smallest of the 
multiplicities of the eigenvalues $\lambda ^0_j$).

\begin{cor}\label{Zariski1}
If ${\cal S}$ is non-empty, then the set $\Sigma (\lambda ^0)$ 
contains a point $\lambda ^1$ for which the sum 
$\sum _{j=1}^{p+1}\lambda ^1_j$ ($\lambda ^1_j$ being an integer eigenvalue 
of $A_j$) is $>1$.
\end{cor}

{\em Proof:} 

$1^0$. Call a point from $\Sigma (\lambda ^0)$ {\em good} ({\em bad}) if 
${\cal S}$ projects on this point (if not). If a line $l_1\subset {\bf C}'$ 
passing through two points from 
$\Sigma (\lambda ^0)$ contains infinitely many bad points, then it contains 
only bad points from $\Sigma (\lambda ^0)$ (see the proposition and the 
corollary and remember that ${\cal S}$ is constructible). 

$2^0$. Suppose that no line parallel to $l_1$ and passing through two points from 
$\Sigma (\lambda ^0)$ contains only finitely many bad points. Hence, 
$\Sigma (\lambda ^0)$ must contain only bad points. The above proposition 
and corollary imply that ${\cal S}$ does not project on any line passing 
through two points of $\Sigma (\lambda ^0)$. Hence, it does not project on 
any point of any affine subspace of ${\bf C}'$ of dimension $k$ 
passing through $k$ points from $\Sigma (\lambda ^0)$ for $k=2,\ldots$, 
dim${\bf C}'$ (proved by induction on $k$). Hence, ${\cal S}$ must be 
empty -- a contradiction. 

$3^0$. This means that for every line $l_1$ passing through two points from 
$\Sigma (\lambda ^0)$ there exists a line $l_1'$ parallel to it,  
passing through two points from 
$\Sigma (\lambda ^0)$ and 
containing only finitely many bad points. Choose an index $j$ such that 
$A_j$ has at least two different eigenvalues ($\lambda '_j$ and 
$\lambda ''_j$, of multiplicities $m'$ and $m''$); one of the two 
eigenvalues is the integer eigenvalue of $A_j$. 

$4^0$. Denote by $t_1$, $t_2$ two points from ${\bf C}'$ such that the CCE 
which changes $t_1$ to $t_2$ is of the form 
$\lambda '_j\mapsto \lambda '_j-m''$, 
$\lambda ''_j\mapsto \lambda ''_j+m'$ (all other 
eigenvalues remaining the same). 
Hence, there exists a line $l$ in ${\bf C}'$ parallel to the one passing 
through $t_1$ and $t_2$, also passing through two points from 
$\Sigma (\lambda ^0)$ and such that $l$ contains only finitely many bad 
points. Hence, the line $l$ contains a point for which the sum 
$\sum _{j=1}^{p+1}\lambda ^1_j$ is a positive integer. 

The corollary is proved.

\subsection{The modified basic technical tool\protect\label{modifiedBTT}}

The {\em modified basic technical tool} is used to prove the existence of 
deformations (depending analytically on $\varepsilon \in ({\bf C},0)$) of 
$(p+2)$-tuples $A'\in {\cal S}$ with trivial centralizers. It will be used 
in different contexts and we explain it here in one of them (namely, the 
proof of Proposition~\ref{MBTT}). The 
basic points in the reasoning in all other contexts will be the same, there 
will be differences only in the technical details. 

{\bf Proof of Proposition \ref{MBTT}:}

$1^0$. Denote by $A_j^0=Q_j^{-1}G_jQ_j$ the matrices from an irreducible 
$(p+2)$-tuple $A'\in {\cal S}$, $G_j$ being Jordan matrices. 
We look for matrices of the form $A_{p+2}\equiv A_{p+2}^0=H$ ($H$ was 
defined in Section~\ref{defnot}), and for $j\leq p+1$ 

\[ A_j(\varepsilon )=(I+\varepsilon X_j(\varepsilon ))^{-1}Q_j^{-1}(G_j+
\varepsilon L_j)Q_j
(I+\varepsilon X_j(\varepsilon ))~,~\varepsilon \in ({\bf C},0)\] 
where $L_j$ are diagonal matrices which are polynomials of the semi-simple 
parts of the corresponding matrices $G_j$. Hence, $[G_j,L_j]=0$ and 
the matrices $G_j$ and $G_j+\varepsilon L_j$ define one and the same Jordan 
normal form for $\varepsilon$ small enough, see 4) -- 5) from 
Theorem~\ref{correspondingJNF}.  We want conditions 1) -- 5) 
from the previous subsection to hold. 

Obviously, one has 
$A_j=A_j^0+\varepsilon ([A_j^0,X_j(0)]+Q_j^{-1}L_jQ_j)+o(\varepsilon )$ for 
$j\leq p+1$. 

$2^0$. Conditions 2) and 3) yield in first approximation w.r.t. 
$\varepsilon$ the following system of equations linear in the unknown 
variables the entries of the matrices $X_j(0)$:

\begin{equation}\label{firstapproximation}
\sum _{j=1}^{p+1}[A_j^0,X_j(0)]=-\sum _{j=1}^{p+1}Q_j^{-1}L_jQ_j~~,~~
(\sum _{j=1}^{p+1}\alpha _j[A_j^0,X_j(0)])|_{\kappa ,\nu }=
-(\sum _{j=1}^{p+1}\alpha _jQ_j^{-1}L_jQ_j)
|_{\kappa ,\nu }
\end{equation}
with $\kappa =m_0+1,\ldots ,n$; $\nu =1,\ldots ,m_0$; 
$\alpha _j=1/(a_j-a_{p+2})$. The left hand-sides 
of these equations are linear forms in the entries of the matrices $X_j(0)$.

\begin{lm}\label{independent}
These linear forms are linearly independent.
\end{lm}

The proof of this lemma occupies the rest of the proof of the proposition. 
It implies the existence of $X_j$ analytic in $\varepsilon \in ({\bf C},0)$. 
Indeed, for $\varepsilon =0$  
equations (\ref{firstapproximation}) are 
solvable and the mapping 

\[ (X_1(0),\ldots ,X_{p+1}(0))\mapsto 
(\sum _{j=1}^{p+1}[A_j^0,X_j(0)],
(\sum _{j=1}^{p+1}\alpha _j[A_j^0,X_j(0)])|_{\kappa ,\nu })\]
is surjective onto $sl(n,{\bf C})\times {\bf C}^{(n-m_0)m_0}$. The existence of 
$X_j$ analytic in $\varepsilon$ small enough  
follows from the implicit function theorem. 
 
{\em Proof of the lemma:}

$1^0$. If the lemma were 
not true, then there should exist a couple of matrices 
$(0,0)\neq (V,W)\in sl(n,{\bf C})^2$ 
such that $W_{i,j}=0$ if $i>m_0$ or if $j\leq m_0$ and 

\[ {\rm tr}(V(\sum _{j=1}^{p+1}[A_j^0,X_j(0)])+
W(\sum _{j=1}^{p+1}\alpha _j[A_j^0,X_j(0)]))=0~~{\rm ~for~all~} 
X_j(0), ~{\rm i.e.}\] 
  
\[ {\rm tr}(([V,A_j^0]+\alpha _j[W,A_j^0])X_j(0))=0~~
{\rm identically~in~the~entries~of~}X_j(0), ~{\rm i.e.~} 
[V,A_j^0]+\alpha _j[W,A_j^0]=0.\]  

Summing up the equalities for $j=1,\ldots p+1$ and making use of 
$\sum _{j=1}^{p+1}A_j^0=-A_{p+2}^0$, one gets 

\begin{equation}\label{sum=0}
[V,A_{p+2}^0]+[W,-\sum _{j=1}^{p+1}\alpha _jA_j^0]=0
\end{equation}
The left lower  
$m_0\times (n-m_0)$-block of the second summand in (\ref{sum=0}) is 0. 
Hence, so is the left lower block of the first summand, and, hence, the 
one of $V$ itself (remember that $A_{p+2}^0={\cal H}$).

$2^0$. Choose $\gamma \in {\bf C}$ such that the matrix $V'=V+\gamma I$ be 
non-degenerate. Hence, for all $t\in {\bf C}\backslash \{ a_{p+2}\}$ the 
matrix $V'+W/(t-a_{p+2})$ is 
non-degenerate, block upper-triangular and 
det$(V'+W/(t-a_{p+2}))\equiv$det$V'$. Hence, 
its inverse is also block upper-triangular, with constant diagonal blocks 
and with constant non-zero determinant. 

$3^0$. Consider the fuchsian system 
$\dot{X}=A(t)X$ with $A(t)=(\sum _{j=1}^{p+2}A_j^0/(t-a_j))$. Perform in it 
the change of variables $X\mapsto (V'+W/(t-a_{p+2}))X$:  

\[ A(t)\rightarrow (V'+W/(t-a_{p+2}))^{-1}(W/(t-a_{p+2})^2)+
(V'+W/(t-a_{p+2}))^{-1}A(t)(V'+W/(t-a_{p+2}))\] 
(gauge transformation). The matrix 
$U(t)\stackrel{{\rm def}}{=}V'+W/(t-a_{p+2})$ being holomorphic and 
holomorphically invertible for 
$t\neq a_{p+2}$, the system remains fuchsian at $a_1$, $\ldots$, $a_{p+1}$ 
and has no singularities other than $a_j$. 

The equalities $[\gamma I+V,A_j^0]+\alpha _j[W,A_j^0]=0$ are equivalent 
to $[U(a_j), A_j^0]=0$ and imply that its 
residua at $a_j$, $j\leq p+1$, don't change. (The residua of the new system 
equal $(U(a_j))^{-1}A_j^0U(a_j)$.)

Check that the system remains fuchsian at $a_{p+2}$ as well. 
To this end observe that the matrix 
$U^{-1}=(V'(I+(V')^{-1}W/(t-a_{p+2}))^{-1}$ equals  
$(V')^{-1}-(V')^{-1}W(V')^{-1}/(t-a_{p+2})$ because 
$(I+(V')^{-1}W/(t-a_{p+2}))^{-1}$$=I-(V')^{-1}W/(t-a_{p+2})$ due to 
$((V')^{-1}W)^2=0$ (recall the block upper-triangular form of 
$V'$ and $W$). This implies that there are no polar terms 
of order higher than 1 at $a_{p+2}$ in the new system. 

Indeed, a priori the matrix $-U^{-1}\dot{U}+U^{-1}A(t)U$ (with $U$ and 
$U^{-1}$ as above) can have at $a_{p+2}$ a pole of order $\leq 3$. Its 
coefficient before $1/(t-a_{p+2})^3$ equals

\[ -((V')^{-1}W)^2-(V')^{-1}W(V')^{-1}A_{p+2}W\]
where each of the two summands is 0 (this follows from the form of the 
matrices $V$, $W$ and $A_{p+2}$). The one before $1/(t-a_{p+2})^2$ equals 

\[ (V')^{-1}W-(V')^{-1}W(V')^{-1}A_{p+2}V'+(V')^{-1}A_{p+2}W=
(V')^{-1}PV'~,~\]
\[ P=W(V')^{-1}-W(V')^{-1}A_{p+2}+A_{p+2}W(V')^{-1}\]
where $P=0$ for the same reason. 
Hence, the residuum at $a_{p+2}$ of the system also doesn't change (because 
the sum of all residua is 0). 

$4^0$. Hence, the solution $UX$ to the new system (which is 
the old one) equals $XC$, $C\in GL(n,{\bf C})$ (some solution to the old 
system). Make the analytic continuations of both solutions along one and the 
same closed contour. Hence, $UXM=XMC$ where $M$ is the 
monodromy operator corresponding to the contour. But 
$UXM=XCM$ which implies $[C,M]=0$. This is true for any 
closed contour. 

The monodromy group of the system being irreducible, the matrix $C$ must be 
scalar. Hence one has $UX\equiv CX$, i.e. 
$U\equiv C$ which implies $V'=C$, $W=0$. Recall that 
$V'=V+\gamma I$, $V\in sl(n,{\bf C})$. Hence, $V=0$.   
 
The lemma is proved, the proposition as well.

\section{Proof of the necessity\protect\label{necessity}}

\subsection{The proof in the case when $\xi$ is primitive}

Given a $(p+1)$-tuple with a trivial centralizer, one can deform it into one 
also with a trivial centralizer, with relatively generic eigenvalues and 
with the same Jordan normal forms,  
therefore we presume the eigenvalues relatively generic. 

We make use of Lemma \ref{diagifnondiag} and Corollary \ref{samesequence} 
and consider only the case when all 
matrices $A_j$ and $M_j$ are diagonalizable. Denote by $\Lambda ^n$ 
the PMVs of their eigenvalues; they are defined by the 
Jordan normal forms $J_j^n$. Denote by $\Lambda ^{n_{\nu}}$ the PMVs 
defined by the Jordan normal forms $J_j^{n_{\nu}}$ (their definition is 
given before Theorem~\ref{generic1}).

Consider the fuchsian system 

\begin{equation}\label{Fuchs1}
\dot{X}=(\sum _{j=1}^{p+2}A_j/(t-a_j))X
\end{equation}
where the $(p+2)$-tuple of matrices $A_j$ belongs to ${\cal S}$ and the 
eigenvalues of $A_j$ satisfy the conclusion of 
Corollary~\ref{Zariski1}. Hence, the $(p+2)$-tuple of matrices $A_j$ 
is irreducible and the eigenvalues of its monodromy 
operators satisfy the only non-genericity relation  

\[ \sigma _1\ldots \sigma _{p+2}=1~~~~~~~~~~~~~~~~(\gamma _0)\] 
where $\sigma _j=1$ (notice that 1 is 
the only eigenvalue of $M_{p+2}$).

\begin{lm}\label{reprtrivcentr1}
The monodromy group of system (\ref{Fuchs1}) with eigenvalues defined as 
above has a trivial centralizer. 
\end{lm}

The eigenvalues of 
system (\ref{Fuchs1}) satisfy the conditions of 
Lemma~\ref{distance} and Lemma~\ref{reprtrivcentr1} follows from 
Corollary~\ref{cordistance}. 

\begin{lm}\label{theMGisreducible1}
The monodromy group of system (\ref{Fuchs1}) with eigenvalues defined as 
in Lemma~\ref{reprtrivcentr1} can be conjugated to the form 
$\left( \begin{array}{cc}\Phi &\ast \\0&I\end{array}\right)$ where $\Phi$ is 
$n_1\times n_1$.
\end{lm}

{\bf Remark:} Notice that the subrepresentation $\Phi$ can be reducible. 

\begin{lm}\label{subreprtrivcentr1}
The centralizer $Z(\Phi )$ of the subrepresentation $\Phi$ is trivial.
\end{lm}

The last two lemmas are proved in Subsection \ref{proofsoflemmas}. We prove 
them also in the case when $\xi$ is not primitive because we need this for 
the next subsection.

The subrepresentation $\Phi$ being of dimension 
$n_1<n$, one can use induction on $n$ to prove the necessity. The 
induction base is the case when condition $(\omega _n)$ holds -- in this 
case there is nothing to prove. 

The PMV of the matrices $M_j'$ defining $\Phi$ equals 
$\Lambda ^{n_1}$. It follows from 
Lemma~\ref{subreprtrivcentr1} that for generic eigenvalues close to the ones 
of the matrices $M_j'$ defining $\Phi$ there exist irreducible 
$(p+1)$-tuples of 
diagonalizable matrices $\tilde{M}_j'\in GL(n_1,{\bf C})$ with PMV 
$\Lambda ^{n_1}$ and satisfying (\ref{M_j}) (this can be proved by using 
the basic technical tool in the multiplicative version).

Hence, if $\Lambda ^n$ is good, then $\Lambda ^{n_1}$ is good. 
Condition $(\omega _n)$ 
doesn't hold by assumption and conditions $(\alpha _n)$ and $(\beta _n)$ 
hold, see the Introduction. Finally, the PMV $\Lambda ^{n_s}$ is the same 
for $\Lambda ^n$ and for $\Lambda ^{n_1}$ (this follows from the definition 
of the PMVs $\Lambda ^{n_{\nu}}$). If 
$\Lambda ^{n_1}$ is good, then $\Lambda ^{n_s}$ satisfies condition 
$(\omega _{n_s})$ (one can't have $n_s=1$, see the remark after 
Theorem~\ref{generic1}). Hence, if the 
PMV $\Lambda ^n$ is good, then it satisfies conditions {\em i)} and 
{\em ii)} of Theorem~\ref{generic1}, i.e. they are necessary.

The necessity is proved. 

\subsection{The proof in the case when $\xi$ is not primitive and in the 
additive version}

$1^0$. If $\xi$ is not primitive, then the proof needs only small 
modifications. The eigenvalues of the matrices $A_j$ of system 
(\ref{Fuchs1}) are only 
relatively generic and satisfy only the non-genericity relation 
$(\gamma ^*)$ defined at the end of Subsection~\ref{qdxm} (and its 
multiples). 

$2^0$. Hence, the $(p+1)$-tuple of matrices-residua $A_j$ might be 
reducible. Suppose that it is in block upper-triangular form. Consider 
instead of the system its restriction to one of the diagonal blocks $P$; 
this restriction is presumed to be irreducible. 

For this restriction Lemma \ref{reprtrivcentr1} holds again. This 
follows from Corollary~\ref{cordistance}. The rest of the proof is the same 
because for all $j$ $J_j^n$ is a multiple of $J(A_j|_P)$. 
Recall that Lemmas~\ref{theMGisreducible1} and \ref{subreprtrivcentr1} 
are proved also in the case when $\xi$ is not primitive. 

$3^0$. In the additive version we proved the necessity in the case 
of generic eigenvalues in \cite{Ko1}. If the eigenvalues are non-generic 
and if $q=1$, then given such an irreducible $(p+1)$-tuple of matrices $A_j$ 
one can deform it by means of the basic technical tool into a nearby one 
with generic eigenvalues and the same Jordan normal forms of the respective 
matrices. Thus the necessity is proved for $q=1$ in the additive version.

$4^0$. Consider the additive version with $q>1$. Given 
such an irreducible $(p+1)$-tuple of matrices $A_j$, one can multiply it by 
$c\in {\bf C}^*$ to make the eigenvalues canonical. Next, by means of the 
basic technical tool one can deform it into a nearby irreducible one with 
the same Jordan normal forms of the matrices $A_j$ and with canonical 
relatively strongly generic eigenvalues. By 
Corollary~\ref{Leveltred1}, the monodromy group must be irreducible. 
By Proposition~\ref{canonicalevs}, for each $j$ one has $J(A_j)=J(M_j)$. 
Conditions {\em i)} and {\em ii)} are necessary in the multiplicative 
version, therefore they will be fulfilled in the additive one as well.   

\subsection{Proofs of Lemmas  
\protect\ref{theMGisreducible1} and \protect\ref{subreprtrivcentr1}
\protect\label{proofsoflemmas}}

{\em Proof of Lemma \ref{theMGisreducible1}:}

We assume that condition $(\omega _n)$ does not hold (otherwise there is 
nothing to prove).

$1^0$. Consider first (in $1^0$ -- $5^0$) the case when $\xi$ is a primitive 
root of unity. The monodromy group can 
be conjugated to a block upper-triangular form. The diagonal blocks define 
either irreducible or one-dimensional representations. The eigenvalues of 
each diagonal block $1\times 1$ 
satisfy the non-genericity relation $(\gamma _0)$; it is the only one 
satisfied by them due to the definition of the eigenvalues. This means that 
there is a single diagonal block of size $>1$. 

$2^0$. The block in the right  
lower corner must be of size 1. Indeed, by 
Lemma~\ref{bolibr} the left upper block can't be of size 1 (because 
the corresponding sum of eigenvalues $\lambda _{k,j}$ is a positive 
integer). Hence, 
it must be the only block of size $>1$ and the matrices $M_j$ 
look like this: 

\[ M_j=\left( \begin{array}{cc}M_j'&L_j\\0&I\end{array}\right) \hspace{3cm}
(M)\] 
where the size of $M_j'$ is $n'\times n'$ and the $(p+1)$-tuple of 
matrices $M_j'$ is irreducible. 

$3^0$. The block $M'$ must be of size $\leq n_1$. 
Indeed, if its size $n'$ is $>n_1$, then we show that the columns of the 
$(p+1)$-tuples of matrices $L_j$ aren't linearly independent modulo the 
space ${\cal W}$ defined below which will imply the non-triviality of the 
centralizer of the monodromy group. 

Denote by ${\cal W}\subset {\bf C}^{n'}$ the space of $(p+1)$-tuples of 
vector-columns of the 
form $(M_j-I)X$, $X\in {\bf C}^{n'}$. These are the vector-columns (right 
upper blocks) obtained 
by conjugating the $(p+1)$-tuple of matrices 
$\left( \begin{array}{cc}M_j'&0\\0&1\end{array}\right)$ by 
$\left( \begin{array}{cc}I&X\\0&1\end{array}\right)$. 

One has dim${\cal W}=n'$. Indeed, if dim${\cal W}<n'$, then the images 
of the linear operators $X\mapsto (M_j-I)X$ belong to a proper subspace of 
${\bf C}^{n'}$. This subspace can be assumed to belong to the space spanned 
by the first $n'-1$ vectors of the canonical basis of ${\bf C}^{n'}$ (which 
can be achieved by conjugating the monodromy operators $M_j$ by a 
block-diagonal matrix with diagonal blocks of sizes $n'$ and $n-n'$, the 
latter equal to $I$). 

But then the matrices $M_j'$ will be of the form 
$\left( \begin{array}{cc}M_j''&\ast \\0&1\end{array}\right)$ which 
contradicts the irreducibility of the matrix group generated by 
$M_1'$, $\ldots$, $M_{p+1}'$. 

$4^0$. The $(p+1)$-tuples of columns of the blocks $L_j$ belong to a 
subspace of ${\bf C}^{(p+1)n'}$ of dimension 
$\Delta :=r_1+\ldots +r_{p+1}-n'$ (each column of $L_j$ belongs to a space 
of dimension $r_j$ and there are $n'$ linear conditions satisfied by these 
columns; these conditions result from (\ref{M_j}) and look like this: 
$\sum _{j=1}^{p+1}M_1'\ldots M_{j-1}'L_j=0$; they are linearly independent 
because the change of variables $L_j'=M_1'\ldots M_{j-1}'L_j$ transforms 
them to $\sum _{j=1}^{p+1}L_j'=0$ and the independence in this form is 
evident).

$5^0$. The columns of the block $L$ must be linearly independent modulo the 
space ${\cal W}$, otherwise the monodromy group will be a direct sum of the 
form $M_j=\left( \begin{array}{cc}M_j'''&0\\0&1\end{array}\right)$, 
$M_j'''\in GL(n-1,{\bf C})$. Hence, 
$\Delta \geq (n-n')+n'$ (because the block $L$ has $n-n'$ columns and 
dim${\cal W}=n'$). Recall that $r_1+\ldots +r_{p+1}=n+n_1$. Hence, 
$n+n_1-n'\geq n-n'+n'$, i.e. $n'\leq n_1$.

$6^0$. Let now $\xi$ be a non-primitive root. 
The diagonal blocks of the monodromy group can be of two types. The first 
are of size 1, the eigenvalues satisfying the non-genericity relation 
$(\gamma _0)$. 

Describe the second type of diagonal blocks. Their sizes are $>1$ and can 
be different. Define the {\em unitary set} of eigenvalues: for each $j$ 
divide by $(m_0,q)$ the multiplicities of all eigenvalues 
$\sigma _{k,j}$ of the ones that are 
equal among themselves and are $\neq 1$. A 
block $F$ of the second type contains $h$ times the unitary set, 
$1\leq h\leq (m_0,q)$, and a certain number of eigenvalues equal to 1. 
(To different matrices $M_j$ there correspond, 
in general, different numbers of eigenvalues from the unitary set; therefore 
one must, in general, add some number of eigenvalues 1 for some values of 
$j$ to make the number of 
eigenvalues of the restrictions of the matrices $M_j$ to $F$  
equal; one then could eventually add one and the same number of eigenvalues 
equal to 1 to all matrices $M_j|_F$.) 
 
The eigenvalues of the blocks of the second type satisfy the 
non-genericity relation $(\gamma ^*)$ (defined in 
Subsection~\ref{qdxm}) and eventually $(\gamma _0)$ as well.

$7^0$. Denote by $\kappa (F)$ the 
ratio "number of eigenvalues $\sigma _{k,j}$ equal to 1"/"number of 
eigenvalues $\sigma _{k,j}$ not equal to 1" (eigenvalues of the 
restriction of the monodromy group to $F$), and by $\kappa _0$ the same 
ratio computed for the entire matrices $M_j$ (in both ratios one takes into 
account the eigenvalues of all matrices $M_j$). Then one must have 
$\kappa (F)<\kappa _0$.

Indeed, one can't have $\kappa (F)\geq \kappa _0$ 
because condition $(\omega _n)$ does not hold, hence, 
the restriction of the monodromy group to $F$ wouldn't satisfy this 
condition either. In the presence of the non-genericity relation 
$(\gamma _0)$ this implies a contradiction with the following 

\begin{lm}(see \cite{Ko1}). 
The following condition is necessary for the existence of irreducible 
$(p+1)$-tuples of matrices $M_j$ from the conjugacy classes $C_j$ and 
satisfying (\ref{M_j}):

\[ \min _{b_1,\ldots ,b_{p+1}\in {\bf C}^{\ast}, b_1\ldots b_{p+1}=1}
({\rm rk}(b_1M_1-I)+\ldots +{\rm rk}(b_{p+1}M_{p+1}-I))\geq 2n\]
\end{lm}

But then the sum of the eigenvalues $\lambda _{k,j}$ corresponding to the 
eigenvalues $\sigma _{k,j}$ from $F$ will be negative. If the block $F$ is 
to be in the right lower corner, then this sum must be positive 
(Lemma~\ref{bolibr}; it can't be 0 because the only non-genericity relation 
satisfied by the eigenvalues $\lambda _{k,j}$ is $(\gamma ^*)$ and its 
multiples, see Subsection~\ref{qdxm}; if the sum is 0, then by 
Corollary~\ref{Leveltred1} the $(p+1)$-tuple of matrices-residua would be 
block upper-triangular up to conjugacy -- a contradiction). Hence, the right 
lower block is of size 1. 

$8^0$. Denote by $\Pi$ the left upper $(n-1)\times (n-1)$-block. Conjugate 
it to make all non-zero rows of the restriction of 
the $(p+1)$-tuple $\tilde{M}$ of matrices $M_j-I$ to $\Pi$ linearly 
independent. After the conjugation some of the rows of the restriction of 
$\tilde{M}$ to $\Pi$ might be 0. In this case conjugate the matrices $M_j$ 
by one and the same permutation matrix which places the zero rows of 
$M_j-I$ in the last (say, $n-n'$) positions (recall that the last row of 
$M_j-I$ is 0, so $n-n'\geq 1$). Notice that if the restriction 
to $\Pi$ of a row of  
$M_j-I$ is zero, then the $n$-th position of the row is 0 as well, 
otherwise $M_j$ is not diagonalizable. 

$9^0$. After this conjugation the monodromy matrices have the form (M) from 
$2^0$, the matrix group generated by the matrices $M_j'$ is not 
necessarily irreducible, but can be conjugated to a block upper-triangular 
form, its restrictions to the diagonal blocks (all of sizes $>1$) being 
irreducible. Hence, the mapping 

\begin{equation}\label{(M-I)X} 
(X_1,\ldots ,X_{p+1})\mapsto \sum _{j=1}^{p+1}(M_j'-I)X_j~,~
X_j\in {\bf C}^{n'}
\end{equation}
is surjective onto ${\bf C}^{n'}$. Indeed, the matrix algebra ${\cal M}$ 
generated 
by the matrices $M_j'-I$ is block upper-triangular, its restrictions to 
each diagonal block (say, of size $u$) is irreducible and, hence, is 
$gl(u,{\bf C})$ (the Burnside theorem). Thus the algebra contains a 
non-degenerate matrix $L$. 

One has $L=\sum _{j=1}^{p+1}(M_j'-I)H_j$, $H_j\in {\cal M}$. For every 
vector-column $X\in {\bf C}^{n'}$ there exists a unique $Y\in {\bf C}^{n'}$ 
such that $X=LY$. Hence, one can set $X_j=H_jY$ which proves the 
surjectivity of the mapping.

$10^0$. If one defines the space ${\cal W}$ as above, then one finds that 
dim${\cal W}=n'$. This is proved like in the case when $\xi$ is a primitive 
root, see $3^0$, but the form 
$\left( \begin{array}{cc}M_j''&\ast \\0&1\end{array}\right)$ of the 
matrices $M_j'$ is forbidden not because the group generated by them 
must be irreducible (which, in general, is not true) but just because by 
definition there are no diagonal blocks of size 1. For the rest the proof 
goes on like in the case when $\xi$ is primitive.
 
This proves the lemma.

{\em Proof of Lemma \ref{subreprtrivcentr1}:} 

$1^0$. Consider first the case when $\xi$ is a primitive root (in $1^0$ -- 
$4^0$). If the lemma is not true, then $Z(\Phi )$ either 
A) contains a diagonalizable matrix $D$ with exactly two distinct 
eigenvalues or B) it contains a nilpotent matrix $N$ with $N^2=0$, see $1^0$ 
of the proof of Corollary~\ref{cordistance}.

$2^0$. In case A) one conjugates the monodromy group to the form 
$\left( \begin{array}{cc}G_j&L_j\\0&I\end{array}\right)$ with 
$G_j=\left( \begin{array}{cc}M_j^1&0\\0&M_j^2\end{array}\right)$.  
The sizes of $M_j^1$, $M_j^2$ equal the multiplicities of the two 
eigenvalues of $D$ and one has 
$D=\left( \begin{array}{cc}\alpha I&0\\0&\beta I\end{array}\right)$, 
$\alpha \neq \beta$, $D$ and $G_j$ are $n_1\times n_1$. 

At least one of 
the blocks $M_j^i$ must equal $I$ (for all $j$) because there is 
a single diagonal block of size $>1$. But this would mean that the 
monodromy group is a direct sum. Indeed, if $M_j^2=I$ for all $j$, then 
in the rows of the block $M_j^2$ and in the last columns (the ones of the 
block $I$) the entries of $M_j$ must be 0, otherwise $M_j$ is not 
diagonalizable. Being a direct sum contradicts 
Lemma~\ref{reprtrivcentr1}.

$3^0$. In case B) one 
can conjugate $N$ to the form 
$N=\left( \begin{array}{ccc}0&0&I\\0&0&0\\0&0&0\end{array}\right)$ (or 
$N=\left( \begin{array}{cc}0&I\\0&0\end{array}\right)$) with $I$ being 
$v\times v$, $v\leq n_1/2$; the second case corresponds to $v=n_1/2$. 

Hence, 
$M_j=\left( \begin{array}{cccc}M_j^1&R_j&T_j&Q_j\\0&M_j^2&S_j&H_j\\
0&0&M_j^1&P_j\\0&0&0&I\end{array}\right)$ where $M_j^1$ is $v\times v$ and 
if $v=n_1/2$, then the blocks $R_j$, $M_j^2$, $S_j$ and $H_j$ are absent. 

$4^0$. Like in case A) one shows that $M_j^1=I$ (this implies 
that in fact the possibility $v=n_1/2$ does not exist because there would be 
no diagonal block of size $>1$ at all). This means that the matrix 
with a single non-zero entry in position $(1,n)$ belongs to the centralizer 
of the monodromy group which contradicts Lemma~\ref{reprtrivcentr1}.

$5^0$. Let $\xi$ not be a primitive root. Suppose that we are in case A). 
We showed in the proof of the previous lemma (see $4^0$ of its proof) that 
the $(p+1)$-tuple of columns of the blocks $L_j$ (see $2^0$ of the present 
proof) belong to a subspace of 
${\bf C}^{(p+1)n_1}$ (denote it by ${\cal V}$) of dimension 
$\Delta =n+n_1-n_1=n$. When this subspace 
is factorized by the space ${\cal W}$ (see $3^0$ of the proof of the 
previous lemma), then the dimension becomes $n-n_1$. 

On the other hand, there are $n-n_1$ columns of the blocks $L_j$. In case 
A) the space ${\cal V}/{\cal W}$ splits into a direct sum of two such spaces 
defined for each of the blocks $M_j^1$ and $M_j^2$. The sum of their 
dimensions equals $n-n_1$, the number of columns of $L_j$; 
this implies that one can conjugate the matrices $M_j$ by 
$G\in GL(n,{\bf C})$ to the form 
$M_j=\left( \begin{array}{cccc}M_j^1&0&L_j'&0\\0&M_j^2&0&L_j''\\
0&0&I&0\\0&0&0&I\end{array}\right)$ which means that the monodromy group is 
a direct sum (one makes a self-evident permutation of the rows and columns 
which results from conjugation to achieve a block-diagonal form of the 
matrices). The matrix $G$ is block-diagonal, with diagonal blocks of 
sizes $n_1$ and $n-n_1$, the former equal to $I$. 

This is a contradiction with Lemma \ref{reprtrivcentr1}.

$6^0$. In case B) 
%the space ${\cal V}/{\cal W}$ splits into a direct 
%sum of two spaces, corresponding to the blocks $M_j^1$, $M_j^2$ and $M_j^1$. 
%Hence, 
a conjugation with 
a matrix $G'$ (defined like $G$ in $5^0$) brings the block 
$\left( \begin{array}{c}Q_j\\H_j\\P_j\end{array}\right)$ to the form 
$\left( \begin{array}{ccc}U_j^1&0&0\\0&U_j^2&0\\0&0&U_j^1
\end{array}\right)$ 
with heights of the blocks $U_j^1$, $U_j^2$ the same as the ones of 
$M_j^1$, $M_j^2$. 

This means that the matrix 
$\left( \begin{array}{cccccc}0&0&I&0&0&0\\
0&0&0&0&0&0\\0&0&0&0&0&0\\0&0&0&0&0&I\\0&0&0&0&0&0\\0&0&0&0&0&0
\end{array}\right)$ belongs to the centralizer of the monodromy group.
This is again a contradiction with Lemma~\ref{reprtrivcentr1}.

The lemma is proved.

\section{Proof of the sufficiency\protect\label{sufficiency}}

\subsection{Plan of the proof}

$1^0$. Two cases are possible:

Case A) $d_1+\ldots +d_{p+1}\geq 2n^2+2$

Case B) $d_1+\ldots +d_{p+1}=2n^2$

The condition $d>1$ excludes the possibility $d_1+\ldots +d_{p+1}=2n^2-2$, 
see the remarks after Theorem~\ref{generic1}. In case B) $\xi$ is presumed 
primitive.

In case A) we construct $(p+1)$-tuples of 
nilpotent matrices $A_j$ with trivial centralizers where for each $j$ 
$J(A_j)$ corresponds to the necessary Jordan normal 
form $J(M_j)$ of the monodromy operator $M_j$, see 
Lemma~\ref{constructnilp}, part I). Such a construction was 
already carried out in \cite{Ko1}. 

$2^0$. After this we deform the $(p+1)$-tuple into a nearby  
$(p+2)$-tuple 
of matrices $\tilde{A}_j$ with $J(\tilde{A}_j)=J(M_j)$ for $j\leq p+1$, with  
$A_{p+2}=\varepsilon {\cal H}$ and with 

\[ \tilde{A}_j=(I+\varepsilon X_j(\varepsilon ))^{-1}Q_j^{-1}
(G_j+\varepsilon L_j)Q_j(I+\varepsilon X_j(\varepsilon ))\] 
for $j=1,\ldots ,p+1$ so that the matrix $G_j+\varepsilon L_j$ be  
conjugate to $\varepsilon L_j$ (see 6) of 
Theorem~\ref{correspondingJNF}). We do this like in the modified basic 
technical tool. The following condition must hold:

\begin{equation}\label{p+2}
(\sum _{j=1}^{p+1}\tilde{A}_j/(a_j-a_{p+2}))_{\mu ,\nu}=0~,~
\mu =m_0+1,\ldots ,n~,~\nu =1,\ldots ,m_0. 
\end{equation}

$3^0$. Then one multiplies the $(p+2)$-tuple by $1/\varepsilon$; hence, the 
difference between the two eigenvalues of $\tilde{A}_{p+2}$ becomes equal to 
1. Condition (\ref{p+2}) implies that the singularity at $a_{p+2}$ of the 
fuchsian system with residua $\tilde{A}_j$ will be apparent (i.e. with 
local monodromy equal to $I$). 

$4^0$. In case B) one constructs $(p+2)$-tuples of nilpotent matrices 
$A_j$ with trivial centralizers, see Lemma~\ref{constructnilp}, part II). 
The Jordan normal form of the 
matrix $A_{p+2}$ has $m_0$ blocks of size 2 and $n-2m_0$ ones of size 1. 
One sets $\tilde{A}_{p+2}=A_{p+2}+\varepsilon {\cal H}$, the other matrices 
$\tilde{A}_j$ are defined like in case A). The matrix $A_{p+2}$ is such that 
for $\varepsilon \neq 0$ the matrix $\tilde{A}_{p+2}$ is conjugate to 
${\cal H}$. 
One has $\tilde{A}_{p+2}|_{\varepsilon =0}\neq 0$ to make possible the 
construction of a $(p+2)$-tuple of matrices $A_j$ with a trivial centralizer.  

Like in case A) one multiplies the residua $\tilde{A}_j$ by $1/\varepsilon$. 
Condition (\ref{p+2}) holds in case B) as well and the fuchsian system 
obtained after this multiplication has an apparent singularity at $a_{p+2}$. 

$5^0$. The deformation of $A_j$ into $\tilde{A}_j$ is possible to be done 
when certain transversality conditions make the 
theorem of the implicit function applicable.

\subsection{The basic lemma}

Denote by $\tilde{J}_j^n$ Jordan normal forms ($j=1,\ldots ,p+1$) satisfying 
conditions {\em i)} and {\em ii)} of Theorem~\ref{generic1}. 
Denote by $J_j^n$ both their corresponding Jordan normal forms  
with a single eigenvalue and 
the nilpotent conjugacy classes defining these Jordan normal forms. 

\begin{lm}\label{constructnilp}
Let $d>1$. Let the $p+1$ Jordan normal forms $J_j^n$ 
satisfy conditions {\em i)} and {\em ii)} of 
Theorem~\ref{generic1}. 

$I)$ Then in case A) there exists a $(p+1)$-tuple of nilpotent 
matrices $A_j$ satisfying (\ref{A_j}) such that 

1) for $j=1,\ldots ,p+1$ the matrix $A_j$ belongs either to $J_j^n$ 
or to its closure;

2) the centralizer ${\cal Z}$ of the $(p+1)$-tuple is trivial;

3) the $(p+1)$-tuple is in block upper-triangular form, the diagonal 
blocks being of sizes $n_s$, $n_{s-1}-n_s$, $\ldots$, $n-n_1$; the 
restriction of the $(p+1)$-tuple only to the first of them is non-zero;

4) the restriction to the diagonal block of size $n_s$ of the matrix 
$B:=\sum _{j=1}^{p+1}\alpha _jA_j$ is in Jordan normal form 
and its first $d$ eigenvalues are non-zero and simple; 

5) the restrictions to 
this block of the matrices $A_j$ define an irreducible representation.

$II)$ In case B) there exists a $(p+2)$-tuple of nilpotent 
matrices $A_j$ satisfying (\ref{A_j}) such that 1), 2) and 3) hold and 

6) the Jordan normal form of the matrix $A_{p+2}$ consists of $m_0$ blocks 
of size 2 and of $n-2m_0$ blocks of size 1; its non-zero entries are all in 
the diagonal block of size $n_s$; 

7) the restrictions to the diagonal block of size $n_s$ of the matrices 
$A_j$ are themselves block upper-triangular (and define a representation 
with a trivial centralizer), with diagonal blocks of equal 
size (which is 2, 3, 4 or 6) defining non-equivalent irreducible 
representations; 

8) the restriction to the diagonal block of size $n_s$ of the matrix $B$ 
is upper-triangular, with distinct non-zero eigenvalues.
\end{lm}

{\bf Remark:} The Jordan normal forms $J_j^{n_s}$ from case B) correspond to 
a triple or quadruple of matrices from the four so-called {\em special cases}, 
see Subsection~\ref{thebasicresult}.

\begin{cor}\label{corconstructnilp}
Let the $(p+1)$-tuple of conjugacy classes $C_j\in GL(n,{\bf C})$  
define Jordan normal forms $J_j^n$ which satisfy conditions {\em i)} and 
{\em ii)} of Theorem~\ref{generic1}. Let $d>1$ and let in case B) $\xi$ 
be primitive. Then 
there exist $(p+1)$-tuples of matrices $M_j\in C_j$ satisfying (\ref{M_j}) 
and with trivial centralizers.
\end{cor}

{\em Proof:}

$1^0$. Construct $(p+1)$- or $(p+2)$-tuples of nilpotent matrices $A_j^0$ 
like in the lemma where for $j=1,\ldots ,p+1$ $A_j^0$ is nilpotent and 
$J(A_j^0)=J_j^n$. Set $A_j^0=Q_j^{-1}G_jQ_j$ where $G_j$ are 
nilpotent Jordan matrices. Then construct matrices 

\[ \tilde{A}_j=(I+\varepsilon X_j(\varepsilon ))^{-1}Q_j^{-1}
(G_j+\varepsilon L_j)Q_j(I+\varepsilon X_j(\varepsilon ))\] 
such that for $j=1,\ldots ,p+1$ the matrix $G_j+\varepsilon L_j$ be  
conjugate to $\varepsilon L_j$ (see 6) of 
Theorem~\ref{correspondingJNF}). 

$2^0$. Choose the 
eigenvalues of the matrices $L_j$ to be canonical for $j\leq p+1$ 
and to have $\exp (2\pi iL_j)\in C_j$ (see Proposition~\ref{canonicalevs}). 
Set $L_{p+2}={\cal H}$.

\begin{lm}\label{abc}
There exist matrices $X_j$ analytic in $\varepsilon \in ({\bf C},0)$ 
such that 

a) $X_{p+2}\equiv 0$;

b) $\tilde{A}_1+\ldots +\tilde{A}_{p+2}=0$;

c) $(\sum _{j=1}^{p+1}\tilde{A}_j/(a_j-a_{p+2}))_{\mu ,\nu}=0$ 
for $\mu =m_0+1,n$, $\nu =1,\ldots ,m_0$.
\end{lm}

The lemma is proved in the next subsection. 

$3^0$. Fix $\varepsilon \neq 0$. Multiply the matrices 
$\tilde{A}_j$ by $1/\varepsilon$. Conditions a) and c) imply that the 
singularity at $a_{p+2}$ of the fuchsian system with poles $a_j$ and residua 
$\tilde{A}_j$ is apparent. Hence, the monodromy operators have the necessary 
eigenvalues and Jordan normal forms. 

Let $\xi$ be primitive. Choose the eigenvalues of the matrices $A_j$ 
generic and satisfying the conditions of Corollary~\ref{cordistance}. Hence, 
the centralizer of the monodromy group is trivial.
If in case A) $\xi$ is not primitive, then the existence of representations 
with trivial centralizers follows from Lemma~\ref{2n2}. 

The corollary is proved.

\subsection{Proof of Lemma \protect\ref{abc}}

The reader is supposed to have read and understood the proof of 
Lemma \ref{caseB} to which we refer.

$1^0$. Like in the proof of Proposition \ref{MBTT} one finds out that the 
existence of the necessary matrices $X_j$ follows from the solvability 
of the system of linear equations (\ref{firstapproximation}). So one has 
to prove Lemma~\ref{independent}. Notice that this time one has $A_{p+2}=0$, 
unlike in the conditions of Lemma~\ref{independent}, so the lemma has to be 
reproved. This is what we do.

$2^0$. In case A) equation (\ref{sum=0}) is reduced to 
$[W,\sum _{j=1}^{p+1}\alpha _jA_j^0]=0$ which implies $W=0$ (recall 
condition 4) from Lemma~\ref{constructnilp} and the fact that only the 
entries in positions $(\kappa ,\nu )$ of $W$ can be non-zero with 
$\kappa =1,\ldots ,m_0<d$, $\nu =m_0+1,\ldots ,n$). 

Hence, $[V,A_j^0]=0$ for $j=1,\ldots ,p+1$ which implies $V=0$ (recall 
condition 2) from Lemma~\ref{constructnilp}). 

$3^0$. Consider case B). Denote by $H_{\mu ,\nu}$ the block in position 
$(\mu ,\nu )$ when the matrices from $gl(n,{\bf C})$ are block-decomposed, 
with sizes of the diagonal blocks equal to 
$n_s$, $n_{s-1}-n_s$, $n_{s-2}-n_{s-1}$, $\ldots$, $n-n_1$. Recall that 
$m_0<d<n_s$. 

Equation (\ref{sum=0}) implies that $V_{H_{i,1}}=0$ for $i>1$. 
This follows from conditions 3) and 8) of Lemma~\ref{constructnilp} and from 
$W_{H_{i,1}}=0$ for $i>1$. The same equation implies that 
the restriction of $W$ only to $H_{1,1}$ can be non-zero. Indeed, one 
has $[V,A_{p+2}^0]_{H_{1,i}}=0$ for $i>1$, hence, 
$[W,\sum _{j=1}^{p+1}\alpha _jA_j^0]_{H_{1,i}}=0$. Conditions 3) and 8) 
of Lemma~\ref{constructnilp} imply that $W_{H_{1,i}}=0$. On the other hand, 
for $j>1$ one has $W_{H_{j,i}}=0$.

$4^0$. We prove (in $5^0$ -- $8^0$) that $V_{H_{1,1}}=W_{H_{1,1}}=0$. 
By $3^0$, this would imply $W=0$ and the
equations $[V,A_j^0]=0$, $j=1,\ldots ,p+1$ would imply $V=0$ (recall 
condition 2) of Lemma~\ref{constructnilp}). Hence, $V=W=0$ which proves 
Lemma~\ref{independent}. As we know that $V_{H_{i,1}}=W_{H_{i,1}}=0$ for 
$i>1$, we can replace in the equations 

\begin{equation}\label{VaW}
[V+\alpha _jW,A_j^0]=0
\end{equation}
the matrices 
$V$, $W$ and $A_j^0$ by their restrictions to $H_{1,1}$. This is what we do.

$5^0$. Block-decompose the matrices from $gl(n_s,{\bf C})$, all diagonal 
blocks being of size 2,3,4 or 6 (recall that the Jordan normal forms 
$J_j^{n_s}$ correspond to one of the {\em special cases} from 
Subsection~\ref{thebasicresult}). Denote the blocks by $U_{\mu ,\nu}$. The 
form of the matrices $A_{p+2}$ and $W$ implies that the restriction of 
equation (\ref{VaW}) to $U_{i,k}$ with $i\geq k$ looks like this: 

\[ V|_{U_{i,k}}A_j|_{U_{k,k}}-A_j|_{U_{i,i}}V|_{U_{i,k}}=0\] 
Hence, if $i>k$, then $V|_{U_{i,k}}=0$, if $i=k$, then 
$V|_{U_{i,k}}=\gamma _iI$. This follows from the non-equivalence of the 
irreducible representations defined by the matrices $A_j|_{U_{k,k}}$ and 
$A_j|_{U_{i,i}}$ and from Schur's lemma. 

$6^0$. The reader should have read the proof of 
Lemma~\ref{caseB} to which we refer. In this proof the blocks $F'$, 
$F''$, $F$ and $^tF$ were defined. 
Like in the proof of Lemma \ref{caseB} one shows that 
$V_{U_{\mu ,\nu}}=0$ if $\mu \neq \nu$, $U_{\mu ,\nu}\not\subset F$, and 
$V_{U_{\mu ,\nu}}=\gamma _{\mu}I$ if $\mu =\nu$. Let 
$U'=U_{\mu ,\nu}\subset F$. Then one has 

\begin{equation}\label{eq*}
A_j|_{U_{\mu ,\mu}}(V+\alpha _jW)|_{U'}-(V+\alpha _jW)|_{U'}
A_j|_{U_{\nu ,\nu}}+(\gamma _{\mu}-\gamma _{\nu})A_j|_{U'}=0
\end{equation}
Consider the fuchsian system with poles at $a_1$, $\ldots$ ,$a_{p+1}$ and 
residua $A_j^0$. Perform in it the change of variables 

\begin{equation}\label{eq**}
X\mapsto \tilde{W}X~,~\tilde{W}=\delta I+V+W/(t-a_{p+2})
\end{equation}
where $\delta \in {\bf C}$ is chosen such that $\Delta :=$det$\tilde{W}$ 
be non-zero (notice that $\Delta$ does not depend on $t$ due to the form 
of $V$ and $W$). The change (\ref{eq**}) transforms the residua as 
follows: $A_j^0\mapsto A_j^1=\tilde{W}(a_j)^{-1}A_j^0\tilde{W}(a_j)$, 
$j=1,\ldots ,p+1$. At $a_{p+2}$ a new singular point appears whose residuum 
is 0 (from the polar part only the term $O(1/(t-a_{p+2})^2)$ is non-zero). 

$7^0$. The block $U'$ of the residuum $A_j^1$ equals up to a non-zero factor 
$A_j|_{U_{\mu ,\mu}}(V+\alpha _jW)|_{U'}-(V+\alpha _jW)|_{U'}
A_j|_{U_{\nu ,\nu}}$. This non-zero factor depends on $\delta$. If 
$\gamma _{\mu}\neq \gamma _{\nu}$, then equation (\ref{eq*}) shows that the 
blocks $A_j|_{U'}$ are obtained as a result of the change (\ref{eq**}) in 
the fuchsian system with residua $A_j^0$. 

This, however, is impossible because the sum of the residua of a meromorphic 
1-form on ${\bf C}P^1$ is 0, so one should have 
$\sum _{j=1}^{p+1}A_j|_{U'}=0$; on the other hand, one has 
$\sum _{j=1}^{p+1}A_j|_{U'}=-A_{p+2}|_{U'}\neq 0$ by construction (see the 
proof of Lemma~\ref{caseB}) -- a contradiction. Hence, one has 
$\gamma _{\mu}=\gamma _{\nu}$ for all $(\mu ,\nu )$. As 
$V\in sl(n_s,{\bf C})$, one must have (for all $i$) $\gamma _i=0$. Thus 
only the entries of $V$ belonging to the block $F$ can be non-zero.

$8^0$. Sum up equations (\ref{eq*}) from $1$ to $p+1$. As 
$\sum _{j=1}^{p+1}A_j^0=0$ and $\sum \alpha _jA_j=-B$, one gets 
$[-A_{p+2},V]+[B,W]=0$. The form of the matrices $A_{p+2}$ and $V$ 
implies that $[-A_{p+2},V]=0$. The matrix $B$ is upper-triangular and has 
distinct eigenvalues. This and the form of $W$ implies $W=0$. 

But then one has $[A_j,V]=0$ for all $j$. Hence, 
$A_j|_{U_{\mu ,\mu}}V|_{U'}-V|_{U'}
A_j|_{U_{\nu ,\nu}}=0$. The $(p+1)$-tuples of matrices $A_j|_{U_{\mu ,\mu}}$ 
and $A_j|_{U_{\nu ,\nu}}$ defining non-equivalent irreducible 
representations, this implies $V=0$.

The lemma is proved.

\subsection{Proof of Lemma \protect\ref{constructnilp}}

$1^0$. Block decompose the matrices from $gl(n,{\bf C})$, the diagonal 
blocks being square, of sizes $n_s$, $n_{s-1}-n_s$, $n_{s-2}-n_{s-1}$, 
$\ldots$, $n-n_1$. Recall that we denote by $H_{\mu ,\nu}$ the block in 
position $(\mu ,\nu )$ (it is of size 
$(n_{s-\mu +1}-n_{s-\mu +2})\times (n_{s-\nu +1}-n_{s-\nu +2})$; we set 
$n_{s+1}=0$). Denote by $L_k$ the left upper block 
$n_{s-k+1}\times n_{s-k+1}$ (one has $L_1=H_{1,1}$).

$2^0$. Construct the matrices $A_j|_{L_1}$. For $j=1,\ldots ,p+1$ one has 
$J(A_j|_{L_1})=J_j^{n_s}$, for $j=p+2$ the Jordan normal form of 
$A_{p+2}|_{L_1}$ consists of $m_0$ blocks of size 2 and of $n_s-2m_0$ ones 
of size 1; the sizes $n_{\nu}$ being divisible by $d>m_0$, one has 
$n_s\geq 2d>2m_0$ (one can't have $n_s=d$ because this would mean that 
the matrices $A_j|_{L_1}$ are scalar -- a contradiction with condition 
$(\omega _{n_s})$).

\begin{lm}\label{csA}
In Case A) there exists an irreducible $(p+1)$-tuple of matrices 
$A_j^s:=A_j|_{L_1}$ satisfying condition (\ref{p+2}), with zero sum and 
with $J(A_j^s)=J_j^{n_s}$ for $j=1,\ldots ,p+1$. The matrix 
$B^s:=\sum _{j=1}^{p+1}\alpha _jA_j^s$ is in Jordan normal form, its first 
$d$ eigenvalues are simple and non-zero. 
\end{lm}

The lemma results from Theorem \ref{threematrices}, see Section 
\ref{existencenilp}.

\begin{lm}\label{caseB}
In Case B) there exists a $(p+2)$-tuple of matrices 
$A_j^s:=A_j|_{L_1}$ satisfying condition (\ref{p+2}), with zero sum, 
with trivial centralizer and 
with $J(A_j^s)=J_j^{n_s}$ for $j=1,\ldots ,p+1$; $J(A_{p+2}^s)$ consists 
of $m_0$ blocks of size 2 and of $n_s-2m_0$ blocks of size 1 and the matrix 
algebra generated by $A_1^s$, $\ldots$, $A_{p+2}^s$ contains a 
non-degenerate matrix. The matrix 
$B^s:=\sum _{j=1}^{p+1}\alpha _jA_j^s$ is upper-triangular and has 
distinct non-zero eigenvalues.  
\end{lm}

The lemma is proved in the next subsection.

$3^0$. Suppose that the matrices $A_j^{s-k+1}:=A_j|_{L_k}$ (whose sum 
is 0) are constructed such that 

a) for $j=1,\ldots ,p+1$ one has $A_j^{s-k+1}\in J_j^{n_{s-k+1}}$ or 
$A_j^{s-k+1}$ belongs to the closure of $J_j^{n_{s-k+1}}$;

b) in case B) $J(A_{p+2})$ consists of $m_0$ blocks of size 2 and of 
$n_{s-k+1}-2m_0$ blocks of size 1; in case A) $A_{p+2}=0$;

c) the matrices $A_j^{s-k+1}$ are block upper-triangular, with 
$A_j^{s-k+1}|_{H_{\mu ,\nu}}=0$ for $\mu >\nu$ and for $\mu =\nu >1$;

d) the columns of the $(p+1)$-tuple of matrices $A_j^{s-k+1}$ 
$j\leq p+1$ are linearly independent.

The last conditions means that if $(A_j^{s-k+1})^i$ denotes the $i$-th 
column of the matrix $A_j^{s-k+1}$ and if 
$\sum _{i=1}^{n_{s-k+1}}\beta _i(A_j^{s-k+1})^i=0$ for some constants 
$\beta _i\in {\bf C}$ and for $j=1,\ldots ,p+1$, then 
$\beta _1=$$\ldots$$=\beta _{n_{s-k+1}}=0$. 

{\bf Remark:} The entries of the matrix $A_{p+2}$ outside $L_1$ are 0. 
Therefore the construction from this moment on goes on like in \cite{Ko1} 
and we omit some details.

$4^0$. Construct the matrices $A_j^{s-k}$. Set $A_j^{s-k}|_{H_{k+1 ,\nu}}=0$ 
for $\nu =1,\ldots ,k+1$. Hence, condition c) will hold. 

Define the linear 
subspaces $S_{j,k}\subset {\bf C}^{n_{s-k+1}}$ of vector-columns as follows. 
Fix matrices $Q_j\in GL(n_{s-k+1},{\bf C})$ such that the matrices 
$T_j:=Q_j^{-1}A_j^{s-k+1}Q_j$ be in upper-triangular Jordan normal forms. 
Then for $j=1,\ldots ,p+1$ the space $S_{j,k}$ is spanned by the vectors 
of the form $Q_jV$ where the non-zero coordinates of $V$ can be only in the 
rows where $T_j$ has units and in the last rows of the smallest 
$r(J_j^{n_{s-k}})-$rk$(A_j^{s-k+1})$ Jordan blocks of $T_j$. (In some cases 
one has to choose part of the blocks of given size; there are no conditions 
imposed on the choice.)

For $j=1,\ldots ,p+1$ one has 
dim$(S_{j,k})=$rk$(A_j^{s-k+1})+r(J_j^{n_{s-k}})-$rk$(A_j^{s-k+1})=
r(J_j^{n_{s-k}})$. 

Denote by $D_k$ the union of blocks $H_{1,k+1}\cup \ldots \cup H_{k,k+1}$. 

\begin{lm}\label{Sjk}
Let for $j=1,\ldots ,p+1$ the matrix $A_j^{s-k+1}$ belong to the conjugacy 
class $J_j^{n_{s-k+1}}$ or to its closure. If the columns of 
$A_j^{s-k}|_{D_k}$ belong to the space $S_{j,k}$, then 
the matrix $A_j^{s-k}$ belongs to the conjugacy 
class $J_j^{n_{s-k}}$ or to its closure.
\end{lm}

The reader can find the proof of this fact in \cite{Ko1}. 

$5^0$. Denote by $\Gamma _k$ the linear space of $(p+1)$-tuples of 
vector-columns $Q_jV_j\in {\bf C}^{n_{s-k+1}}$ with zero sum where 
$V_j\in S_{j,k}$. One has dim$(\Gamma _k)\geq n_{s-k}$. Indeed,

\[ {\rm dim}(\Gamma _k)\geq {\rm dim}(S_{1,k})+\ldots +
{\rm dim}(S_{p+1,k})-n_{s-k+1}=\] 
\[ =r(J_1^{n_{s-k}})+\ldots +r(J_{p+1}^{n_{s-k}})-n_{s-k+1}=
n_{s-k}+n_{s-k+1}-n_{s-k+1}=n_{s-k}\]
We admit that the $n_{s-k+1}$ conditions arizing from 
$Q_1V_1+\ldots +Q_{p+1}V_{p+1}=0$ 
might not be linearly independent which explains why there is an inequality. 

Define $\Delta _k\subset {\bf C}^{n_{s-k+1}}$ as the space spanned by the 
$(p+1)$-tuples of vector-columns of the matrices $A_j^{s-k+1}$, $j\leq p+1$. 
Hence, dim$(\Delta _k)=n_{s-k+1}$ and 
dim$(\Gamma _k/\Delta _k)\geq n_{s-k}-n_{s-k+1}$. 

$6^0$. Hence, one can choose $n_{s-k}-n_{s-k+1}$ columns of the 
$(p+1)$-tuple of matrices $A_j^{s-k}|_{D_k}$ which belong to $\Gamma _k$ and 
are linearly independent modulo the space $\Delta _k$. With this choice 
the columns of the $(p+1)$-tuple of matrices $A_j^{s-k}$ will be linearly 
independent, i.e. condition c) from $3^0$ will hold. Condition b) holds by 
construction. Condition a) follows from Lemma~\ref{Sjk}. Thus conditions 
a) -- d) hold for every $k$ (and, in particular, for $k=s+1$). 

$7^0$. Prove that if the matrices $A_j=A_j^0$ are constructed like this, 
then the centralizer ${\cal Z}'$ of their $(p+2)$-tuple is trivial. 
Let $X\in {\cal Z}'$. Denote by $F$ a matrix from the 
algebra generated by the matrices $A_j$ whose restriction $F'$ to 
$L_1$ is non-degenerate. Such a matrix exists -- 
in Case A) this follows from the irreducibility of the algebra generated by 
the matrices $A_j^s$, in Case B) this follows from Lemma~\ref{caseB}.

The commutation relation $[X,F]|_{H_{s+1,1}}=0$ implies $X|_{H_{s+1,1}}=0$ 
(because $F|_{H_{\mu ,\nu }}=0$ if $\mu >\nu$ and if $\mu =\nu >1$, and 
one has $[X,F]|_{H_{s+1,1}}=(XF)|_{H_{s+1,1}}=X|_{H_{s+1,1}}F'=0$ and $F'$ 
is non-degenerate). 

The commutation relation $[X,F]|_{H_{s,1}}=0$ implies $X|_{H_{s,1}}=0$ 
(because $X|_{H_{s+1,1}}=0$ and 
one has $[X,F]|_{H_{s,1}}=(XF)|_{H_{s,1}}=X|_{H_{s,1}}F'=0$). 
Similarly $X|_{H_{k,1}}=0$ for $k>1$. 

For $k=1$ the commutation 
relations $[X,A_j]|_{H_{1,1}}=0$ are equivalent to $[X|_{H_{1,1}},A_j^s]=0$, 
hence, $X|_{H_{1,1}}=\alpha I$ by triviality of the centralizer of the 
$(p+2)$-tuple of matrices $A_j^s$, see Lemmas~\ref{csA} and \ref{caseB}.
Assume that $\alpha =0$.

$8^0$. After $7^0$ the commutation relations $[X,A_j]|_{H_{k,2}}=0$ 
($k=1,\ldots ,s+1$) become equivalent to $(A_jX)|_{H_{k,2}}=0$. They imply 
$X|_{H_{k,2}}=0$, otherwise the columns of the $(p+2)$-tuple of matrices 
$A_j$ would be linearly dependent -- a contradiction with condition d). 
In the same way one deduces that $X|_{H_{k,l}}=0$ for $k=1,\ldots ,s+1$; 
$l\geq 2$. Hence, $X=0$. Without the assumption $\alpha =0$, see $7^0$, this 
means that $X=\alpha I$, i.e. ${\cal Z}'$ is trivial. 

Conditions 1) -- 8) of the lemma follow immediately from the 
construction of the matrices $A_j$. The lemma is proved. 

\subsection{Proof of Lemma \protect\ref{caseB}}

$1^0$. Block-decompose the matrices from $gl(n_s,{\bf C})$, all diagonal 
blocks being of size $\chi =2$,3,4 or 6 (recall that the Jordan normal forms 
$J_j^{n_s}$ correspond to one of the {\em special cases} from 
Subsection~\ref{thebasicresult}). Denote the blocks by $U_{\mu ,\nu}$. 
Construct a block-diagonal $(p+1)$-tuple of matrices $A_j^0$ with diagonal 
blocks defining irreducible non-equivalent representations. If they are 
constructed like in examples (ex0) -- (ex7) from 
Subsection~\ref{someexamples}, see Lemma~\ref{secondmethod1} as well, then 
the matrix 
$B$ can have distinct non-zero eigenvalues (recall that one has the right to 
multiply the diagonal blocks by non-zero constants -- this preserves the 
conjugacy classes of $A_j$ and changes the eigenvalues of $B$). 

$2^0$. Define $A_{p+2}$. Denote by $F'$ the right upper block 
$m_0\times (n_s-m_0)$ (recall that $m_0<d\leq n_s/2$). If $m_0$ is divisible 
by the size $\chi$, then set $F=F'$ and $n'=n_s-m_0$. If not, then define 
$F''$ as obtained from $F'$ by deleting its first (i.e. left) $\chi -\psi$ 
columns where $\psi$ is the rest of the division of $m_0$ by $\chi$; set 
$n'=n_s-m_0-\chi +\psi$. 

Notice that the block $F''$ is a union of $[m_0/\chi ]$ rows each of 
$[n'/\chi ]$ entire blocks $U_{\mu ,\nu}$ and (when $\psi >0$) of a row of 
$[n'/\chi ]$ non-entire blocks $U_{\mu ,\nu}$ (of which only the first 
$\psi$ rows belong to $F''$). Denote by $F$ the block which is the union of 
all blocks $U_{\mu ,\nu}$ which entirely or partially belong to $F''$. 
Denote by $^tF$ the transposed to $F$. 

Require the restriction of $A_{p+2}$ to $F''$ to be of maximal rank and 
its restriction to every block $U_{\mu ,\nu}$ which entirely or 
partially belongs to $F''$ to be non-zero. 
Hence, $J(A_{p+2})$  
consists of $m_0$ blocks of size 2 and of $n-2m_0$ ones of size 1.

$3^0$. Define the matrices $A_j$, $j\leq p+1$. Their restriction to every 
diagonal block $U_{\mu ,\mu}$ equals the one of $A_j^0$ to it, their 
restriction to every block $U_{\mu ,\nu}\subset F$ is of the form 
$A_j^0|_{U_{\mu ,\mu}}D_{j;\mu ,\nu}-D_{j;\mu ,\nu}A_j^0|_{U_{\nu ,\nu}}$, 
$D_{j;\mu ,\nu}\in gl(\chi ,{\bf C})$, their other blocks are 0. 

The representations defined by $A_j^0|_{U_{\mu ,\mu}}$ and 
$A_j^0|_{U_{\nu ,\nu}}$ being non-equivalent for $\mu \neq \nu$, the map 

\[ (D_1,\ldots ,D_{p+1})\mapsto \sum _{j=1}^{p+1}
(A_j^0|_{U_{\mu ,\mu}}D_j-D_jA_j^0|_{U_{\nu ,\nu}})\] 
is surjective onto $gl(\chi ,{\bf C})$. Hence, one can choose the blocks 
$D_{j;\mu ,\nu}$ such that the sum of the matrices $A_j$ to be 0. 

Notice that for all $j$ the matrix $A_j$ is conjugate to $A_j^0$. 

$4^0$. There remains to be proved that the centralizer ${\cal Z}$ of the 
$(p+1)$-tuple of matrices $A_j$ is trivial. Denote by $U$ a block 
$U_{\mu _0,\nu _0}$ from $^tF$. For $X\in {\cal Z}$ the commutation relation 
$[A_{p+2},X]=0$ restricted to $U_{\mu _0,\nu _0}$ yields $X_U=0$. 

One has 
$A_j|_{U_{\mu _0,\mu _0}}X|_{U_{\mu _0,i}}-X|_{U_{\mu _0,i}}A_j|_{U_{i,i}}=0$, 
hence, if $i\neq \mu _0$, then $X|_{U_{\mu _0,i}}=0$, if $i=\mu _0$, 
then $X|_{U_{\mu _0,i}}=\gamma _iI$. In the same way one deduces that 
$X|_{U_{i,\nu _0}}=0$, if $i\neq \nu _0$, and 
$X|_{U_{i,\nu _0}}=\gamma _iI$ if $i=\nu _0$. 

Let $U'=U_{\mu _1,\nu _1}\subset F$. Set $X_{U'}=X'$, $A_j|_{U'}=A_j'$ and 
recall that $X_{U_{i,i}}=\gamma _iI$. The commutation relations restricted 
to $U'$ yield

\[ (\gamma _{\mu _1,\mu _1}-\gamma _{\nu _1,\nu _1})A_j'+
A_j|_{U_{\mu _1,\mu _1}}X'-X'A_j|_{U_{\nu _1,\nu _1}}=0\] 
Sum them up from 1 to $p+1$. One has $\sum _{j=1}^{p+1}A_j|_{U_{i,i}}=0$ for 
all $i$. Hence, the sum is 
$(\gamma _{\mu _1,\mu _1}-\gamma _{\nu _1,\nu _1})A_{p+2}|_{U'}$ which is 0 
only if $\gamma _{\mu _1,\mu _1}=\gamma _{\nu _1,\nu _1}$. But then 
one must have $A_j|_{U_{\mu _1,\mu _1}}X'-X'A_j|_{U_{\nu _1,\nu _1}}=0$ 
for all $j$. The non-equivalence of the representations defined by 
the blocks $A_j|_{U_{i,i}}$ for the different values of $i$ implies that 
$X'=0$. Hence, the centralizer is trivial.

The lemma is proved.

\section{On the existence of irreducible representations defined by 
nilpotent matrices\protect\label{existencenilp}}

\subsection{The basic result\protect\label{thebasicresult}}

In this section we consider the question of existence of irreducible 
$(p+1)$-tuples of nilpotent matrices from given conjugacy classes $c_j$, 
with given quantities $r_j=r(c_j)$ and with zero sum. This question was 
considered in \cite{Ko2}. We have to repeat some of the reasoning 
from \cite{Ko2} for two reasons:

-- there is an inexactitude in the formulation of the basic result from 
\cite{Ko2}, so we explain here what is correct and what is wrong;

-- we want to prove something more than the existence of such 
$(p+1)$-tuples.
 
Define as {\em special} the following cases. In each case the Jordan normal 
form defined by the class $c_j$ has Jordan blocks of one and the same size 
$l_j$. The special cases are

\[ \begin{array}{lllll}a)&p=3,&n=2g,&g>1,&l_1=l_2=l_3=l_4=2;\\
b)&p=2,&n=3g,&g>1,&l_1=l_2=l_3=3;\\ 
c)&p=2,&n=4g,&g>1,&l_1=l_2=4, ~l_3=2;\\ 
d)&p=2,&n=6g,&g>1,&l_1=6,~l_2=3,~l_3=2\end{array}\]
 
Define also some {\em  almost special} cases. They are obtained from one of 
the special ones (it is understood from the notation from which) by 
replacing 
in one of the Jordan normal forms a couple of blocks of size $l_j$ by two 
blocks of sizes $l_j+1$ and $l_j-1$. For these cases we give the sizes of 
the blocks:

\[ \begin{array}{lllll} a1)&J_1:3,1,2,\ldots ,2&J_2, J_3, J_4:2,\ldots ,2&
&n=2g, ~g>1\\
b1)&J_1:4,2,3,\ldots ,3&J_2:3,\ldots ,3&J_3:3,\ldots ,3&n=3g, ~g>1\\
c1)&J_1:4,\ldots ,4&J_2:4,\ldots ,4&J_3:3,1,2,\ldots ,2&n=4g, ~g>1\\ 
c2)&J_1:5,3,4,\ldots ,4&J_2:4,\ldots ,4&J_3:2,\ldots ,2&n=4g, ~g>1\\
d1)&J_1:6,\ldots ,6&J_2:3,\ldots ,3&J_3:3,1,2,\ldots ,2&n=6g, ~g>1\\
d2)&J_1:6,\ldots ,6&J_2:4,2,3,\ldots ,3&J_3:2,\ldots ,2&n=6g, ~g>1\\
d3)&J_1:7,5,6,\ldots ,6&J_2:3,\ldots ,3&J_3:2,\ldots ,2&n=6g, ~g>1
\end{array}\] 

{\bf Definition.} A representation defined by a $(p+1)$-tuple of 
matrices $A_j$ with zero sum is called {\em nice} if it is either 
irreducible or is reducible, with a trivial centralizer and the matrices 
$A_j$ admit a simultaneous conjugation to a block upper-triangular form 
in which the restrictions of the $(p+1)$-tuple to all diagonal blocks 
(they are all of sizes $>1$) define non-equivalent  
irreducible representations. (Hence, the algebra defined by the matrices 
$A_j$ contains a non-degenerate matrix.)

We prove in this section the following : 

\begin{tm}\label{threematrices}
1) Let the $(p+1)$-tuple of nilpotent conjugacy classes  
$c_j\in gl(n,{\bf C})$ be given 
with $r_1+\ldots +r_{p+1}\geq 2n$ and not 
corresponding to any of the special cases. 
Then there exists a $(p+1)$-tuple of matrices $A_j\in c_j$ defining a nice 
representation and there exists a 
$(p+1)$-tuple of distinct non-zero complex numbers $\alpha _j$ 
such that the matrix $B:=\alpha _1A_1+\ldots +\alpha _{p+1}A_{p+1}$ has 
a simple non-zero eigenvalue.

2) In the conditions of 1), if the $(p+1)$-tuple of nilpotent conjugacy 
classes $c_j\in gl(n,{\bf C})$ does not correspond to any of the 
almost special cases a1), b1), c2), or d3), then there exists an irreducible 
$(p+1)$-tuple of matrices $A_j\in c_j$ satisfying (\ref{A_j}). In all 7 
almost special cases one can obtain a matrix $B$ with distinct non-zero 
eigenvalues. 

3) If the conjugacy classes $C_j$ are unipotent, satisfying 
$r_1+\ldots +r_{p+1}\geq 2n$ and not 
corresponding to any of the special cases,  
then there exists a $(p+1)$-tuple of matrices $M_j\in C_j$ defining a nice 
representation and satisfying (\ref{M_j}). If they do not correspond to any 
of the almost special cases a1), b1), c2), or d3), then there exists an 
irreducible $(p+1)$-tuple of matrices $M_j\in C_j$ satisfying (\ref{M_j}).
\end{tm}

{\bf Remarks:} 1) Condition $(\omega _n)$ is necessary for the existence of 
nice representations -- it must hold for the restrictions of the matrices to 
each of the diagonal blocks and the rank of a nilpotent matrix is $\geq$ the 
sum of the ranks of its restrictions to these blocks.

2) In \cite{Ko2} the existence of irreducible representations in cases 
a1), b1), c2) and d3) is claimed which is not true. The proof of the basic 
result from \cite{Ko1} (which uses the result from \cite{Ko2}) needs only 
the existence of nice representations (not necessarily irreducible) and is 
performed with the corrected result in exactly the same way as it is done in 
\cite{Ko1}. The rest of the results from \cite{Ko2} are correct.

The theorem is proved in Subsection \ref{proofofthreematrices}. We don't 
prove part 3) of it which follows from parts 1) and 2) when one considers 
as matrices $M_j$ the monodromy operators of a fuchsian system with 
matrices-residua $A_j$ satisfying the conditions of part 1) or 2) 
(remember Proposition~\ref{canonicalevs}). Condition 
$(\omega _n)$ is necessary for the existence of irreducible 
$(p+1)$-tuples of nilpotent matrices satisfying (\ref{A_j}) or of 
unipotent matrices $M_j$ satisfying (\ref{M_j}), see \cite{Ko2}. 

We precede the 
proof of the lemma by the construction of examples of irreducible  
triples or quadruples of matrices $A_j$ in some particular cases, see the 
next subsection, and by deducing 
Corollary~\ref{threematrices1} from the theorem. We refer to the examples as 
to (ex1), (ex2) etc. 
The irreducibility of the triples or quadruples from each example is 
checked by proving that the matrix algebra generated 
by the matrices is $gl(n,{\bf C})$. To this end one first finds a matrix $S$ 
from the algebra with a single non-zero entry (by multiplying in a suitable 
order the matrices) and then again by suitable multiplications of the matrix 
$S$ by the matrices $A_j$ one 
obtains all other elements of the canonical basis of $gl(n,{\bf C})$, see 
examples of applications of this technique in \cite{Ko2}.

\begin{cor}\label{threematrices1}
Let $d>1$ and let the $(p+1)$-tuple of nilpotent conjugacy classes  
$c_j\in gl(n,{\bf C})$ be given 
with $r_1+\ldots +r_{p+1}\geq 2n$ and not corresponding to 
any of the special cases. Then there exists an irreducible $(p+1)$-tuple of 
matrices $A_j\in c_j$ satisfying 
(\ref{A_j}) and a $(p+1)$-tuple of distinct non-zero complex numbers 
$\alpha _j$ such that the matrix 
$B:=\alpha _1A_1+\ldots +\alpha _{p+1}A_{p+1}$ has 
at least $d$ simple non-zero eigenvalues.
\end{cor}

{\em Proof:} 

$1^0$. Denote by $c_j^1\in gl(n/d,{\bf C})$ the conjugacy classes obtained 
from $c_j$ by 
reducing the size of the matrices and the number of Jordan blocks of a given 
size $d$ times. Denote by $c_j^l$ the conjugacy classes obtained from 
$c_j^1$ by increasing the size of the matrices and the number of Jordan 
blocks of a given size $l$ times. In particular, $c_j^d=c_j$. 

$2^0$. Suppose first that the conjugacy classes $c_j^1$ don't correspond to 
an almost special case (they don't correspond to a special case, even with 
$g=1$, because in this case $c_j$ would also correspond to such a case). We 
show that there exist irreducible $(p+1)$-tuples of matrices 
$A_j^l\in c_j^l$. For $l=1$ this follows from the above theorem. Suppose 
that it is true for $l\leq l_0$. Then there exists an irreducible 
triple of matrices $A_j^{l_0}\in c_j^{l_0}$ whose sum is zero. 

$3^0$. Consider the 
matrices $\left( \begin{array}{cc}A_j^{l_0}&0\\0&A_j^1\end{array}\right)$. 
The irreducible representations $L^{l_0}$ and $L^1$ defined by the triples 
of matrices $A_j^{l_0}$ and $A_j^1$ satisfy the condition 
dim Ext$^1(L^{l_0},L^1)\geq 2$. Indeed, one has 
dim Ext$^1(L^{l_0},L^1)=l_0(\sum _{j=1}^{p+1}d(c_j^1)-2(n/d)^2)\geq 2$ 
because $\sum _{j=1}^{p+1}d(c_j^1)\geq 2n^2+2$ (see \cite{Ko2}, Lemma 3). 
One subtracts $(n/d)^2$ twice -- once to factor out conjugations with 
block upper-triangular matrices and once because the sum of the matrices 
is 0.

\begin{lm}\label{Ext>=2}
1) If for the non-equivalent irreducible representations $L'$, $L''$  
defined by the matrices 
$A_j'\in c_j'$, $A_j''\in c_j''$ (each satisfying (\ref{A_j})) one has 
dim Ext$^1(L',L'')\geq 1$, then there exists a representation 
defined by matrices $A_j'''$ satisfying (\ref{A_j}), with 
$A_j'''\in c_j'\oplus c_j''$, which is their semi-direct sum. If 
dim Ext$^1(L',L'')\geq 2$, then there exists an irreducible such 
representation. 
The conjugacy classes $c_j'$, $c_j''$ are 
arbitrary (not necessarily nilpotent).

2) Let the conjugacy classes $c_j'$, $c_j''$ be nilpotent. Let the matrices 
$B^{(i)}=\sum _{j=1}^{p+1}\alpha _jA_j^{(i)}$, $i=1,2$ have respectively 
$m'$ and $m''$ simple non-zero eigenvalues. If the matrices $A_j^{(i)}$ are 
like in 1) and if 
dim Ext$^1(L',L'')\geq 2$, then there exists an irreducible 
representation defined by matrices $A_j$ satisfying (\ref{A_j}), with 
$A_j\in c_j'\oplus c_j''$ and such that the matrix $B$ (defined like 
$B'$, $B''$) has at least $m'+m''$ simple non-zero eigenvalues. 
\end{lm}

The lemma is proved after the corollary. It implies the existence of 
irreducible triples of matrices $A_j^{l_0+1}\in c_j^{l_0+1}$ with at least 
$l_0+1$ simple non-zero eigenvalues. 

$4^0$. If the conjugacy classes $c_j^1$ correspond to an almost special 
case, then for $l>1$ the classes $c_j^l$ correspond to a $(p+1)$-tuple of 
conjugacy classes which contains in its closure a {\em neighbouring} case 
(their definition is given in {\bf 8.} of Subsubsection~\ref{plan}). 
For such cases we prove that there exist irreducible $(p+1)$-tuples 
with $B$ having distinct non-zero eigenvalues, see 
Subsubsection~\ref{NC}. Hence, one can choose an irreducible $(p+1)$-tuple 
of matrices $A_j\in c_j^l$ which is close to a neighbouring case and 
$B$ will still be with distinct non-zero eigenvalues.

The corollary is proved. 

{\em Proof of Lemma \ref{Ext>=2}:}

$1^0$. As dim Ext$^1(L',L'')\geq 1$, there exists a semi-direct sum of 
$L'$ and $L''$ which is 
not reduced to a direct one, i.e. there exists a $(p+1)$-tuple of matrices 
$A_j'''=\left( \begin{array}{cc}A_j'&G_j\\0&A_j''\end{array}\right)$ which 
by simultaneous conjugation can't be transformed into a block-diagonal one. 

$2^0$. The representations $L'$, $L''$ are not equivalent which means 
that the centralizer ${\cal Z}$ of the $(p+1)$-tuple of matrices $A_j'''$ is 
trivial (if 
$\left( \begin{array}{cc}U&M\\N&P\end{array}\right) \in {\cal Z}$, 
then $NA_j'-A_j''N=0$; as $L'$ and $L''$ are irreducible and 
non-equivalent, one has $N=0$; then $[A_j',U]=0$ and $[A_j'',P]=0$, 
i.e. $U=\alpha I$, $P=\beta I$ (Schur's lemma); then 
$A_j'M-MA_j''=(\beta -\alpha )G_j$ which means that $M=0$, $\alpha =\beta$, 
otherwise the sum of $L'$ and $L''$ must be direct). 

$3^0$. The condition dim Ext$^1(L',L'')\geq 2$ implies the existence of  
infinitesimal conjugations of $A_j'''$ by matrices $I+\varepsilon X_j$, 
$X_j=\left( \begin{array}{cc}Y_j&0\\Z_j&T_j\end{array}\right)$ which are not 
tantamount to a simultaneous infinitesimal conjugation of the $(p+1)$-tuple 
of matrices $A_j'''$. Indeed, one has 

\[ \tilde{A}_j:=(I+\varepsilon X_j)^{-1}A_j'''(I+\varepsilon X_j)=
A_j'''+\varepsilon 
\left( \begin{array}{cc}[A_j',Y_j]+G_jZ_j&0\\
A_j''Z_j-Z_jA_j'&[A_j'',T_j]-Z_jG_j\end{array}\right) +o(\varepsilon )\]
Set $\Phi =\{ (P_1,\ldots ,P_{p+1})|P_j=A_j''Z_j-Z_jA_j', 
\sum _{j=1}^{p+1}P_j=0\}$, 
$\Psi =\{ (Q_1,\ldots ,Q_3)|Q_j=A_j''Z-ZA_j'\}$ (as 
$\sum _{j=1}^{p+1}A_j^{(i)}=0$, 
$i=1,2$, one has $\sum _{j=1}^{p+1}Q_j=0$). The condition 
dim Ext$^1(L',L'')\geq 2$ implies that dim$\Phi -$dim$\Psi \geq 2$. 

$4^0$. Denote by $\Xi$ the codimension 1 subspace of $\Phi$ satisfying the 
condition tr$(\sum _{j=1}^{p+1}G_jZ_j)=0$. Hence, if the $(p+1)$-tuple of 
blocks $Z_j$ of the matrices $X_j$ belongs to $\Xi /\Psi$, then one can find blocks 
$Y_j$ and $T_j$ such that 
$\sum _{j=1}^{p+1}\tilde{A}_j=o(\varepsilon )$ (i.e. one can solve 
the equations $\sum _{j=1}^{p+1}[A_j',Y_j]+G_jZ_j=0$ and 
$\sum _{j=1}^{p+1}([A_j'',T_j]-Z_jG_j)=0$ w.r.t. $Y_j$ and $T_j$). 

$5^0$. Hence, one can find 
a conjugation of $A_j'''$ by matrices 
$I+\varepsilon X_j+\varepsilon ^2S_j(\varepsilon )$ analytic in 
$\varepsilon$, with $X_j$ as above and such that the sum of the conjugated 
matrices (denoted by $A_j$) to be 0 (identically in $\varepsilon$). The 
existence of $S_j$ 
is proved by complete analogy with the basic technical tool. 

As $(X_1,\ldots ,X_{p+1})\in (\Xi /\Psi )$, the $(p+1)$-tuple is irreducible. 
(Indeed, 
the matrix algebra ${\cal A}$ generated by the matrices $\tilde{A}_j$ 
contains matrices of the form $Y'=Y+O(\varepsilon )$ for every matrix $Y$ 
blocked as $A_j'''$; this follows from the main result of \cite{Ko4}. In 
particular, for $Y=H=\left( \begin{array}{cc}I&0\\
0&0\end{array}\right)$. 
Conjugate $Y'$ (with $Y=H$) by a matrix $I+O(\varepsilon )$ to make its 
left lower block 0 (identically in $\varepsilon$). 
This conjugation can't annihilate the left lower blocks of all matrices 
from ${\cal A}$ due to $(X_1,\ldots ,X_{p+1})\in (\Xi /\Psi )$. If 
$S\in {\cal A}$ 
has a non-zero such block, then $SH-HSH$ has the same left lower block and 
all its entries are $O(\varepsilon )$. Multiplying $S$ by matrices of 
the form $Y'$ and 
adding the matrices $Y'$, one obtains a basis of $gl(n,{\bf C})$.)
 
This proves 1). 

$6^0$. Prove 2) One can assume that the matrices $B'$, $B''$ have no 
non-zero eigenvalue in common. This can be achieved by multiplying one of 
the $(p+1)$-tuples of matrices $A_j^{(i)}$ by $h\in {\bf C}^*$ (we use here 
the fact that the matrices are nilpotent -- such a multiplication does not 
change the conjugacy classes). Hence, the matrix 
$B'''$ has at least $m'+m''$ simple non-zero eigenvalues. For 
$\varepsilon \neq 0$ small enough the eigenvalues of the matrix 
$B=\alpha _1A_1+\ldots +\alpha _{p+1}A_{p+1}$ will be close to the 
ones of $B'''$, hence, it will have at least $m'+m''$ non-zero simple 
eigenvalues. 

The lemma is proved.

\subsection{Some examples\protect\label{someexamples}}

We explain two methods for constructing irreducible triples or quadruples 
of nilpotent matrices. The first 
one places the non-zero entries of 
the matrices $A_j$ in positions $(k,k+1)$, $k=1,\ldots ,n-1$, and $(n,1)$. 
In all examples except (ex0) one has $p=2$.  

{\bf Example (ex0):} Let $n=2$, $p=3$. Set $A_1=-A_2=\left( \begin{array}{cc}
0&1\\0&0\end{array}\right)$, $A_3=-A_4=\left( \begin{array}{cc}
0&0\\1&0\end{array}\right)$. The quadruple is irreducible and for almost 
all values of $\alpha _j$ the matrix $B$ has two different non-zero 
eigenvalues. 

{\bf Example (ex1):} Let $n\geq 4$. Let 

\[ A_1=\left( \begin{array}{cccccccc}
0&1&0&0&0&\ldots &0&0\\0&0&1&0&0&\ldots &0&0\\0&0&0&1&0&\ldots &0&0\\
0&0&0&0&1&\ldots &0&0\\0&0&0&0&0&\ldots &0&0\\
\vdots &\vdots &\vdots &\vdots&\vdots &\ddots &\vdots &\vdots \\
0&0&0&0&0&\ldots &0&1\\
0&0&0&0&0&\ldots &0&0\end{array}\right) ~,~ 
A_2=\left( \begin{array}{cccccccc}
0&-1&0&0&0&\ldots &0&0\\0&0&0&0&0&\ldots &0&0\\0&0&0&-1&0&\ldots &0&0\\
0&0&0&0&-1&\ldots &0&0\\0&0&0&0&0&\ldots &0&0\\
\vdots &\vdots &\vdots &\vdots&\vdots &\ddots &\vdots &\vdots \\
0&0&0&0&0&\ldots &0&-1\\-1&0&0&0&0&\ldots &0&0\end{array}\right) ~,\]  
$A_3=-A_1-A_2$. One checks directly that the matrices are nilpotent, that 
$r_1=r_2=n-1$, $r_3=2$ and $(A_3)^2=0$. 

{\bf Example (ex2):} Let $n=3$. Let $A_1=\left( \begin{array}{ccc}0&1&0\\
0&0&1\\0&0&0\end{array}\right)$, $A_2=\left( \begin{array}{ccc}0&-1&0\\
0&0&0\\1&0&0\end{array}\right)$, $A_3=\left( \begin{array}{ccc}0&0&0\\
0&0&-1\\-1&0&0\end{array}\right)$. Hence, each matrix is nilpotent, of rank 
2. 

\begin{lm}\label{firstmethod}
The characteristic polynomial of a matrix having non-zero entries in 
positions $(k,k+1)$, $k=1,\ldots ,n-1$, and $(n,1)$ and zeros elsewhere is 
of the form $\lambda ^n+a$, $a\neq 0$. Hence, the eigenvalues of such a 
matrix are all non-zero and distinct.
\end{lm}

The lemma is to be checked directly. It implies that the eigenvalues of the 
matrices $B$ defined after examples (ex1) and (ex2) are non-zero and 
distinct. 

Another method is the non-zero entries to be in positions 
$(k,k+1)$, $k=1,\ldots ,n-1$, and $(n-1,1)$, $(n,2)$. The idea to give such 
examples -- the first of the matrices is initially in Jordan normal form and 
then one conjugates it with $I+E_{n,1}$. 
As our matrices will be sparce, we'll list only the entries of these 
positions, in this order. E.g., we write 
$A_1~:~1~1~0~|~1~-1$ instead of 
$A_1=\left( \begin{array}{cccc}0&1&0&0\\0&0&1&0\\1&0&0&0\\0&-1&0&0
\end{array}\right)$ (the vertical line separates the last two entries just 
for convenience). 
The reader is advised to draw in each example the matrices oneself. 

{\bf Example (ex3):} Let $n=6$. Let 

\[ \begin{array}{rrrrrrrrr}A_1~:&1&1&1&1&1&|&1&-1\\
A_2~:&-1&-1&0&-1&-1&|&0&0\\
A_3~:&0&0&-1&0&0&|&-1&1
\end{array}\]
Hence, $J(A_1)$ consists of a single block of size 6, $J(A_2)$ consists 
of two blocks of size 3 and $J(A_3)$ consists 
of three blocks of size 2. One can give examples 
when $n$ is even, $J(A_2)$ consists of two blocks of 
size 3 and of $(n-6)/2$ blocks of size 2; $J(A_3)$ consists of $n/2$ 
blocks of size 2; $J(A_1)$ consists of a single 
block of size $n$. To this end one adds to the right of the 
fourth from the left column of numbers in the above example a pack of 
$n-6$ units in the row of $A_1$, a pack of $(n-6)/2$ groups of the form 
$0, -1$ in the row of $A_2$ and a pack of $(n-6)/2$ groups of the form 
$-1, 0$ in the row of $A_3$.

{\bf Example (ex4):} Let $n=9$. Let 

\[ \begin{array}{rrrrrrrrrrrr}A_1~:&1&1&1&0&1&1&1&1&|&1&-1\\
A_2~:&-1&-1&0&-1&-1&0&-1&-1&|&0&0\\A_3~:&0&0&-1&1&0&-1&0&0&|&-1&1
\end{array}\] 
One checks directly that $J(A_1)$ consists of two Jordan blocks, of sizes 4 
and 5, $J(A_2)$ consists of three Jordan blocks of size 3 and $J(A_3)$ 
consists of three Jordan blocks of size 2 and of one of size 3. 

{\bf Example (ex5):} Let $n=10$. Let 

\[ \begin{array}{rrrrrrrrrrrrr}A_1~:&1&1&1&1&0&1&1&1&1&|&1&-1\\
A_2~:&-1&-1&0&-1&-1&-1&0&-1&-1&|&0&0\\A_3~:&0&0&-1&0&1&0&-1&0&0&|&-1&1
\end{array}\]
In this example $J(A_1)$ consists of two Jordan blocks of size 5, $J(A_2)$ 
consists of two Jordan blocks of size 3 and of one of size 4, $J(A_3)$ 
consists of five Jordan blocks of size 2. 

{\bf Example (ex6):} Let $n=12$. Let 

\[ \begin{array}{rrrrrrrrrrrrrrr}A_1~:&1&1&1&0&1&1&1&0&1&1&1&|&1&-1\\
A_2~:&-1&-1&0&-1&-1&0&-1&-1&0&-1&-1&|&0&0\\
A_3~:&0&0&-1&1&0&-1&0&1&-1&0&0&|&-1&1
\end{array}\]
Hence, $J(A_1)$ consists of three Jordan blocks of size 4, $J(A_2)$ consists 
of four Jordan blocks of size 3 and $J(A_3)$ consists of two Jordan blocks of 
size 3 and of three of size 2. 

{\bf Example (ex7):} Let $n=5$. Let 

\[ \begin{array}{rrrrrrrr}A_1~:&1&1&1&1&|&1&-1\\
A_2~:&-1&-1&0&-1&|&0&0\\
A_3~:&0&0&-1&0&|&-1&1
\end{array}\]
$J(A_1)$ consists of a single block of size 5, $J(A_2)$ and $J(A_3)$ consist 
each of a block of size 3 and of a block of size 2. One can give examples 
when $n>5$ is odd, each of $J(A_2)$ and $J(A_3)$ consists of one block of 
size 3 and of $(n-3)/2$ blocks of size 2; $J(A_1)$ consists of a single 
block of size $n$. To this end one adds to the left of the 
vertical line in the above example a pack of $n-5$ units 
in the row of $A_1$, a pack of $(n-5)/2$ groups of the form 
$0, -1$ in the row of $A_2$ and a pack of $(n-5)/2$ groups of the form 
$-1, 0$ in the row of $A_3$.

\begin{lm}\label{secondmethod}
The characteristic polynomial of a matrix having non-zero entries only in 
positions $(k,k+1)$, $k=1,\ldots ,n-1$, and $(n-1,1)$, $(n,2)$ is of the 
form $\lambda ^n+b\lambda$. Hence, when $b\neq 0$ its roots are all distinct 
and one of them equals 0.
\end{lm}

The lemma is to be checked directly. 

\begin{lm}\label{secondmethod1}
For the triples of matrices $A_j$ from each of the examples (ex3) -- (ex7) 
there exist conjugations of the matrices $A_j$ with matrices 
$I+O(\varepsilon )$ analytic in $\varepsilon \in ({\bf C},0)$ such 
that for $\varepsilon \neq 0$ the eigenvalues of the matrix $B$ are distinct 
and non-zero. 
\end{lm}

{\em Proof:}

$1^0$. In examples (ex4) -- (ex6) an 
infinitesimal conjugation of $A_1$ with $I+\varepsilon (E_{2,1}-E_{n,n-1})$, 
of $A_2$ with $I+\varepsilon E_{2,1}$ and of $A_3$ with 
$I-\varepsilon E_{n,n-1}$ creates a new triple (of matrices $A_j^1$) 
satisfying (\ref{A_j}) in first approximation w.r.t. $\varepsilon$. Hence, 
there exist true conjugations differing from the above ones by terms 
$O(\varepsilon ^2)$ after which for the deformed matrices $A_j^2$ 
(\ref{A_j}) holds. One has $A_j^1-A_j^2=o(\varepsilon )$. 

$2^0$. Choose $\alpha _j$ such that the entries of $B^1|_{\varepsilon =0}$ 
in positions $(k,k+1)$, $k=1,\ldots ,n-1$, and $(n-1,1)$, $(n,2)$ to be 
$\neq 0$. The entries in the following positions (and only they) of the 
matrix $B^1$ (hence, $B^2$ as well) are non-zero and $O(\varepsilon )$: 
$(1,1)$, $(2,2)$, $(n-1,n-1)$, $(n,n)$ and $(n,1)$. Hence, 
det$B^1=v\varepsilon +o(\varepsilon )$, $v\neq 0$ and for 
$\varepsilon \neq 0$ small enough all eigenvalues of $B^1$ and $B^2$ are 
non-zero and distinct, the last of them is $O(\varepsilon )$. 

$3^0$. In examples (ex7) and (ex3) the triple of matrices $A_j^1$ is created 
by an infinitesimal conjugation of $A_1$ and $A_2$ with 
$I-\varepsilon E_{n,n-1}$ and of $A_3$ with $I+\varepsilon E_{2,1}$. 
The entries in the following positions (and only they) of the 
matrix $B^1$ (hence, $B^2$ as well) are non-zero and $O(\varepsilon )$: 
$(n-1,n-1)$, $(n,n)$ and $(n,1)$. For the rest the reasoning is the same. 

The lemma is proved. 

\begin{lm}\label{Extexamples}
1) Denote by $\Phi _1$, $\Phi _2$ two irreducible representations defined by 
any of examples (ex1) -- (ex7) (with $n=6$ in (ex3) and $n=5$ in (ex7)). 
Then Ext$^1(\Phi _1,\Phi _2)\geq 2$.

2) The same is true if each of $\Phi _1$, $\Phi _2$ corresponds to one 
of the examples (ex3) or (ex7) without the restrictions 
respectively $n=6$ and $n=5$.  
\end{lm}

The proof is left for the reader.

\subsection{Proof of Theorem \protect\ref{threematrices}
\protect\label{proofofthreematrices}}

\subsubsection{Simplification and plan of the proof\protect\label{plan}}

{\bf 1.} The first observation to be made is that if there exists a 
$(p+1)$-tuple of matrices $A_j\in c_j$ satisfying 1) or 2) of the theorem, 
then there exist such $(p+1)$-tuples for every $(p+1)$-tuple of conjugacy 
classes $c_j'$ where for each $j$ either $c_j=c_j'$ or $c_j$ is 
subordinate to $c_j'$. To this end one has to apply the basic technical tool.

{\bf 2.} Call {\em operation} $(s,l)$, $s\geq l\geq 1$ the changing of a 
given nilpotent conjugacy class $c$ containing two Jordan blocks of sizes 
$s$ and $l$ to $c'$ in which these blocks are replaced by two blocks of 
sizes $s+1$ and $l-1$. The class $c$ lies in the closure of $c'$. 
If the nilpotent conjugacy class $c$ lies in the 
closure of the nilpotent conjugacy class $c''$, then there exists a 
sequence of conjugacy classes $c_0=c$, $c_1$, $\ldots$ ,$c_{\mu}=c''$ such 
that $c_i$ is obtained from $c_{i-1}$ as a result of some operation $(s,l)$. 
This is proved in \cite{Ko1}. 

{\bf 3.} Let $A\in c$, $A'\in c'$, $c$, $c'$ being nilpotent conjugacy 
classes. Then $c$ lies in the closure of $c'$ if and only if for all $i$ 
one has rk$(A)^i\leq$rk$(A')^i$ (proved in \cite{Ko1}). 

{\bf 4.} For given $r\in {\bf N}$ denote by $\Omega _0(r)$ the nilpotent 
orbit of rank $r$ of least dimension. It is unique and has Jordan blocks 
either of one and the same size or of two consecutive sizes. This follows 
from {\bf 2}. If the size is one and the same, then 
there exists a single orbit 
$\Omega _1(r)$ containing $\Omega _0(r)$ in its closure and contained in 
the closure of any orbit with the same value of $r$ and different from 
$\Omega _0(r)$ and $\Omega _1(r)$. It is obtained from $\Omega _0(r)$ by 
an operation $(k,k)$, $k$ being the size of the blocks of $\Omega _0(r)$. 

{\bf 5.} Let $r_1+\ldots +r_{p+1}>2n$. Then it is possible to change some 
of the conjugacy classes to subordinate ones of smaller rank to get the 
condition $r_1+\ldots +r_{p+1}=2n$ and the new $(p+1)$-tuple of ranks not to 
correspond to any of the ones of special or almost special cases.

{\bf 6.} Call {\em merging} of two nilpotent conjugacy classes $c$ and $c'$ 
the following procedure defined for $r(c)+r(c')\leq n-1$ and when 
$c=\Omega _0(r(c))$, $c'=\Omega _0(r(c'))$ (hence, at least one of the 
two Jordan normal forms has only blocks of size $\leq 2$). If $A\in c$, 
$A'\in c'$ are Jordan matrices with decreasing order of the sizes $b_i$, 
$b_i'$ of the blocks, rk$A\geq$rk$A'$ (hence, $b_i'\leq 2$),  
then construct the nilpotent matrix $A''$ as follows: insert on the first 
superdiagonal (it comprises the positions $(k,k+1)$) between 
the packs of $b_i-1$ and $b_{i+1}-1$ units from $A$ a unit from 
$A'$ as long as this is possible. These packs are units 
from Jordan blocks of sizes $b_i$ and $b_{i+1}$. 
When inserting a  
unit, it takes the place of the 0 and the units of $A$ do not change their 
positions. 

Given the matrix $A''$ it is self-evident how to represent it in the form 
$A+A'$ with $A\in c$, $A'\in c'$. If there exists an irreducible or nice 
representation in which one of the matrices is from $c''$, then there exists 
such a representation with one more matrix, two of the matrices being 
from $c$ and $c'$, the other conjugacy classes remaining the same.  

{\bf 7.} We prove the theorem only in the case $r_1+\ldots +r_{p+1}=2n$ 
making use of {\bf 5.} Making use of {\bf 4.}, when $p\geq 4$, one might 
restrict oneself to the case when $c_j=\Omega _0(r_j)$ for all $j$. 

It is also possible to restrict oneself to the cases 
$p=2$ and $p=3$ due to the possibility to merge Jordan normal forms. 
Indeed, if $r_1+\ldots +r_{p+1}=2n$, if $c_j=\Omega _0(r(c_j))$ for all $j$ 
and if $p\geq 4$, then a merging is 
always possible. Moreover, it is possible to be done with avoiding to come 
to the special case a) or to the almost special case a1). Indeed, when 
passing from $p=4$ to $p=3$ by a merging, the sum of three of the quantities 
$r_j$ is $\leq n$ and one can choose a couple to be merged such that one 
of the four quantities $r_j$ which remain after the merging to be $<n/2$. 

When passing from $p=3$ to $p=2$ by merging (the passage is not defined for 
the cases a) and a1)), one can avoid to come to any of 
the other special or almost special cases. Indeed, merging results either in 
a Jordan normal form with a single Jordan block of size $m>1$ or with 
greatest difference $h$ between the sizes of two of the Jordan blocks 
$\geq 2$; one has $m\geq 3$ and the cases $m=3$ and $h=2$ are possible only 
if one of the Jordan normal forms to be 
merged has only one block of size $2$, the rest of size 1. In this case 
one can merge other two of the Jordan normal forms to avoid the special and 
almost special cases.

{\bf 8.} Finally, after having restricted oneself to the cases $p=2$ and 
$p=3$, one can again use {\bf 1.}, {\bf 2.} and {\bf 4.} to replace the 
given conjugacy classes by the triple or quadruple of conjugacy classes 
$\Omega _0(r_j)$. This is not always possible because one could come to 
a special or almost special case. 

Let this be not a special or an almost special case. 
Use the notation "$A_1=2$" or "$A_1=(2,3)$" in the sense "the Jordan 
normal form of $A_1$ has only Jordan blocks of size 2" ("of sizes 2 and 3"). 
The following cases of triples $\Omega _0(r_1)$, $\Omega _0(r_2)$, 
$\Omega _0(r_3)$ are possible (see \cite{Ko2}):

\[ \begin{array}{lcclc}(A)&A_1=(1,2)&&(B)&A_1=2~{\rm or}~A_1=(2,3); A_2=(2,3)\\ 
(C)&A_1=(2,3); A_2=3; A_3=(3,4)&&(D)&A_1=(2,3); A_2=3; A_3=4\\ 
(E)&A_1=(2,3); A_2=3; A_3=(4,5)&&(F)&A_1=(2,3); A_2=3; A_3=5\\ 
(G)&A_1=(2,3); A_2=3; A_3=(5,6)&&
(H)&A_1=2~{\rm or}~(2,3); A_2=(3,4); A_3=4\\ 
(I)&A_1=2~{\rm or}~(2,3); A_2=(3,4); A_3=(4,5)&&
(J)&A_1=2~{\rm or}~(2,3); A_2=(3,4); A_3=5\\
(K)&A_1=2~{\rm or}~(2,3); A_2=(3,4); A_3=(5,6)&&&\end{array}\] 

Case (C) is considered in Subsubsection \ref{caseC}, case (F) is 
considered in Subsubsection~\ref{caseF}. The other cases are considered 
in Subsubsection~\ref{caseA}. 

If the triple or quadruple of classes $\Omega _0(r_j)$ 
is a special case, then we use {\bf 2.} and {\bf 4.} and 
prove part 1) of the theorem for all almost special cases stemming from it -- 
this will imply that part 1) holds for the initial triple of nilpotent 
orbits, see Subsubsection~\ref{ASC2}. For the almost special cases 
c1), d1) and d2) we prove that 2) holds as well, see 
Subsubsection~\ref{ASC1}. 

Finally, we consider all cases obtained from one of the almost special 
ones a1), b1), c2) and d3) by an operation $(s,l)$ on one of the three 
orbits (call such cases {\em neighbouring}), see 
Subsubsection~\ref{NC}.  We prove that 
part 2) of the theorem holds for all these cases. Hence, part 2) will 
hold for all cases when $p=2$ or 3, when $r_1+\ldots +r_{p+1}=2n$ and when 
all special and almost special cases are avoided. Hence, it will hold 
when $r_1+\ldots +r_{p+1}\geq 2n$ and when 
all special and almost special cases are avoided.

\subsubsection{Proof of the theorem for $p=2$, in case (C)
\protect\label{caseC}}

We construct triples of matrices satisfying 
the requirements of the lemma. The construction is done by induction on $n$. 
The induction base are the examples from Subsection~\ref{someexamples}. 

The size $n$ must be 
$>3$. Decrease the size of the matrices by 3 and delete a block of size 3 
from each Jordan normal form. Denote the new Jordan normal forms by $J_j'$. 
This can lead only to case (C) again or to 
case (D) (with $n$ replaced by $n-3$). By inductive assumption, there exists 
an irreducible triple of nilpotent matrices $A_j'$ of size $n-3$ with 
$J(A_j')=J_j'$ and such that the matrix 
$B'=\alpha _1A_1'+\alpha _2A_2'+\alpha _3A_3'$ have a simple non-zero 
eigenvalue.

Denote by $A_j''$ the triple of matrices $A_j$ from example (ex2). 
One can assume that the matrices $B'$ and 
$B''=\alpha _1A_1''+\alpha _2A_2''+\alpha _3A_3''$ have no eigenvalue in 
common (to achieve this one can, if necessary, multiply one of the triples  
by $h\in {\bf C}^*$). By Lemma~\ref{Ext>=2} (the reader should check that 
it is applicable), there exist irreducible triples 
of nilpotent matrices $A_j$ with $J(A_j)=J_j$, the matrix $B$ having at 
least two simple non-zero eigenvalues (hence, at least one). This proves the 
induction step in case (C).

\subsubsection{Proof of the theorem for $p=2$, in cases 
(A), (B), (D), (E) and (G) -- (K)\protect\label{caseA}}

$1^0$. The cases (B), (D), (E), (G) and (H) -- (K) are considered by 
analogy with (C) and we define only the matrices $A_j''$. As we make use 
of examples (ex3) -- (ex7), one should keep in mind 
Lemma~\ref{secondmethod1}. Lemma~\ref{Ext>=2} is applicable in 
all cases because there holds Lemma~\ref{Extexamples}. (One can represent 
the triple of Jordan normal forms as a direct sum of triples of 
Jordan normal forms corresponding to one of the examples (ex1) -- (ex7) 
and then use Lemma~\ref{Extexamples}.)

In case (B) $J(A_3)$ must contain a block of size $\kappa \geq 5$. The 
matrices $A_j''$ of size $\kappa$ are defined by examples (ex7) and (ex3) 
depending on whether 
$n$ is even or odd. If $J(A_3)$ consists of a single block, then these 
examples define directly the matrices $A_j$. If not, then the matrices 
$A_j'$ correspond again to case (B). 

In case (D) the matrices $A_j''$ are of size 12; they are  
the triple from example (ex6) (if $n=12$, then 
case (D) is proved directly by example (ex6)). The matrices $A_j'$ are from 
case (D) again.

In case (E) the matrices $A_j''$ are of size 9; they are defined in example 
(ex4). The matrices $A_j'$ are from case (E) or from case (D) or from case 
(F). 

In case (G) there are at least three blocks of size 2 in $J(A_1)$. The 
matrices $A_j''$ are defined by example (ex3). The matrices $A_j'$ are either 
from case (G) or from case (F). 

In cases (H) and (I) the matrices $A_j''$ are of size 4; they are defined by 
example (ex1) with $n=4$. One has to prove that there are at least two 
Jordan blocks 
of size 2 in $J(A_1)$, see \cite{Ko2}. In case (H) the matrices $A_j'$ are 
from case (H) or from case (D), in case (I) they are from case (I), (J), (H), 
(E), (F) or (D). 

In case (J) the matrices $A_j''$ are of size 10; they are defined by 
example (ex5) ($J(A_2)$ contains at least five blocks of size 3, 
otherwise $r_1+r_2+r_3>2n$). The matrices $A_j'$ are either from case (J) 
or from case (F). 

In case (K) there are at least three blocks of size 2 in $J(A_1)$ and at 
least two blocks of size 3 in $J(A_2)$ for the same reason, see \cite{Ko2}, 
and the matrices $A_j''$ are of size 6, defined by example 
(ex3). The matrices $A_j'$ are from case (K), (J), (G) or (F).  

In all these cases the triple of Jordan normal forms of the matrices $A_j'$ 
corresponds to neither of the special or almost special cases. 

$2^0$. Consider case (A) (we follow the same ideas as 
in \cite{Ko2}). Decrease by 1 the sizes of the matrices and 
decrease by 1 the sizes of two Jordan blocks respectively of $J(A_2)$ and 
$J(A_3)$. Delete a Jordan block of size 1 from $J(A_1)$. One can choose the 
diminished blocks of $J(A_2)$, $J(A_3)$ such that the triple of Jordan 
normal forms of size $n-1$ obtained like this not to correspond to any of 
the special or almost special cases. (It suffices to leave in each of the 
two Jordan normal forms of size $n-1$ a couple of blocks of different size.) 
This defines the Jordan normal forms of the 
matrices $A_j'$. The sum of their ranks equals $2(n-1)$. 

Set $A_j^1=\left( \begin{array}{cc}A_j'&G_j\\0&A_j''\end{array}\right)$. 
The matrices $A_j''$ are of size 1 and equal 0. One constructs the blocks 
$G_j$ such that the Jordan normal form of $A_2^1$ to be 
the necessary one (i.e. $\Omega _0(r_2)$). One sets $G_3=0$, $G_2=-G_1$. 
(The condition $G_3=0$ can be achieved by conjugating the triple.) 
Hence, the matrices $A_j^1$ {\em do not} define a semi-direct 
sum of $L'$ and $L''$ because $J(A_2^1)$ is not a direct sum of 
$J(A_2')$ and $J(A_2'')$. 

Thus the matrix $B^1$ has a simple non-zero eigenvalue. It is an eigenvalue 
of $B'$.

After this one deforms the triple of matrices $A_j^1$ into an 
irreducible one, with the necessary Jordan normal forms, see \cite{Ko2}. 
One can also use the same reasoning as the one from case (F), see the next 
subsubsection. 
Hence, for the deformed triple the analog of the matrix $B^1$ still 
has a simple non-zero eigenvalue.

\subsubsection{Proof of the theorem for $p=2$, in case (F)
\protect\label{caseF}}

There exists an irreducible triple of nilpotent $9\times 9$ matrices $A_1'$, $A_2'$, 
$A_3'$, $A_1'+A_2'+A_3'=0$ where $J(A_1')$ ($J(A_2')$; $J(A_3')$) contains 
1 block $3\times 3$ and 3 blocks $2\times 2$ (3 blocks $3\times 3$; 
1 block $4\times 4$ and 1 block $5\times 5$), see case (E). Hence, the 
matrix $B'$ has a non-zero simple eigenvalue. 

There exists a triple of 
nilpotent $15\times 15$-matrices $A_1^0$, $A_2^0$, 
$A_3^0$, $A_1^0+A_2^0+A_3^0=0$ with $J(A_j^0)={\cal J}_j$ where 
${\cal J}_1$ (${\cal J}_2$; ${\cal J}_3$) consists of 1 block $3\times 3$ 
and 6 blocks $2\times 2$ (of 5 blocks $3\times 3$; of 3 blocks $5\times 5$) 
and the matrices $A_1^0$, $A_2^0$ look like this (recall that $A_j'$ are 
$9\times 9$):

\[ \left( \begin{array}{ccccccc}A_1'&\eta _1&\eta _2&\eta _3&\eta _4&
\eta _5&\eta _6\\0&0&1&0&0&0&0\\0&0&0&0&0&0&0\\0&0&0&0&1&0&0\\
0&0&0&0&0&0&0\\0&0&0&0&0&0&0\\0&0&0&0&0&0&0\end{array}\right) ~~,~~
\left( \begin{array}{ccccccc}A_2'&\varphi _1&\varphi _2&\varphi _3&
\varphi _4&
\varphi _5&\varphi _6\\0&0&-1&0&0&0&0\\0&0&0&-1&0&0&0\\0&0&0&0&0&0&0\\
0&0&0&0&0&-1&0\\0&0&0&0&0&0&-1\\0&0&0&0&0&0&0\end{array}\right)\]

The matrix $B^0$ has the same non-zero eigenvalues as $B'$. Hence, 
it has a simple non-zero eigenvalue. One can deform the triple of matrices 
$A_j^0$ into an irreducible triple of matrices $A_j$ with the necessary 
Jordan normal forms, see \cite{Ko2}. Hence, the matrix $B$ for small values 
of the deformation parameter still has a simple non-zero eigenvalue. 

We explain the details of this deformation to show why this method is not 
applicable to some of the almost special cases (as was claimed in 
\cite{Ko2}). 

One looks for matrices $A_j$ of the form 
$A_1=A_1^0+\varepsilon L$ where only the last row of $L$ is non-zero, 
with $L_{15,14}=1$; we assume that $L_{15,j}=0$ for $j=10,11,12,13,15$ and 
that for $\varepsilon \neq 0$ the sizes of the Jordan blocks of $J(A_1)$ 
equal 3,2,2,2,2,2,2; 
one sets $A_j=(I+\varepsilon X_j(\varepsilon ))^{-1}A_j^*
(I+\varepsilon X_j(\varepsilon )$, $j=2,3$, where $X_j$ are analytic in 
$\varepsilon \in ({\bf C},0)$ (their existence is justified by the 
basic technical tool). 
 
The matrix algebra ${\cal A}^0$ generated by the matrices $A_j$ contains 
matrices of the 
form $\left( \begin{array}{cc}P&Q\\O(\varepsilon )&O(\varepsilon )
\end{array}\right)$ for all $P$, $Q$ where $P\in gl(9,{\bf C})$. The 
centralizer of ${\cal A}^0$ is trivial, see \cite{Ko2}. 

Fix a matrix $F\in {\cal A}^0$ with $P=I$, $Q=0$. By taking powers of $F$, 
one can assume that the 6 right columns of $F$ are identically 0. 

Conjugate the matrices $A_j$ by a matrix 
$\left( \begin{array}{cc}I&0\\O(\varepsilon )&I\end{array}\right)$ so that 
after the conjugation the matrix $F$ become equal to 
$\left( \begin{array}{cc}I&0\\0&0\end{array}\right)$. Hence, 
${\cal A}^0$ contains all matrices of the form 
$Y=\left( \begin{array}{cc}P&Q\\0&0\end{array}\right)$ (if $S\in {\cal A}^0$, 
then $Y=FS\in {\cal A}^0$). Two cases are possible after the conjugation:

1) the three matrices become block upper-triangular, with diagonal 
blocks of sizes 9 and 6;

2) there exists an entry in the left lower block 
$6\times 9$ which is $\neq 0$ for $\varepsilon \neq 0$. 

Eliminate case 1). If this were true, then the restrictions of $A_j$ to the 
right lower block $6\times 6$ would be nilpotent, with sizes of the blocks 
equal to 2,2,2; 3,3; 5,1. They must define a nilpotent algebra ${\cal A}^1$. 

Indeed, it 
is impossible to have an irreducible triple with such sizes of the blocks 
because one would have $r_1+r_2+r_3=11<12$. 
It is impossible to conjugate the algebra to a 
block upper-triangular form with at least one diagonal block $\tilde{P}$ 
irreducible, 
of size $m>1$. Indeed, the restrictions of the matrices to such a block 
would be nilpotent, with sizes of the blocks not greater respectively 
than 2, 3 and 5. One checks directly that for $m=2,3,4$ and 5 it is 
impossible to have $r_1'+r_2'+r_3'\geq 2m$, $r_j'$ being the ranks of the 
restrictions of the matrices to the block $\tilde{P}$. (It is this part of 
the reasoning which 
is not applicable to the proof of the almost special cases a1), b1), c2) and 
d3); e.g., in case b1) there exist block upper-triangular triples 
$6\times 6$ with irreducible diagonal blocks $3\times 3$.)

On the other hand if ${\cal A}^1$ is nilpotent, then it can be conjugated to 
an upper-triangular form. The restriction to the right lower block 
$6\times 6$ of the matrix $(A_1+A_2/2)|_{\varepsilon =0}$ has non-zero 
entries in positions $(k,k+1)$, $k=10,\ldots ,14$. Hence, the conjugation 
can be carried out by a matrix $I+O(\varepsilon )$. Such a conjugation 
cannot annihilate the entry $A_{1;15,14}$. Hence, the algebra ${\cal A}^1$ 
is not nilpotent. 

Consider case 2). Let the algebra ${\cal A}^0$ contain a matrix $S$ with 
$S_{10,j}|_{\varepsilon =0}\neq 0$ for some $j\leq 9$ (if one has 
$S_{i,j}|_{\varepsilon =0}\neq 0$ for $i>10$ and for some $j\leq 9$, then 
one can multiply $S$ by $(A_1+A_2/2)^{i-10}$ to have 
$S_{10,j}|_{\varepsilon =0}\neq 0$). 

Then one can assume that only the 
left lower block $6\times 9$ of $S$ is non-zero (one can consider 
instead of $S$ the matrix $SF-FSF$). Multiply the matrix $S$ by matrices $Y$ 
defined above. Hence, ${\cal A}^0$ contains matrices of the form 
$\left( \begin{array}{cc}P&Q\\O(\varepsilon )&O(\varepsilon )
\end{array}\right)$ with $P\in gl(10,{\bf C})$, for all $P$ and $Q$ and 
one can repeat the reasoning which led us to cases 1) and 2), but this time 
the size of the block $P$ has increased by 1. Continuing like this, we see 
that ${\cal A}^0=gl(15,{\bf C})$, i.e. the triple $A_1$, $A_2$, $A_3$ is 
irreducible.

\subsubsection{Proof of the theorem in the almost special cases c1), d1) 
and d2)\protect\label{ASC1}}

There exist irreducible triples of nilpotent matrices $A_j$ 
satisfying (\ref{A_j}) with sizes of the Jordan blocks like in cases 
c1), d1) or d2) but with $g=1$. These triples can be obtained by 
deforming respectively the ones from examples (ex1) with $n=4$ for case 
c1) and (ex3) for cases d1) and d2). 

Consider the direct sum of such a triple and of a triple from example 
(ex1) with $n=4$ for case c1) and of one from example (ex3) for cases 
d1) and d2). Lemma~\ref{Ext>=2} is applicable to such a direct sum which 
provides the existence of irreducible triples from cases c1), d1) and d2) 
for $g=2$. In the same way one constructs such triples for all $g>1$ -- 
by deforming the direct sum of a triple for $g-1$ and of one from 
example (ex1) with $n=4$ for case 
c1) or (ex3) for cases d1) and d2) and by using Lemma~\ref{Ext>=2}. 

The irreducible representations thus obtained can be considered as 
deformations of certain direct or semi-direct sums of representations 
whose diagonal blocks are of sizes 4 or 6 and whose matrices $B$ have 
distinct non-zero eigenvalues, see Lemma~\ref{secondmethod1}. This property 
persists under small deformations.

\subsubsection{Proof of the theorem in the almost special cases a1), b1),  
c2) and d3)\protect\label{ASC2}}

We consider only case b1) in detail. The other cases are treated by analogy 
and we explain the differences at the end of the subsubsection. Recall that 
in all these cases we prove the existence of nice representations.

Construct the matrices 
$A_j=\left( \begin{array}{ccccc}A_j^1&0&\ldots&0&H_j^1\\
0&A_j^2&\ldots&0&H_j^2\\\vdots&\vdots&\ddots&\vdots&\vdots\\
0&0&\ldots&A_j^{g-1}&H_j^{g-1}\\0&0&\ldots&0&A_j^g\end{array}\right)$ where 
$A_j^k$ are $3\times 3$, nilpotent, of rank 2 and the representations 
defined by the triples $A_j^k$ are irreducible for all $k$. Moreover, 
they are presumed to be non-equivalent (this can be achieved by 
multiplying them by constants $g_k\in {\bf C}^*$) and the matrices $B^k$ to 
have non-zero distinct eigenvalues; the eigenvalues of the matrix $B$ can be 
presumed non-zero and distinct as well; see example (ex2). 

Assume that the matrix $A_1$ is in upper-triangular Jordan normal form. 
Then we set 
$H_1^k=\left( \begin{array}{ccc}0&0&0\\0&0&0\\0&0&1\end{array}\right)$ for 
all $k$. Hence, the Jordan normal form of $A_1$ consists of $g-2$ blocks of 
size 3, of a block of size 4 and of a block of size 2 (to be checked 
directly). 

The blocks $H_j^k$ for $j=2,3$ are defined such that 
$H_j^k=A_j^kD_j^k-D_j^kA_j^d$, $D_j^k\in gl(3,{\bf C})$ and 
$H_1^k+H_2^k+H_3^k=0$. Such a choice of $H_j^k$ is possible because the 
representations defined by the triples $A_j^k$ for different values of 
$k$ are irreducible non-equivalent and the mapping 

\[ (D_2^k,D_3^k)\mapsto A_2^kD_2^k-D_2^kA_2^g+A_3^kD_3^k-D_3^kA_3^g\] 
is surjective onto $gl(3,{\bf C})$. Notice that the blocks $H_2^k$, $H_3^k$ 
result from conjugation of $A_j$ with $I+D_j$ where $D_j$ has non-zero 
entries only in the last 3 columns and first $3(g-1)$ rows (its restriction 
to the rows with indices $3k-2$, $3k-1$, $3k$ and to the last three 
columns equals $D_j^k$). Hence, $J(A_2)$ and $J(A_3)$ 
consist each of $g$ blocks of size 3. 

Prove that the centralizer ${\cal Z}$ of the triple of matrices 
$A_j$ is trivial. Block-decompose a matrix from $gl(n,{\bf C})$ into blocks 
$3\times 3$ 
(denoted by $X^{\mu, \nu}$). For $Y\in {\cal Z}$ set 
$Y^{\mu ,\nu}=Y|_{X^{\mu ,\nu}}$. One has (first for $\nu <\mu =g$ and then 
for $\mu \geq \nu$, $\nu \leq g-1$) 
$Y^{\mu ,\nu}A_j^{\nu}-A_j^{\mu}Y^{\mu ,\nu}=0$. Hence, if $\mu \neq \nu <g$, 
then $Y^{\mu ,\nu}=0$, if $\mu =\nu$, then $Y^{\mu ,\nu}=\alpha _{\mu}I$ 
(we use the non-equivalence of the representations and Schur's lemma).

For $\nu =g$ one has 
$(\alpha _k-\alpha _g){H_1^k}+A_1^kY^{k,g}-Y^{k,g}A_1^g=0$. 
Hence, $\alpha _k=\alpha _g$ because $H_1^k$ is not of the form 
$A_1^kY^{k,g}-Y^{k,g}A_1^g$. But then 
$A_j^kY^{k,g}-Y^{k,g}A_j^g=0$ for $j=2,3$ which implies 
$Y^{k,g}=0$ for $k<g$. Hence, the centralizer is trivial. 

In all other almost special cases one similarly constructs triples 
or quadruples of matrices $A_j$ satisfying the conclusions from 1) of 
the theorem.

In all cases the matrix $A_j$ whose Jordan form has to be changed (from 
equal sizes of the Jordan blocks to one obtained by replacing a couple 
$l_j$, $l_j$ of sizes 
by $l_j+1$, $l_j-1$) has equal blocks $H_j^k$ which have a unit in the right 
lower corner and zeros elsewhere; the restriction of the matrix $A_j$ to 
the diagonal blocks is in Jordan normal form. The diagonal blocks are 
of sizes 2, 4 or 6. We leave the details for the reader.

\subsubsection{Proof of the theorem in the neighbouring  
cases\protect\label{NC}}

There are two types of neighbouring cases. Recall that an almost special 
case 
is obtained from a special one by an operation $(s,l)$ performed on one of 
the three or four Jordan normal forms and a neighbouring case is obtained by 
performing another such operation. In the first type the second operation is 
performed on one of the other Jordan normal forms, in the second type it is 
performed on the same one. We consider only cases neighbouring to 
a1), b1), c2) and d3), in the cases neighbouring to c1), d1) and d2) there 
exist irreducible representations, see Subsubsection~\ref{ASC1}.

{\bf Remark:} In all neighbouring cases the triples or quadruples can be 
considered as 
deformations of triples or quadruples from almost special cases, in which 
the matrix $B$ has distinct non-zero eigenvalues. Hence, this is so in 
all neighbouring cases as well.

{\bf Neighbouring cases of the first type.}

Construct an irreducible triple or quadruple in a neighbouring case 
of the first type. Let $A_j^*$ be the matrices from the triple or quadruple 
of the corresponding almost special case as constructed in the previous 
subsubsection. 

Consider the only case of first type 
neighbouring to b1). (It is the only one up to permutation of the three 
Jordan normal forms. The rest of neighbouring cases of the first type are 
considered by analogy.) 
Suppose that it is the orbit of $A_2^*$ to be changed and that $A_2^*$ is in 
upper-triangular Jordan normal form. We assume that the triple $A_1^*$, 
$A_2^*$, $A_3^*$ is obtained from the triple of matrices $A_j$ constructed 
in the previous subsubsection by conjugation with $(I+D_2)^{-1}$. 
Set $A_2=A_2^*+\varepsilon A_2^0$ where the matrix $A_2^0$ has non-zero 
entries only in the last three rows and first $n-3$ columns. The restriction 
of $A_2^0$ to the block $X^{g,k}$ defined in the previous subsubsection  
equals $H_1^k$ (defined also there). Hence, for 
$\varepsilon \neq 0$ the Jordan normal form of $A_2$ is the required one. 

After this look for $A_j$ ($j=1,3$) in the form 
$A_j=(I+\varepsilon X_j(\varepsilon ))^{-1}A_j^*
(I+\varepsilon X_j(\varepsilon )$ where $X_j$ are analytic in 
$\varepsilon \in ({\bf C},0)$ (their existence is justified by the 
basic technical tool). 

Show that the algebra ${\cal A}'$ generated by the matrices $A_j$ is 
$gl(n,{\bf C})$ 
from where part 2) of the theorem follows. The algebra ${\cal A}$ 
generated by the 
matrices $A_j^*$ is the one of all matrices $W$ having arbitrary entries in 
the diagonal blocks $X^{i,i}$ and in the blocks $X^{g,i}$. This follows from 
the basic result in \cite{Ko4}. Hence, ${\cal A}'$ contains matrices of the 
form $S+\varepsilon T$ for all $S\in {\cal A}$.

The algebra ${\cal A}'$ 
contains the matrix 
$(A_2)^3/\varepsilon =\sum _{\kappa =1}^{g-1}E_{3g-2,3\kappa }$.   
By multiplying and postmultiplying it by matrices from 
${\cal A}'$, one can obtain matrices of the form $V+\varepsilon Z$ for $V$ 
having any restriction to the block $X^{i,j}$, for all $i$ and for 
$j\leq g-1$. 

The matrices $W$ and $V$ contain a basis of $gl(n,{\bf C})$ for 
$\varepsilon \neq 0$ small enough. Hence, the triple $A_1$, $A_2$, $A_3$ 
is irreducible.  

{\bf Neighbouring cases of the second type.}

{\bf $1^0$.} The cases of second type neighbouring to a given 
almost special case can be characterized by the sizes of the blocks of 
the Jordan normal form which changes w.r.t. the corresponding special case. 
There are four possibilities:

\[ \begin{array}{lllll}
1)&l_j+2, l_j-2,l_j,\ldots ,l_j&&2)&l_j+2,l_j-1,l_j-1,l_j,\ldots ,l_j\\
3)&l_j+1,l_j+1,l_j-2,l_j,\ldots ,l_j&&4)&
l_j+1,l_j+1,l_j-1,l_j-1,l_j,\ldots ,l_j\end{array}\]

Possibilities 1), 2), 3) and 4)  
appear for the first time respectively for $g=2$, $g=3$, $g=3$ and $g=4$. 
We explain the construction for these minimal values of $g$, for all others 
the existence is proved by induction on $g$, when considering direct 
sums of triples or quadruples constructed for $g-1$ and triples or 
quadruples defined by examples (ex0), (ex1) with $n=4$ or (ex3). We deform 
such direct sums into irreducible representations by means of 
Lemma~\ref{Ext>=2}.

{\bf $2^0$. Possibility 1).}

Let $g=2$. Explain in details the case neighbouring to b1). 
Denote by 
$A_j^*=\left( \begin{array}{cc}A_j^1&H_j^1\\0&A_j^2\end{array}\right)$ 
matrices defining triples from case b1) with $g=2$ and $H_j^1$ defined like 
in the previous subsubsection. The centralizer of the triple is trivial and 
the representations defined by the matrices $A_j^1$ and $A_j^2$ are 
non-equivalent. $A_1^1$ and $A_1^2$ are upper-triangular Jordan blocks of 
size 3. Hence, the matrix algebra ${\cal A}$ 
generated by the matrices $A_j^*$ contains all block upper-triangular 
matrices with blocks $3\times 3$, see \cite{Ko4}. In particular, it contains 
the matrix $S=\left( \begin{array}{cc}0&0\\0&I\end{array}\right)$.  
%$T=\left( \begin{array}{cc}I&0\\0&0\end{array}\right)$.

Set $A_1=A_1^*+\varepsilon Y$, $Y=E_{4,1}$. Hence, $A_1$ has for 
$\varepsilon \neq 0$ Jordan blocks of sizes 5 and 1. Set 
$A_j=(I+\varepsilon X_j(\varepsilon ))^{-1}A_j^*
(I+\varepsilon X_j(\varepsilon ))$, $X_j$ being analytic in 
$\varepsilon \in ({\bf C},0)$ 
and such that $\sum _jA_j=0$ (the existence of $X_j$ follows from the 
basic technical tool). 

The algebra ${\cal A}'$ contains for all $F\in {\cal A}$ a matrix 
$F+O(\varepsilon )$. Hence, ${\cal A}'$ contains a matrix 
$S'=S+O(\varepsilon )$. 
Hence, it contains the matrix $Q=S'(A_1)^3$ which is of the form 
$\varepsilon \left( \begin{array}{cc}\ast &\ast \\G&\ast \end{array}\right) 
+o(\varepsilon )$ 
with $G=E_{4,3}\neq 0$. By multiplying and postmultiplying $Q/\varepsilon$ 
by matrices from ${\cal A}'$, one can obtain matrices of the form of 
$Q/\varepsilon$ with any block $G$. These matrices together with the 
matrices $F+O(\varepsilon )$ form a basis of $gl(6,{\bf C})$. Hence, 
the triple of matrices $A_j$ is irreducible. 

In all other neighbouring cases with possibility 1) 
the left lower block of the matrix $Y$ has a single unit in its left upper 
corner and zeros elsewhere.

{\bf $3^0$. Possibility 2).}

Consider the case neighbouring to b1) (the ones neighbouring to a1), 
c2) and d3) are considered by analogy). Consider a block 
upper-triangular triple of matrices $A_j^1$ with a trivial centralizer 
like in case b1) with $g=2$; the diagonal blocks are $3\times 3$, they 
define non-equivalent representations. Consider its direct sum with an 
irreducible triple of matrices $A_j^2$ defined by example (ex2). There 
exists a semi-direct sum of such triples (one can use arguments like 
the ones from the proof of Lemma~\ref{Ext>=2}; 
the matrices are block upper-triangular, with diagonal blocks 
$3\times 3$). After this 
one deforms the triple into a nearby irreducible one like in the previous 
example -- one sets 

\[ A_1=\left( \begin{array}{ccccccccc}0&1&0&0&0&0&0&0&0\\0&0&1&0&0&0&0&0&0\\
0&0&0&1&0&0&0&0&0\\0&0&0&0&0&0&0&0&0\\0&0&0&0&0&1&0&0&0\\0&0&0&0&0&0&0&0&0\\
\varepsilon &0&0&0&0&0&0&1&0\\0&0&0&0&0&0&0&0&1\\0&0&0&0&0&0&0&0&0
\end{array}\right) \] 
and $A_j=(I+\varepsilon X_j(\varepsilon ))^{-1}A_j^*
(I+\varepsilon X_j(\varepsilon ))$ for $j=2,3$. The matrix algebra 
${\cal A}$ contains all block upper-triangular matrices from 
$gl(9,{\bf C})$ with blocks $3\times 3$, see \cite{Ko4}. 
Hence, ${\cal A}'$ contains a matrix of the form 
$F=E_{7,7}+O(\varepsilon )$. 

One has $(A_1)^3=E_{1,4}+\varepsilon E_{7,3}$. The matrix 
$R=F(A_1)^3/\varepsilon$ belongs 
to ${\cal A}'$. One has $R_{7,3}\neq 0$ for $\varepsilon =0$.  
Like in the previously considered case one concludes that 
${\cal A}'=gl(9,{\bf C})$ and that the triple $A_1$, $A_2$, $A_3$ is 
irreducible.

{\bf $4^0$. Possibility 3).}

In the case neighbouring to b1) (the ones neighbouring to a1), c2) and d3) 
are considered by analogy) one has $g\geq 3$. For 
$g=3$ one sets $A_j=A_j^*+\varepsilon A_j^0$ where the 
triple of matrices $A_j^*$ is block upper-triangular, with blocks 
$3\times 3$, the diagonal blocks defining non-equivalent irreducible 
representations. Set 

\[ A_1^*+\varepsilon A_1^0=
\left( \begin{array}{ccccccccc}0&1&0&0&0&0&0&0&0\\0&0&1&0&0&0&0&0&0\\
0&0&0&0&0&0&0&0&1\\0&0&0&0&1&0&0&0&0\\0&0&0&0&0&1&0&0&0\\0&0&0&0&0&0&0&0&1\\
\varepsilon &0&0&-\varepsilon &0&0&0&1&0\\0&0&0&0&0&0&0&0&1\\
0&0&0&0&0&0&0&0&0\end{array}\right) \]
For $j=2,3$ one sets 
$A_j^*=\left( \begin{array}{ccc}A_j^1&0&A_j^1D_j^1-D_j^1A_j^3\\
0&A_j^2&A_j^2D_j^2-D_j^2A_j^3\\0&0&A_j^3\end{array}\right)$ (i.e. $A_j^*$ are 
the matrices $A_j$ from the previous subsubsection). The blocks 
$A_j^i$, $i=1,2,3$ are nilpotent rank 2 matrices. Set 
$A_j=(I+\varepsilon X_j(\varepsilon ))^{-1}A_j^*
(I+\varepsilon X_j(\varepsilon ))$ for $j=2,3$. 
We let the reader check oneself that the matrices $A_j$ have the necessary 
Jordan normal forms. 

One shows next that the algebra ${\cal A}$ contains all matrices with 
arbitrary blocks in positions (1,1), (1,3), (2,2), (2,3), (3,3) (using 
\cite{Ko4}). Hence, the algebra ${\cal A}'$ contains a matrix 
$P=E_{7,7}+O(\varepsilon )$ and the matrix $P(A_1)^3/\varepsilon$ has (for 
$\varepsilon =0$) non-zero entries in positions (7,3) and (7,6). Like in the 
previous case one concludes that ${\cal A}'=gl(9,{\bf C})$, 
i.e. the triple $A_1$, $A_2$, $A_3$ is irreducible.

{\bf $5^0$. Possibility 4).} 

One must have $g\geq 4$. Let $g=4$. 
Consider the case neighbouring to b1).  
One constructs a triple of matrices $A_j^*$ with 

\[ A_1^*=\left(Ê\begin{array}{cccc}A_1^0&\varphi _1&\varphi _2&\varphi _3\\
0&0&0&0\\0&0&0&0\\0&0&0&0\end{array}\right) ~,~
A_2^*=\left(Ê\begin{array}{cccc}A_2^0&\eta _1&\eta _2&\eta _3\\
0&0&1&0\\0&0&0&1\\0&0&0&0\end{array}\right)Ê\] 
and $A_3^*=-A_1^*-A_2^*$, $A_j^0$ are $9\times 9$. The blocks of $J(A_j^0)$ 
are of sizes 
4,4,1; 3,3,3; 3,3,3, the triple of nilpotent matrices $A_j^0$ is 
irreducible. Its 
existence follows from the case neighbouring to b1) from possibility 3). 
One chooses the 
vector-columns $\varphi _j$ and $\eta _j$ such that $J(A_2^*)$ and 
$J(A_3^*)$ to have each four blocks of size 3 and $J(A_1^*)$ to have 
blocks of sizes 4,4,2,1,1.

Assume that $A_1^*$ is in Jordan normal form and that $\varphi _1$ has 
a unit in its last position and zeros elsewhere. 
Set $A_1=A_1^*+\varepsilon L$ where only the last row of $L$ is non-zero, 
with $L_{12,11}\neq 0$ and $L_{12,10}=L_{12,12}=0$. 
One chooses $L$ such that for $\varepsilon \neq 0$ the 
sizes of the blocks of $J(A_1)$ to be 4,4,2,2. 
Set for $j=2,3$ $A_j=(I+\varepsilon X_j(\varepsilon ))^{-1}A_j^*
(I+\varepsilon X_j(\varepsilon ))$. After this the irreducibility of the 
triple $A_1$, $A_2$, $A_3$ is proved by analogy with case (F).

In the case neighbouring to c2) one sets 

\[ A_1^*=\left(Ê\begin{array}{ccccc}
A_1^0&\varphi _1&\varphi _2&\varphi _3&\varphi _4\\
0&0&0&0&0\\0&0&0&1&0\\0&0&0&0&1\\0&0&0&0&0\end{array}\right) ~,~
A_2^*=\left(Ê\begin{array}{ccccc}A_2^0&\eta _1&\eta _2&\eta _3&\eta _4\\
0&0&-1&0&0\\0&0&0&-1&0\\0&0&0&0&-1\\0&0&0&0&0\end{array}\right)Ê\] 
where $A_j^0\in gl(12,{\bf C})$. The sizes of the blocks of $J(A_j^*)$ are 
5,5,3,3; 4,4,4,4; 2,2,2,2,2,2,2,1,1. The sizes of the blocks of 
$J(A_j^0)$ are 5,5,2; 4,4,4; 2,2,2,2,2,2. The existence of such an 
irreducible triple of nilpotent matrices $A_j^0$ follows from the 
case neighbouring to b1) from possibility 3). Set $A_3=A_3^*+\varepsilon L$, 
$A_j=(I+\varepsilon X_j(\varepsilon ))^{-1}A_j^*
(I+\varepsilon X_j(\varepsilon ))$, $j=1,2$. One has $L_{16,15}\neq 0$ 
and $L_{16,13}=L_{16,14}=L_{16,16}=0$. For the rest of the 
reasoning is like in the previous case. 

In the case neighbouring to d3) one sets 

\begin{equation}\label{cc} 
A_1^*=\left(Ê\begin{array}{cc}
A_1^0&T_1\\0&A_1'\end{array}\right) ~,~
A_2^*=\left(Ê\begin{array}{cccc}A_2^0&T_2\\0&A_2'\end{array}\right)Ê
\end{equation} 
where $A_j^0\in gl(18,{\bf C})$ and 

\[ A_1'=\left(Ê\begin{array}{cccccc}0&0&0&0&0&0\\0&0&1&0&0&0\\0&0&0&1&0&0\\
0&0&0&0&1&0\\0&0&0&0&0&1\\0&0&0&0&0&0\end{array}\right) ~,~
A_2'=\left(Ê\begin{array}{cccccc}0&-1&0&0&0&0\\0&0&-1&0&0&0\\0&0&0&0&0&0\\
0&0&0&0&-1&0\\0&0&0&0&0&-1\\0&0&0&0&0&0\end{array}\right) ~;\] 
$J(A_2^*)$ has blocks of size 3; $J(A_3^*)$ has eleven blocks of size 2 and 
two of size 1; $J(A_1^*)$ has blocks of sizes 7,7,5,5. The sizes 
of the blocks of $J(A_j^0)$ are respectively 7,7,4; six times 3; 
nine times 2. 

Set $A_3=A_3^*+\varepsilon L$, 
$A_j=(I+\varepsilon X_j(\varepsilon ))^{-1}A_j^*
(I+\varepsilon X_j(\varepsilon ))$, $j=1,2$. One 
has $L_{24,23}\neq 0$ and $L_{24,j_0}=0$ for $1\leq j_0\leq 24$, 
$j_0\neq 23$. 
Irreducible triples of such nilpotent matrices $A_j^0$ exist by the case 
neighbouring to d3) from possibility 3). The rest of the 
reasoning is like in the previous two cases.

In the case neighbouring to a1) one represents $A_j^*$ in the form (\ref{cc}) 
with  

\[ A_1'=A_4'=0~,~A_2'=-A_3'=\left(Ê\begin{array}{cc}0&1\\0&0\end{array}\right) \]
with $A_j^0\in gl(6,{\bf C})$. The sizes of the blocks of $J(A_j^*)$ are 
3,3,1,1; 2,2,2,2; 2,2,2,2; 
2,2,2,1,1. The ones of $J(A_j^0)$ are 3,3; 2,2,2; 2,2,2; 2,2,2. The 
existence of an irreducible quadruple of nilpotent matrices $A_j^0$ follows 
from the case neighbouring to a1) from possibility 3).  
Set $A_4=A_4^*+\varepsilon L$, 
$A_j=(I+\varepsilon X_j(\varepsilon ))^{-1}A_j^*
(I+\varepsilon X_j(\varepsilon ))$, $j=1,2,3$. The matrix $L$ 
has a single non-zero entry in position $(8,7)$. The rest of the 
reasoning is like in the previous three cases.

The theorem is proved.

\end{document}